\documentclass[12pt,twoside,reqno,openany]{amsart}

\pdfoutput=1

\usepackage{amsbsy,amscd,amsfonts,amsmath,amssymb,amsthm,color,
fancybox,fancyhdr,footmisc,graphics,graphicx,ifthen,mathrsfs,
multicol,pdfpages,rotating,times,wasysym}
\usepackage[dvipsnames,svgnames,x11names]{xcolor}

\usepackage[all]{xy}
\usepackage[french]{babel}
\usepackage[utf8]{inputenc}
\usepackage[T1]{fontenc}
\usepackage{mathtools}
\newtagform{EngelLie}[\scriptstyle]{$}{$}
\makeatletter\newcommand{\leqnomode}{\tagsleft@true}
\newcommand{\reqnomode}{\tagsleft@false}\makeatother

\newtheorem{Theoreme}[equation]{Th\'eor\`eme}

\newtheorem{These}[equation]{Th\`ese}
\newtheorem{Proposition}[equation]{Proposition}

\theoremstyle{definition}

\newtheorem{Question}[equation]{Question}

\newcommand{\C}{\mathbb{C}}

\newcommand{\R}{\mathbb{R}}

\newcommand{\Z}{\mathbb{Z}}

\newcommand{\GG}{\text{\sc g}}

\newcommand{\LL}{\text{\sc l}}
\newcommand{\MM}{\text{\sc m}}

\newcommand{\TT}{\text{\sc t}}

\newcommand{\maux}{{\text{\usefont{T1}{qcs}{m}{sl}m}}}

\definecolor{blue}{cmyk}{1.,1.,0.,0.63}
\definecolor{red}{cmyk}{0.,1.,1.,0.63}
\definecolor{green}{cmyk}{1.,0.,1.,0.63}
\definecolor{black}{cmyk}{1.,1.,1.,1.}

\newcommand{\blue}{\textcolor{blue}}
\newcommand{\green}{\textcolor{green}}
\newcommand{\red}{\textcolor{red}}

\makeatletter
\renewcommand{\@fnsymbol}[1]
{\ensuremath{\ifcase#1\or $*$\or $**$\or $***$\or $****$\or $*****$
\else\@ctrerr\fi}}
\makeatother

\newcommand{\HEAD}[2]{%
\pagestyle{fancy}
\fancyhead[RO]{\tiny\sf\thepage}
\fancyhead[CO]{{\tiny\sf #1}}
\fancyhead[LE]{\tiny\sf\thepage}
\fancyhead[CE]{{\tiny\sf #2}}
\fancyfoot{}}

\newcommand{\CITATION}[1]{\smallskip\hfill
\begin{minipage}[t]{13.5cm}\baselineskip=0.43cm\parindent=0.91cm
\blue{\footnotesize{\sf \!\!\!\!\!\!\!#1}}\end{minipage}\medskip}

\numberwithin{equation}{section}

\newcommand{\Section}[1]{
\renewcommand{\thesection}{\bf\arabic{section}}
\section{#1}
\renewcommand{\thesection}{\arabic{section}}}

\newcommand{\style}[1]{\text{\footnotesize{\sf #1}}}

\newcommand{\stylesmall}[1]{{\sf #1}}

\newcommand{\Angle}{\style{Angle}}

\renewcommand{\cos}{\style{cos}}

\renewcommand{\det}{\style{det}}

\renewcommand{\exp}{\style{exp}}

\renewcommand{\lim}{\style{lim}}

\newcommand{\longueur}{\style{longueur}}

\newcommand{\modsmall}{\stylesmall{mod}}

\newcommand{\proj}{\style{proj}}

\newcommand{\rad}{\style{rad}}

\newcommand{\Rotation}{\style{Rotation}}

\renewcommand{\sin}{\style{sin}}

\newcommand{\Bignorm}{\Big\vert\!\Big\vert}

\newcommand{\centersmallbullet}{{}_{{}^{{}^{
\scriptscriptstyle{\bullet\!}}}}}

\newcommand{\linestop}{\medskip\centerline{\bf 
-----------------}\medskip}

\newcommand{\smallbullet}{{\scriptscriptstyle{\bullet}}}

\newcommand{\vf}{\vfill\end{document}}

\begin{document}

\setcounter{section}{0}

\begin{center}

{\large\bf
Th\'eor\`eme de Gauss-Bonnet pour les surfaces:}
\label{Gauss-Bonnet-surfaces}

\medskip

{\large\bf
vers une philosophie intrinsèque de la géométrie différentielle}

\medskip

{\bf (Partie I)}

\bigskip

Jo\"el {\sc Merker}\footnote{Laboratoire de Mathématiques d'Orsay,
Université Paris-Sud, CNRS, Université Paris-Saclay, 91405 Orsay Cedex,
France. {\bf joel.merker@math.u-psud.fr}} 
et 
Jean-Jacques {\sc Szczeciniarz}\footnote{
Laboratoire SPHERE, CNRS, Université Paris Diderot (Paris 7),
75013 Paris, France. 
{\bf jean-jacques.szczeciniarz@univ-paris-diderot.fr}}

\smallskip

{\large\footnotesize\sf 
Universit\'e Paris-Sud et Universit\'e Paris Diderot}

\end{center}

\begin{center}
\begin{minipage}[t]{12.5cm}
\parindent 0.53cm
\scriptsize
\noindent
{\sc Résumé}.
\`A juste titre, Hegel et Schopenhauer ont critiqué sévèrement les
diverses démonstrations du théorème $\widehat{A} + \widehat{B} +
\widehat{C} = \pi$ d’Euclide, leur adressant le reproche d’introduire
de façon arbitraire des artifices de tracé qui conduisent le
philosophe spéculatif à un sentiment de malaise en présence de
<<\,tours d’escamotage\,>> au cours desquels la vérité <<\,s’introduit
par la petite porte dérobée\,>>.

En élaborant des analyses philosophiques motivées et {\sl 
conceptuellement
génétiques} qui seront internes à la contexture métaphysique même de la
pensée mathématique\,\,---\,\,trop souvent oblitérée par des pratiques
réductrices\,\,---, nous argumenterons au contraire 
que la vérité <<\,entre
par la grande porte\,>> dans les théories mathématiques, et ce, en nous
référant au développement ultime du théorème d’Euclide, l’un des plus
importants de toute la géométrie différentielle contemporaine, à
savoir le {\sl Théorème de Gauss-Bonnet}, emblème paradigmatique d’une
liaison, au sens lautmanien du terme, entre 
{\sl Topologie} et {\sl Géométrie}.

Effectivement, l'énoncé du {\sl Théorème de Gauss-Bonnet} fait surgir
une forme inattendue de {\sl réflexivité}\,\,---\,\,concept majeur et
universel de philosophie des mathématiques\,\,---, de sorte que la
géométrie se contemple elle-même, atteignant une double polarité
intuitive en elle-même. C’est donc le concept révolutionnaire et
protéiforme de {\sl courbure gaussienne} qui déclenche une
conceptualité nouvelle au-dessus de la géométrie euclidienne.
 
Ici, l’égalité entre intégrale de courbure totale
et caractéristique d'Euler indique qu’un concept de nature topologique
est {\em égal} à un nombre qui exprime un concept de nature
géométrique.  Ceci démontre encore que les mathématiques se
développent par intervention de disciplines différentes les unes sur
les autres, comme outils d’observation, structuration formelle,
nouveaux points de vue unificateurs.  Grâce aux trajectoires multiples
que permet l’architecture mobile des mathématiques, nous disposons
toujours
d'une démultiplication des significations des concepts mathématiques.

Toute une force d’embrassement des êtres rationnels est à l'{\oe}uvre
dans les mathématiques, qui expriment et réalisent un désir continué
de {\sl puissance synthétique}. La conceptualisation est certes une
première force d’embrassement, mais les êtres, toujours, s’y
soustraient, en partie, car le concept, en mathématiques, n’est
souvent qu’une vue partielle, d’autres sous-concepts demeurant
non-vus, d’autres aspects intrinsèques demeurant non-conceptualisés.

Seules les synthèses complètes telles que le {\sl Théorème
de Gauss-Bonnet} embrassent réellement les êtres mathématiques.

Cette force d’entraînement hors de soi que manifeste la pensée
mathématique dans sa puissance d’expansion est ici propulsée à
travers cette unité grâce aux multiples ressorts de la topologie et de
la géométrie.  Ainsi, la question des fondements de la géométrie
(euclidienne) se trouve-t-elle entièrement transformée dans les
nouveaux cadres théoriques de la géométrie différentielle, fondée par
le moyen de ses extensions synthétiques, qui en sont autant
d’accroissements de connaissances.

\end{minipage}
\end{center}

\Section{\bf \'Enonc\'e du th\'eor\`eme de Gauss-Bonnet classique}
\label{enonce-Gauss-Bonnet-classique}
\HEAD{\ref{enonce-Gauss-Bonnet-classique}.~{\sf 
\'Enonc\'e du th\'eor\`eme de Gauss-Bonnet classique}
}{
Jo\"el {\sc Merker} et Jean-Jacques {\sc Szczeciniarz}}

\CITATION{
\scriptsize
One of the star theorems of differential geometry is the
Gauss-Bonnet theorem, which for a compact oriented surface
$S$ states that:
\[
\int_S\,K_\GG\,dA
\,=\,
2\pi\,
\chi(S).
\eqno
\text{\green{Michael {\sc Spivak}}}
\]
}

Le th\'eor\`eme de Gauss-Bonnet propose une forme aboutie et d'une
certaine mani\`ere ultime d'un th\'eor\`eme classique d'Euclide.

\begin{center}
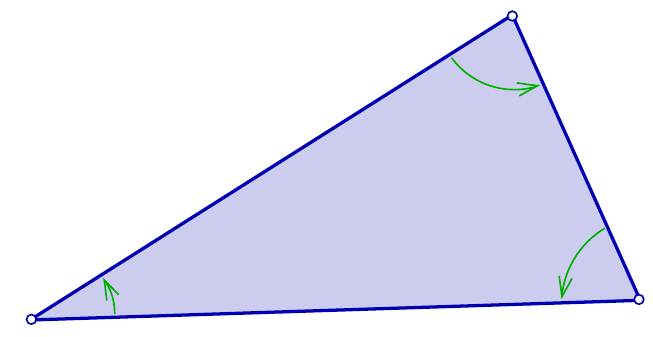
\end{center}

\'Etant donn\'e un triangle $T = ABC \subset \mathcal{P}$ dans le
plan euclidien $\mathcal{P} \cong \R^2$ ayant trois sommets distincts
$A$, $B$, $C$, ce th\'eor\`eme maintenant consid\'er\'e par les g\'eom\`etres
contemporains comme si \'el\'ementaire qu'il s'inscrit presque comme
imm\'ediatet\'e perceptuelle dans leur intuition constituante, \'enonce que
la somme des angles de $T$ en ses sommets vaut toujours:
\[
\pi
\,=\,
\widehat{A}
+
\widehat{B}
+
\widehat{C},
\]
la mesure des angles s'effectuant en {\sl radians}, avec la 
normalisation standard $\pi\, \rad = 180^{\rm o}$.

Cette d\'ecouverte ancienne constitue
une des racines fondamentales de la pens\'ee math\'ematique, en tant que
cette pens\'ee est un rapport au monde qui cherche \`a capturer
visuellement et mentalement des formes, des invariances,
des structures, quelle que soit l'\'etendue et la complexit\'e
du champ ontologique sous-jacent.

Ici, au-del\`a de l'\'enonciation universelle figeante, et au-del\`a 
de la quantification universelle r\'esumable par l'emploi
de quelque symbolisme logique qui soit:
\[
\forall\,
A\,\in\,\mathcal{P},
\ \ \ \ \
\forall\,
B\,\in\,\mathcal{P},
\ \ \ \ \
\forall\,
C\,\in\,\mathcal{P},
\]
il faut s'imaginer une variabilit\'e mobile de tous les triangles
possibles $ABC \subset \mathcal{P}$ dans le plan, faire voir \`a son
intuition math\'ematique interne trois points libres bouger en m\^eme
temps que les trois c\^ot\'es $AB$, $BC$, $CA$ du triangle $T$ se
d\'eplacent comme trois membres articul\'es, et se repr\'esenter ainsi
par l'{\sl entendement d'appropriation de contenu g\'eom\'etrique}
qu'une constance quantitative s'observe dans une variabilit\'e {\em a
priori} inorganis\'ee. Cette variation se fait dans (le) un plan. Si
bien qu'elle peut \^etre vue comme une exploration du plan, de ce qui
le structure comme \'etendue `plate'. Et une compr\'ehension
intellectuelle, appropriation d'un contenu qui est dans cette
appropriation m\^eme g\'eom\'etrique et qui en l'occurrence met en
\'evidence cette constance. On doit pouvoir dire que la base d'une
appropriation de contenu g\'eom\'etrique est une une constance
(quantitative) dans une variabilit\'e inorganis\'ee.

\begin{Question}
La mentalisation math\'ematique est-elle transmissible
par le texte? par le langage? par la diagrammatisation?
\end{Question}
Cette int\'egration par la pens\'ee qui se traduit dans une expression
de la conscience se construit dans la nature m\^eme de la constitution
qui met en \'evidence ces invariants premiers. Il s'agit de la 
fa\c{c}on 
dont l'\'etendue (concept travaill\'e par tous les philosophes)
vient \`a l'espace, comme quantit\'e \'etendue. Nous reprendrons cette
question plus bas en situant le travail de Kant.

Il faut rappeler que le th\'eor\`eme d'Euclide informe toute
l'histoire de la philosophie. Hegel en critique la d\'emonstration et
lui adresse le reproche d'introduire de fa\c{c}on arbitraire un
artifice ext\'erieur, \`a savoir le trac\'e d'une parall\`ele \`a la
base du triangle, laquelle est ensuite utilis\'ee pour faire jouer le
th\'eor\`eme de Thal\`es au moyen de paires d'angles alternes-internes
{\em \'egaux}. 
consid\'erons une critique repr\'esentative de l'histoire de la philosophie, celle de Schopenhauer.

\CITATION{
Nous sommes certainement forc\'es de reconna\^itre, en vertu du
principe de contradiction que ce 
qu'Euclide d\'emontre est bien tel qu'il
le d\'emontre; mais nous n'apprenons pas pourquoi il en est
ainsi. Aussi \'eprouve-t-on presque le m\^eme sentiment de malaise
qu'on \'eprouve apr\`es avoir assist\'e \`a des tours d'escamotage,
auxquels en effet, la plupart des d\'emonstrations d'Euclide
ressemblent \'etonnamment. Presque toujours, chez lui, la v\'erit\'e
s'introduit par la petite porte d\'erob\'ee, car elle r\'esulte,
par accident, de quelque circonstance accessoire; dans certains cas la
preuve par l'absurde ferme successivement toutes les portes, et n'en
laisse ouverte qu'une seule, par laquelle nous sommes contraints de
passer par ce seul motif. Dans d'autres, comme dans le th\'eor\`eme
de Pythagore, on tire des lignes, on ne sait pour quelle raison; on
s'aper\c{c}oit pus tard que c'\'etaient des n{\oe}uds coulants qui
serrent \`a l'improviste, pour surprendre le comportement du curieux
qui cherchait \`a s'instruire; celui-ci tout saisi, est oblig\'e
d'admettre une chose dont la contexture intime lui est encore
parfaitement incomprise, et cela \`a tel point qu'il qu'il pourrra
\'etudier Euclide en entier sans avoir une compr\'ehension effective
des relations de l'espace.\footnotemark}

\footnotetext{Schopenhauer,
{\em Le monde comme volont\'e et comme repr\'esentation}, Livre I, 
15.}

D'ailleurs, il existe de nombreuses autres d\'emonstrations de ce
th\'eor\`eme de g\'eom\'etrie du triangle, comme celle qui consiste
\`a glisser le long des c\^ot\'es en effectuant une rotation \`a
chaque angle, mais \`a toutes ces d\'emonstrations, avec Hegel et
Schopenauer, on pourrait aussi adresser le m\^eme reproche
d'articifialit\'e impr\'evisible des constructions. Certainement, il
faudrait examiner plus en profondeur cette question de
l'<<\,arbitraire\,>> de la d\'emonstration, lequel, chez Hegel, est
simplement une preuve (ou une confirmation) du caract\`ere arbitraire
et dogmatique des math\'ematiques\,\,---\,\,\`a la diff\'erence,
selon lui, de la philosophie. Mais cela nous engagerait dans des
interrogations sans limites sur l'instantan\'eit\'e myst\'erieuse du
caract\`ere {\em a posteriori} de tout acte math\'ematique guid\'e
par l'intelligence des contenus. 

Pour l'instant, il nous faut dire que
l'arbitraire est ce dont tout mth\'ematicien essaie de se d\'efaire
en réalisant des d\'emonstrations les plus intrins\`eques possible, et de
l\`a, nous allons tenter d'apporter une r\'eponses aux objections
philosophiques formul\'ees plus haut qui critiquent la nature des math\'ematiques. 
De ce point de vue il peut \^etre dit que la
g\'eom\'etrie que nous allons examiner <<\,tend vers la philosophie\,>> et
pr\'esente dans son mouvement m\^eme une expression philosophique et
un appel \`a la philosophie. C'est pourquoi il est souvent ais\'e
d'attribuer aux math\'ematiques elles-m\^emes les propri\'et\'es que
Hegel rep\`ere dans la philosophie elle-m\^eme.

En se pla\c{c}ant d'un point de vue sup\'erieur qui embrasse une
multiplicit\'e de g\'eom\'etries possibles, il s'av\`ere que le
th\'eor\`eme d'Euclide est d'abord et avant tout une expression de la
nature <<\,euclidienne\,>> du plan. Et que
l'introduction\,\,---\,\,consid\'er\'ee comme arbitraire ou
ext\'erieure\,\,---\,\,d'une parall\`ele \`a un c\^ot\'e du triangle
correspond \`a une autre expression de cette nature euclidienne. Comme
nous allons le voir, les m\'etriques riemanniennes, les connexions
sym\'etriques, et le th\'eor\`eme de Gauss-Bonnet permettent en effet de
mettre en \'evidence des modulations th\'ematiques \'eclairantes sur la
diversit\'e ontologique des g\'eom\'etries possibles.

Par cons\'equent, en se pla\c{c}ant donc d'un point de vue math\'ematique
<<\,{\sl relev\'e}\,>>\,\,---\,\,dans un sens
proche de l'{\sl Aufhebung} h\'egelienne\,\,---\,\,qui 
tient et doit tenir compte
d'ascensions ontologiques existantes, puisque le th\'eor\`eme
d'Euclide peut et doit \^etre vu comme l'expression d'une
propri\'et\'e essentielle de la nature euclidenne du plan, et m\^eme
de l'espace tridimensionnel, et puisque la d\'emonstration qui fait
appel \`a la notion de parall\`ele et \`a des propri\'et\'es purement
euclidiennes que l'on d\'esigne sous le nom de {\sl Th\'eor\`eme de
Thal\`es}, d\'eveloppe une compr\'ehension de cette propri\'et\'e
essentielle, alors le reproche h\'eg\'elien cesse d'avoir lieu
d'\^etre. Car effectivement, il faut
pouvoir se dire et se repr\'esenter non seulement 
qu'une parall\`ele \`a la base passant par le sommet du triangle
existe et qu'elle est unique, mais surtout que 
c'est l\`a une propri\'et\'e caract\'eristique de
la nature euclidienne d'un espace. 

Une question surgit alors: quel statut accorder au fameux axiome des
parall\`eles dans ce cadre de r\'eflexion? Si m\^eme les
propri\'et\'es euclidiennes famili\`eres sont devenues quasiment des
propri\'et\'es de notre perception de l'espace et du plan, comment en
sommes-nous arriv\'es l\`a? L'axiome des parall\`eles d\'efinit
l'espace euclidien, et il nous fournit aussi les modalit\'es de notre
d\'eplacement dans cet espace, car nous pouvons effectuer des
translations, des rotations, des r\'eflexions.

En tout cas, il faut selon nous distinguer entre une r\'eflexion
math\'ematique en profondeur sur l'espace et une pr\'esentation
axiomatique de la g\'eom\'etrie de l'espace comme nous la trouvons
dans l'{\oe}uvre d'Euclide, qui s'est transmise pendant plus de deux
mill\'enaires. Toute pr\'esentation axiomatique positionne les
r\'esultats qui en sont issus comme effet d'une n\'ecessit\'e
d\'eductive. Cette n\'ecessit\'e se pr\'esente comme interne \`a la
formation de concepts, mais souvent sans \'eclairages g\'en\'etiques
explicites. La r\'eflexion conceptuelle est l'effet de la construction
du contenu math\'ematique m\^eme. Elle est de l'ordre de la synth\`ese
constructive. Et la synth\`ese constructive suppose elle-m\^eme un
examen conceptuel de l'espace sous diff\'erentes caract\'eristiques.

Un tel examen contribue \`a accro\^itre une connaissance de l'espace
qui va se trouver li\'ee \`a des \'etudes produites. Nous montrerons
que dans la forme de ce devenir historique appara\^it comme une
n\'ecessit\'e de la r\'eflexion de la g\'eom\'etrie sur elle-m\^eme,
et nous soutiendrons que les formes de sa transmission m\^eme sous
forme axiomatique en sont encore une manifestation. Et cette
r\'eflexion interne des math\'ematiques sur elles-m\^emes peut adopter
des formes extr\^emement diverses.

Or dans le destin de la g\'eom\'etrie, deux aspects essentiels
nouveaux sont apparus: les propri\'et\'es m\'etriques, puis les
propri\'et\'es topologiques, d'une nature conceptuelle plus
d\'elicate. \`A ce point de notre analyse, nous sommes rendus proche
de la position de Hilbert devant l'axiomatique euclidienne\footnote
{\ Hilbert, {\em Grundlagen der Geometrie,} Mit Supplementen von Dr
Paul Bernays. Stuttgart, B. G. Teubner, 1968.} Certains \'enonc\'es
supposent des pr\'ecisions de nature topologique, d'autres de nature
m\'etrique, et donc, il faut tenir compte de propri\'et\'es
sup\'erieures de l'espace qui ont \'et\'e con\c{c}ues au cours de
l'histoire de la g\'eom\'etrie. La liaison entre g\'eom\'etrie et
topologie, longtemps rest\'ee implicite, s'est révélée
corr\'elativement avec l'\'emergence autonome de la topologie. Ce qui
signifie que toute propri\'et\'e g\'eom\'etrique est (peut-\^etre) en
liaison avec une propri\'et\'e topologique. Ce qui n'emp\^eche pas de
faire une distinction entre ces deux types de propri\'et\'es. Cela
change la question fondationnelle \`a laquelle les philosophes des
math\'ematiques dominants sont tant attach\'es. Tout \'enonc\'e
g\'eom\'etrique doit rester ouvert sur ses corr\'elations
\'eventuelles avec la topologie ou m\^eme avec l'alg\`ebre et
l'analyse qui en pr\'ecisent la signification. C'est l\`a tout le
travail de Hilbert, qui a cependant continu\'e de raisonner en termes
axiomatiques.
 
Notre th\`ese est que dans cette r\'eflexion sur ces liens se
d\'eploie une r\'eflexion philosophique en elle-m\^eme et pour
elle-m\^eme, au contact de la g\'eom\'etrie quand elle est relev\'ee \`a un
niveau math\'ematique sup\'erieur. Et comme nous allons
le voir, ce credo donne sens \`a une \'etude sp\'eculative
autonome sur la forme
g\'en\'erale
du th\'eor\`eme d'Euclide bien connue des math\'ematiciens: 
le {\em th\'eor\`eme de Gauss-Bonnet}.
 
\`A pr\'esent, revenons \`a notre repr\'esentation, comme nous l'avons
dit plus haut, d'une variabilit\'e mentalis\'ee, au sens
d'int\'egr\'ee \`a la pens\'ee de tous les triangles possibles d'un
plan. En refusant de souscrire aux tentatives heideg\'eriennes de
rel\'egation de l'activit\'e scientifique \`a une forme de confinement
au seul d\'eveloppement du principe de raison, il s'agissait
d'insister ici, en guise de pr\'eliminaire m\'etaphysique, sur
l'intrication fondamentale entre le mouvement et la quantification
universelle du multiple, ce \`a quoi seule une activit\'e mentale
libre et impr\'evisible peut et doit donner acc\`es. L'intuition {\em
 est} mouvement, et donc, puisque la math\'ematique {\em est}
synth\`ese d'intuitions, l'activit\'e c\'er\'ebrale intuitive demeure
absolument d\'eterminante dans l'appr\'ehension des contenus
math\'ematiques, fussent-ils axiomatiques. C'est-à-dire que
l'intuition comme saisie imm\'ediate, comme vision intellectuelle, est
li\'ee \`a l'appr\'ehension du mouvement.

Examinons quelque peu ce rapport \`a l'intuition, d'abord en suivant
Kant. Une d\'efinition de 1770 est pr\'ecis\'ee dans la {\em Critique
de la raison pure}.

\CITATION{
Celle-ci [la math\'ematique] ne peut rien \'etablir par simples
concepts mais vole aussit\^ot vers l'intuition dans laquelle elle
consid\`ere le concept {\em in concreto}, non pas empiriquement
pourtant, mais dans une intuition qu'elle repr\'esente {\em a priori}
c'est-\`a-dire a construite et dans laquelle ce qui suit des
conditions g\'en\'erales de la construction doit \'egalement valoir
pour l'objet construit.\footnotemark}

\footnotetext{Kant, {\em Kritik der reinen Vernunft},
Methodenlehre, I, 1, Cassirer, III, p. 486.}

\noindent
En reprenant l'expression de Cavaill\`es,

\CITATION{
[l'intuition] est non pas contemplation d'un tout fait, mais
appr\'ehension dans l'\'epreuve de l'acte des conditions m\^emes qui
la rendent possible.\footnotemark}

\footnotetext{\ 
{\em M\'ethode axiomatique et formalisme. 
Essai sur le probl\`eme du fondement des math\'ematiques},
Hermann, Paris, 1981.}

Et ensuite Kant, dans le m\^eme mouvement de
pens\'ee, distingue philosophie et math\'ematiques de mani\`ere tout
\`a fait diff\'erente de la conception de Hegel.

\CITATION{
Le philosophe peut r\'efl\'echir sur le concept de triangle autant
qu'il veut [\dots],
y distinguer le concept de ligne droite ou d'un angle
ou du nombre trois et les \'eclaircir, sans pour cela parvenir \`a
d'autres propri\'et\'es qui ne sont pas renferm\'ees dans ces
concepts.\footnotemark}

\footnotetext{\ Kant, {\em Kritik der reinen Vernunft},
Methodenlehre, I, 1, Cassirer, III, p. 487.}

Si le math\'ematicien doit \^etre toujours guid\'e par l'intuition
dans l'encha\^inement de ses raisonnements, c'est que l'intuition
poss\`ede une structure, ou une r\'ealit\'e propre, de quelque ordre
soit-elle. Bien plus, la dualit\'e des deux formes de l'intuition
rend le probl\`eme insoluble.

Sans doute y a- t -il subordination de l'espace au temps, puisque toute synth\`ese s'accomplit dans le temps; 

\CITATION{
Je ne peux me repr\'esenter une ligne, 
si petite soit-elle, sans la tirer par la pens\'ee.\footnotemark}

\footnotetext{\ Axiomes de l'intuition, {\em op. cit.},
p. 27.}

\noindent
Et Kant explique l'application du nombre \`a l'espace.

\CITATION{
Le pur sch\'ema de la grandeur comme d'un concept de l'entendement est
le nombre qui est une repr\'esentation embrassant l'addition
successive de l'unit\'e \`a l'unit\'e (homog\`ene). Donc le nombre
n'est rien d'autre que l'unit\'e de la synth\`ese du divers 
{\em d'une intuition homog\`ene en g\'en\'eral} 
par le fait que j'engendre le
temps dans l'appr\'ehension de l'intuition.\footnotemark}

\footnotetext{\ 
{\em Sch\'ematisme}, {\em op. cit.}, p. 144.} 

Rappelons ici, au plus profond de la philosophie de Kant, que la
th\'eorie du sch\'ematisme est le moyen (conceptuel) qu'a invent\'e
Kant pour faire le lien entre d'une part l'intuition qui fournit des
perceptions singuli\`eres, et d'autre part les concepts de
l'entendement qui eux sont universels et actifs. Les sch\'emas sont
alors un mixte intuition\big/concept produit par l'imagination.

C'est ici que les difficult\'es de l'analyse kantienne
apparaissent. Il faut passer de cette <<\,intuition homog\`ene en
g\'en\'eral\,>> 
aux deux intuitions particuli\`eres de notre facult\'e de
conna\^itre.

\CITATION{
[Le temps] contient les rapports de la succession, 
du simultan\'e et de ce qui est simultan\'e avec 
le successif (le persistant) 
[\dots], mani\`ere dont l'esprit s'affecte lui-m\^eme,\footnotemark}

\footnotetext{\ 
Esth\'etique transcendantale, {\em op. cit.}, p. 75} 

\noindent
tandis que l'espace est

\CITATION{
la propri\'et\'e formelle (de l'esprit) 
d'\^etre affect\'e par les objets.\footnotemark}

\footnotetext{\ {\em ibid.}, p.~76.}

La difficult\'e de cette double intuition est soulign\'ee par Kant.

\CITATION{
Ici, toute la difficult\'e r\'eside en ceci: comment le sujet peut-il
avoir une intuition int\'erieure de lui-m\^eme? Mais cette
difficult\'e est commune \`a toutes les th\'eories.\footnotemark}

\footnotetext{\ {\em op. cit.}, p. 127.}

Nous sommes tellement un objet par nous-m\^emes que nous avons besoin
de l'espace pour repr\'esenter notre vie int\'erieure, ajoute
Cavaill\`es. 
 
\CITATION{
Nous ne pouvons nous repr\'esenter le temps, qui n'est pourtant pas un
objet de l'intuition ext\'erieure, que sous l'image d'une ligne, en
tant que nous la tra\c{c}ons, mode de repr\'esentation sans lequel nous
ne pourrions pas reconna\^itre qu'il n'a qu'une dimension; de m\^eme,
la d\'etermination de la longueur du temps ou des moments pour toutes
nos perceptions int\'erieures doit \^etre tir\'ee de ce que les choses
ext\'erieures nous pr\'esentent de changeant, par suite nous devons
ordonner les d\'eterminations du sens intime [interne] dans le
temps comme nous ordonnons celle du sens externe dans
l'espace.\footnotemark}

\footnotetext{\ {\em op. cit.}, p. 129.}

L'analyse de Cavaill\`es est celle que nous 
reprendrons \`a notre compte 
pour lui faire emprunter une autre direction.

\CITATION{
Mais alors le temps s'\'evanouit dans l'espace, 
le temps est grandeur extensive\footnotemark}

\footnotetext{\ Cavaill\`es, {\em op. cit.}, p. 29.} 

\CITATION{
dans laquelle la repr\'esentation des parties rend possible la repr\'esentation du tout\footnotemark}

\footnotetext{\ {\em Axiomes de l'intuition} p. 157.}

\CITATION{
ou plut\^ot il n'y a qu'une seule intuition, celle du mouvement qui
trace les lignes, d\'ecrit les cercles, tire la ligne droite, image du
temps.\footnotemark}

\footnotetext{\ {\em Ibidem.}}

C'est l\`a que se situe la difficult\'e selon Cavaill\`es. Nous
effectuons une distinction entre nous comme sujet connaissant et un
syst\`eme d'objets qui affectent notre sensibilit\'e 
et d\`es lors, dit
encore Cavaill\`es, nous pouvons isoler ce qui se rapporte uniquement
au moi. La repr\'esentation du temps n'est obtenue que par abstraction
de l'espace dans la synth\`ese sensible, dit Cavaill\`es. Mais il a
fallu se donner pr\'ealablement l'espace avant la synth\`ese.

\CITATION{
Ce qui produit d'abord le concept de la succession, c'est le mouvement
comme acte du sujet, si {\em nous faisons abstraction de celui-ci} et
faisons attention seulement \`a l'action par laquelle nous
d\'eterminons le sens intime conform\'ement \`a sa
forme.\footnotemark}

\footnotetext{\ D\'eduction transcendantale, 
{\em op. cit.}, p. 128.}

L'espace est <<\,l'image pure de toutes les
quantit\'es\,>>\footnote{\ Sch\'ematisme, {\em op. cit.}, p. 144.}, et
il est de fait li\'e de mani\`ere interne \`a la quantit\'e. Et la
g\'eom\'etrie travaille et expose ces formes quantitatives. Et c'est
le mouvement qui trace ces formes quantitatives.

En effet, le mouvement est li\'e \`a la quantification universelle du
multiple\,\,---\,\,mais sous quelle forme? On sait depuis les
\'El\'eates que le mouvement est un parcours du multiple. Mais comment
intervient la quantification? Tout point, sans pr\'ecision de cette
notion, est travers\'e ou support\'e par le mouvement. le mouvement
s'\'eprouve comme effectuation du trac\'e de n'importe quelle
figure. Et comme l'espace est l'image pure de toutes les quantit\'es,
pour le sens ext\'erieur, le mouvement est condition de l'effectuation
de la spatialit\'e quantitative. 

Les difficult\'es
li\'ees \`a cette pr\'esentation sont les suivantes. 
La question de fond reste la question du 
temps qui a gliss\'e du temps-sch\`eme au temps-image
spatialis\'ee. Nous pouvons travailler sur ce temps quand par exemple
nous l'envisageons comme param\`etre, mais le temps que Kant
consid\`ere est celui qui n'est pas image, il est

\CITATION{
m\'ethode pour procurer \`a un concept son image.\footnotemark}

\footnotetext{\ {\em Sch\'ematisme}, {\em op. cit.}, p. 144.}
 
En second lieu mais la difficult\'e est moindre, {\em La
quantification universelle n'est pas li\'ee \`a une
empiricit\'e}. Tout comme la notion d'universalit\'e. Le mouvement
supporte l'universalit\'e sans pour autant la produire sous sa forme
id\'eelle. Si le mouvement est l'effectuation du quantitatif
g\'eom\'etrique, il n'est pas empirique mais condition de toute
empiricit\'e repr\'esent\'ee. C'est dans ce cadre que nous proposons
une r\'eflexion sur la g\'eom\'etrie, dont l'axe principal devient le
mouvement. On sait que l'une des difficult\'es rencontr\'ees par la
g\'eom\'etrie a \'et\'e de donner un statut au mouvement
g\'eom\'etrique, vu que ce dernier est toujours d'abord pens\'e comme
un mouvement physique, ce qui est le cas chez Aristote. Le mouvement
cette fois explicit\'e dans des concepts g\'eom\'etriques fait partie
de la r\'eflexion g\'eom\'etrique et de la mise en forme de la
g\'eom\'etrie euclidienne.

Dans notre repr\'esentation d'une variabilit\'e de tous les triangles
possibles d'un plan, il nous faut discerner trois points qui se
meuvent, consid\'erer leurs positions respectives les uns par rapport
aux autres, les joindre virtuellement par des segments de droite. Une
figure de triangle appara\^it d\`es lors que ces trois points ne sont
pas align\'es. La discernabilit\'e du triangle est universelle. C'est
l\`a un travail propre du spatial quantitatif. Mais pourquoi?

C'est une perception de nature topologique qui est premi\`ere: je
sais former une figure en joignant trois points non align\'es. Les
conditions de la connaissance sont-elles ainsi qu'elles nous 
`donnent'
un triangle? Puis-je ne pas voir un triangle quand je consid\`ere
trois points non align\'es qui sont joints? En quel sens? Si je
d\'efinis le triangle comme une figure form\'ee par trois points non
align\'es que j'ai reli\'es, \'evidemment non. L'id\'ealit\'e triangle
est n\'ecessairement structur\'ee ainsi: qui se donne trois points 
non align\'es
dans le plan se donne un triangle. Ce qui veut dire
qu'il se donne trois angles solidaires. Cependant, l'objet triangle
concentre plus que le descriptif que j'en ai donn\'e. On conna\^it une
quantit\'e imposante de propri\'et\'es du triangle. En d\'ecrivant
l'objet triangle, je me donne ainsi un ensemble important de
propri\'et\'es possibles. C'est l'ensemble de ces propri\'et\'es\,\,---\,\,nous dirions 
virtualit\'es\,\,---\,\,que porte l'id\'ealit\'e triangle.

Nous n'avons pas r\'epondu \`a la question de sa discernabilit\'e
universelle. Tout le monde peut tracer un triangle en pens\'ee en
joignant trois points non align\'es, ou tout le monde peut parcourir
en pens\'ee un triangle. Dans ces deux, cas nous avons affaire \`a un
mouvement en pens\'ee. C'est ce mouvement qui est universel, imagin\'e
par tout un chacun.

Dans l'\'elaboration philosophique dont nous avons donn\'e quelques
\'el\'ements, Kant a tent\'e de r\'epondre \`a cette question par sa
th\'eorie du sch\'ematisme. Nous allons pr\'eciser ce qui est
n\'ecessaire pour notre propos. Kant suppose que la constitution de
l'objet et les conditions de l'acc\`es cognitif \`a l'objet sont
identiques. Il en est de m\^eme pour le sch\`eme qui a permis la
production du triangle. Caveing\footnote{\ {\em Le
probl\`eme des objets dans la pens\'ee math\'ematique}, Vrin, Paris,
2004.} pr\'ef\`ere user du terme d'{\sl omnisubjectivit\'e}.
Tout sujet
ayant acc\`es \`a la figure triangle passe par les m\^emes conditions
d'accession qui correspondent \`a la structure de l'objet,
pr\'ecis\'ement parce que ces conditions sont celles qui constituent
l'objet. Il faut voir dans l'omnisubjectivit\'e plus qu'une
g\'en\'eralisation \`a tout sujet. Nous dirons que c'est constitutif
de tout sujet. Cette supposition suffit-elle \`a nous faire sortir
d'une simple g\'en\'eralit\'e empirique? Et à passer \`a une
quantification universelle?

Il faudrait consentir \`a fonder empiriquement l'universalit\'e du
quantificateur. Tout triangle est ainsi construit en lui-m\^eme de
sorte que je puisse y avoir acc\`es. Tout triangle renvoie \`a tout
sujet percevant un triangle. De sorte que cette empiricit\'e
fondatrice (un sujet percevant) est de port\'ee
universelle. L'id\'ealit\'e triangle n'est aucun des triangles
trac\'es sur une feuille de papier. Elle concentre en elle la double
polarit\'e: celle de l'objet, et celle du sujet, dont la corr\'elation
d\'efinit l'universalit\'e. L'empiricit\'e l'a port\'ee au-del\`a
d'elle-m\^eme.

En adoptant une position proche de celle de Kant, nous ne
pr\'esupposons rien de la structure de cette omnisubjectivit\'e. Nous
ob\'eissons \`a une injonction th\'eorique impliqu\'ee par
l'universalit\'e que comporte un objet g\'eom\'etrique.

Le mouvement, celui qui engendre un triangle, ou celui qui trace un
triangle poss\`ede une port\'ee universelle.

Et c'est par l'appel au mouvement que nous pouvons comprendre la
nature d'un objet g\'eom\'etrique comme le triangle, nous ne devons
jamais en disjoindre sa relation \`a son engendrement. C'est l\`a ce
qu'il faut entendre par port\'ee universelle. En restant port\'es par
ce mouvement d'engendrement du triangle nous allons en faire sortir
une r\'eflexion sur soi de nature g\'eom\'etrique.

Donc si comme dit Kant, le philosophe ne sait pas quoi faire avec le
triangle, le g\'eom\`etre en produisant une r\'eflexion à partir de
sa construction va monter vers la philosophie, et transformer \`a
rebours la r\'eflexion philosophique. En approfondissant la
g\'eom\'etrie d'Euclide il la r\'efl\'echit et lui donne une
signification qui transforme et \'etend ses concepts. Le travail du
g\'eom\`etre permet d'extraire certains concepts qui donnent une
signification beaucoup plus ample au th\'eor\`eme d'Euclide. Ce
dernier est mis en position r\'eflexive gr\^ace \`a l'intervention du
concept de courbure. C'est \`a travers ce dernier que le mouvement
accomplit sa fonction constituante, et concentre de nombreux
\'el\'ements qui permettent une extension de la g\'eom\'etrie.

De ce point de vue ils permettent aussi de r\'epondre \`a certaines
questions qui faisaient difficult\'e chez Kant. En simplifiant, la
question de la dualit\'e des intuitions que nous avons rencontr\'ee
plus haut peut prendre la forme suivante. D\`es lors que nous pouvons
\'etablir qu'une forme de r\'eflexivit\'e surgit dans le
d\'eveloppement de la g\'eom\'etrie, que d'une certaine fa\c{c}on la
g\'eom\'etrie peut se voir elle-m\^eme, elle atteint cette double
polarit\'e en elle-m\^eme. D'autre part, nous allons le voir, le
th\'eor\`eme de Gauss-Bonnet nous permet d'op\'erer des abstractions
d'espace. Nous restons malgr\'e ces abstractions dans
l'ext\'eriorit\'e. C'est le concept de courbure qui nous place
au-dessus de la g\'eom\'etrie euclidienne.

De l'autre c\^ot\'e, la difficult\'e kantienne tient \`a
l'indiff\'erence \`a l'objet.

\CITATION{
Le type de la science, son effort vers un type rationnel classique
aussi bien que extension math\'ematique n\'egligent compl\`etement
l'apport de l'objet pour la structure de la th\'eorie.\footnotemark}

\footnotetext{\ Cavaill\`es,
{\em Sur la logique et la th\'eorie de la science},
Paris, Vrin, 1976.}

Le fondement par l'unit\'e formelle de l'espace tel que le propose
Kant, permettant de sortir du concept n'est pas possible. Cette
unit\'e varie quant \`a sa forme et son caract\`ere formel se
transforme dans la construction m\^eme de la g\'eom\'etrie. Se met en
place de la sorte une variation de l'intuition elle-m\^eme qui
emprunte des formes argumentatives.

Nous allons argumenter dans le sens ci-dessus pour donner une
signification philosophique\,\,---\,\,essentiellement issue de
l'approfondissement du travail g\'eom\'etrique sur l'espace 
euclidien\,\,---\,\,au d\'eveloppement ultime 
du th\'eor\`eme d'Euclide, à savoir l'un des plus
importants de la g\'eom\'etrie diff\'erentielle: le 
{\sl th\'eor\`eme de Gauss-Bonnet}.

\begin{center}
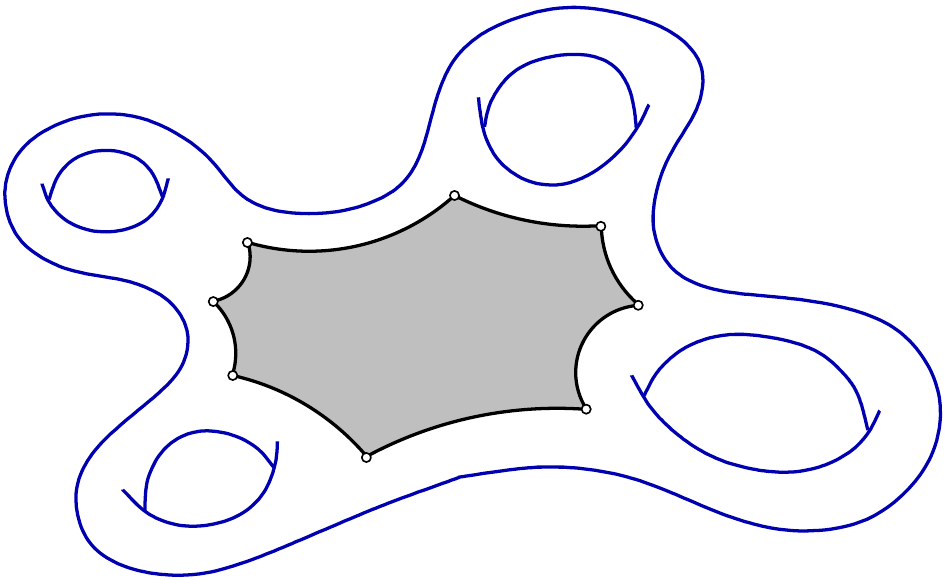
\end{center}

Directement maintenant, {\em i.e.} sans pr\'esentation conceptuelle 
anticipante\,\,---\,\,mais des \'eclaircissements commentatifs 
suivront\,\,---, et puisqu'il s'articule
au th\'eor\`eme d'Euclide, voici le

\begin{Theoreme}
\label{theoreme-Gauss-Bonnet-general}
{\bf [de Gauss-Bonnet g\'en\'eral]}
Sur une surface orient\'ee lisse $S$ munie
d'une m\'etrique gaussienne $ds^2$, soit un domaine relativement 
compact non vide:
\[
D
\,\Subset\,
S,
\]
dont le bord lisse par morceaux:
\[
\partial D
\,=\,
C_1
\cup\cdots\cup
C_\LL
\]
est constitu\'e d'un nombre fini $\LL \geqslant 1$ de courbes
ferm\'ees simples orient\'ees $C_1, \dots, C_\LL$ lisses par morceaux.
Alors sa caract\'eristique d'Euler $\chi(D)$ est donn\'ee par
la formule:
\[
2\pi\,\chi(D)
\,=\,
\int\!\!\int_D\,
K_\GG\,dA
+
\sum_{1\leqslant\ell\leqslant\LL}\,
\int_{C_\ell}
k_g\,ds
+
\sum_{1\leqslant\ell\leqslant\LL}\,
\sum_{1\leqslant\maux\leqslant\MM_\ell}\,
\alpha_{\ell,\maux},
\] 
o\`u $K_\GG$ est la courbure de Gauss 
de $S$, o\`u $dA$ est l'\'el\'ement d'aire
infinit\'esimale sur $S$, o\`u $k_g$ est la courbure g\'eod\'esique, et o\`u
$\alpha_{\ell,1}, \dots, \alpha_{\ell,\MM_\ell}$ 
sont les angles ext\'erieurs entre deux
morceaux lisses cons\'ecutifs de $C_\ell$.
\end{Theoreme}

Dans cette \'egalit\'e, le membre de gauche est de nature purement
topologique, amorphe, qualitative, celui de droite de nature
g\'eom\'etrique, diff\'erentielle, quantitative. Ainsi le th\'eor\`eme de
Gauss-Bonnet est-il embl\`eme paradigmatique d'une {\sl liaison}, au
sens lautmanien du terme, entre:
\[
\text{\footnotesize\sf Topologie}
\,\,\,\longleftrightarrow\,\,\,
\text{\footnotesize\sf G\'eom\'etrie}.
\]

Cette relation entre topologie et g\'eom\'etrie, explicitant
une synth\`ese entre ces deux disciplines, demande quelques
commentaires que nous esquissons ici. Comment un point de vue
topologique peut-il se mettre en liaison avec un point de vue
g\'eom\'etrique? Avant de donner des \'el\'ements de r\'eponse \`a
cette question, insistons sur le fait que les math\'ematiques se
d\'eveloppent par intervention de disciplines diff\'erentes les unes
sur les autres, \`a diff\'erents titres, outils d'observation,
structuration formelle, nouveau point de vue unifiant. Science en
quelque sorte qualitative de l'espace, la toplogie appara\^it comme
plus primitive que la g\'eom\'etrie, m\^eme si elle vient
historiquement apr\`es. La topologie prend ici la forme de recherche
d'invariants. C'est \`a travers des r\'esultats qui mettent en
\'evidence des invariants que la g\'eom\'etrie est rejointe. Mais le
th\'eor\`eme d'Euler est une relation de nature quantitative entre des
constituants poly\'edriques. Par ce biais, elle se rapproche de
certains aspect de la g\'eom\'etrie. Une quantification minimale de la
vis\'ee topologique va permettre la synth\`ese que propose le
th\'eor\`eme de Gauss-Bonnet. Ce n'est pas ainsi que la
caract\'eristique d'Euler a \'et\'e vue, mais la synth\`ese fournie
par le th\'eor\`eme en laisse voir un nouvel aspect.

C'est un th\'eor\`eme tr\`es avanc\'e dont on ne peut remobiliser
instantan\'ement les divers concepts: $ds^2$, $D$, $C_\ell$,
$\chi(D)$, $K_\GG$, $dA$, $k_g$, $\alpha_{\ell,\maux}$,
lesquels seront pr\'esent\'es en temps voulu.

Tout d'abord, dans le cas planaire $S = \mathcal{P}$, lorsque $D = T$
est simplement un 
triangle comme ci-dessus, la caract\'eristique d'Euler
se calcule ais\'ement:
\[
\aligned
\chi(T)
&
\,=\,
\text{\footnotesize{\sf nombre}}
\big(\text{\footnotesize{\sf sommets}}\big)
-
\text{\footnotesize{\sf nombre}}
\big(\text{\footnotesize{\sf ar\^etes}}\big)
+
\text{\footnotesize{\sf nombre}}
\big(\text{\footnotesize{\sf faces}}\big)
\\
&
\,=\,
3-3+1
\\
&
\,=\,
1,
\endaligned
\]
parce que $T$ constitue une triangulation en elle-m\^eme de
lui-m\^eme, il auto-incarne le quantum en lequel sa forme peut
d\'ecomposer toute forme topologique abstraite, et la complexit\'e
int\'erieure potentielle en laquelle il pourrait s'auto-fragmenter se
resynth\'etise par trivialit\'e homologique: le caract\`ere global de
son essence locale est principe m\^eme de son caract\`ere de
monade. 

\begin{center}
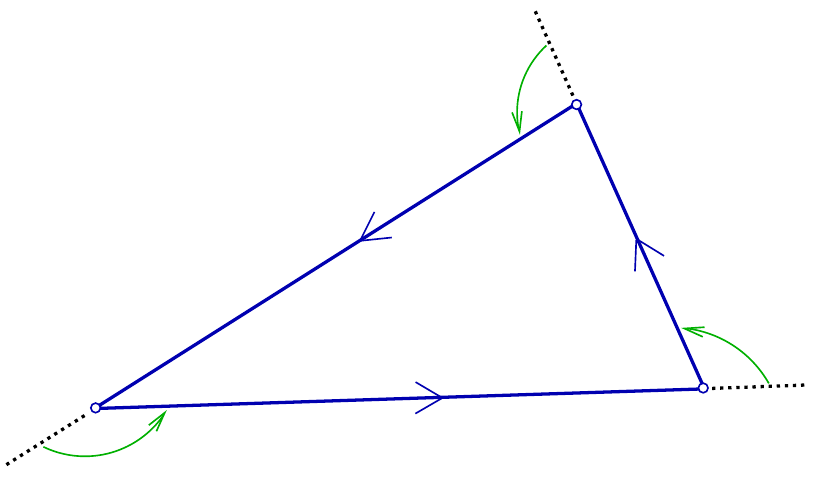
\end{center}

Ensuite, par simple r\`egle des angles compl\'ementaires, les
angles ext\'erieurs en les trois sommets valent:
\[
\pi
-
\widehat{A},
\ \ \ \ \ \ \ \ \ \ \ \ \ \ \ 
\pi
-
\widehat{B},
\ \ \ \ \ \ \ \ \ \ \ \ \ \ \ 
\pi
-
\widehat{C},
\]
et puisque la courbure de Gauss d'un plan $K_\GG \equiv 0$, ainsi
que la courbure g\'eod\'esique de segments de droite
$k_g \equiv 0$, sont toutes deux identiquement nulles, 
la formule de Gauss-Bonnet se <<\,{\sl projette}\,>>:
\[
2\pi
\cdot
1
\,=\,
0
+
0
+
\pi
-
\widehat{A}
+
\pi
-
\widehat{B}
+
\pi
-
\widehat{C},
\] 
pour redonner une formulation \'equivalente au th\'eor\`eme d'Euclide,
\`a un `micro-calcul' pr\`es.

L'id\'ealisation g\'en\'eratrice de contenus math\'ematiques s'articule
\`a une circulation virtuelle possible multiple entre \'etages
de g\'en\'eralit\'e, les deux relations les plus fondamentales \'etant
celle de {\sl v\'erification d'accrochage coh\'erent \`a l'acquis}, 
et celle de {\sl protension \`a l'extension}.

En effet, quelle signification donner \`a un th\'eor\`eme que l'on
retrouve comme corollaire d'un th\'eor\`eme g\'en\'eral dans lequel on
a annul\'e deux termes, ceux qui ne figurent pas dans le th\'eor\`eme
particulier? La r\'eponse est que l'annulation de la valeur d'un terme
ne signifie pas que ce terme n'existe pas. Dire que la courbure est
nulle poss\`ede un sens conceptuel positif. Car la figure
consid\'er\'ee, un plan dont la courbure est nulle, appartient \`a la
cat\'egorie des objets g\'eom\'etriques dont la courbure prend
diff\'erentes valeurs, depuis des valeurs n\'egatives jusqu'\`a des
valeurs positives, en s'annulant au milieu. Le passage au point de
vue de la courbure autorise une extension des objets g\'eom\'etriques
au-del\`a de la g\'eom\'etrie euclidienne. En cons\'equence le
th\'eor\`eme d'Euclide appartient \`a un \'enonc\'e plus g\'en\'eral
qui n'est plus celui confin\'e \`a la g\'eom\'etrie euclidienne. Comme
nous l'avons \'enonc\'e la courbure est donc le concept qui permet
l'extension de la g\'eom\'etrie, il est en quelque sorte un
op\'erateur de r\'eflexion. Il rend possible \`a la mani\`ere d'un
transcendantal le concept qui d\'etermine diff\'erentes 
g\'eom\'etries.

\smallskip

Dans la formule de Gauss-Bonnet, tout contribue \`a l'harmonie
et \`a la v\'erit\'e d'une \'egalit\'e. Entre Euclide et Gauss-Bonnet,
c'est une ascension g\'en\'etique qui {\sl engendre} trois types
de termes: 

\smallskip\noindent$\bullet$\,
la courbure de Gauss $K_\GG$, invariant $2$-dimensionnel; 

\smallskip\noindent$\bullet$\,
la courbure g\'eod\'esique $k_g$, invariant $1$-dimensionnel; 

\smallskip\noindent$\bullet$\,
les discontinuit\'es d'angles ext\'erieurs $\alpha_{\ell,m}$, 
donn\'ees $0$-dimensionnelles.

\smallskip

M\^eme s'il n'y a pas, entre ces trois quantit\'es g\'eom\'etriques:
\[
\sum_{1\leqslant\ell\leqslant\LL}\,
\sum_{1\leqslant\maux\leqslant\MM_\ell}\,
\alpha_{\ell,\maux},
\ \ \ \ \ \ \ \ \ \ \ \ \ \ \
\int_{C_1\cup\cdots\cup C_\LL}\,
k_g\,ds,
\ \ \ \ \ \ \ \ \ \ \ \ \ \ \ 
\int\!\!\int_D\,
K\,dA,
\] 
de correspondance terme \`a terme avec les trois d\'enombrements
(trop variables!) d'une triangulation de $D$ 
({\em voir} ci-dessous): 
\[
\text{\footnotesize{\sf nombre}}
\big(\text{\footnotesize{\sf sommets}}\big),
\ \ \ \ \ \ \ \ \ \ \ \ \ \ \ 
\text{\footnotesize{\sf nombre}}
\big(\text{\footnotesize{\sf ar\^etes}}\big),
\ \ \ \ \ \ \ \ \ \ \ \ \ \ \ 
\text{\footnotesize{\sf nombre}}
\big(\text{\footnotesize{\sf faces}}\big),
\]
il y a {\em \'egalit\'e entre leurs sommes altern\'ees}, \`a un facteur $2\pi$
pr\`es: telle est la nature d'une {\sl d\'ecouverte d'invariant}, par
l'action renouvelable d'une facult\'e universelle de l'entendement
math\'ematique, toujours expos\'e \`a une forme d'insatisfaction dans le fait
de chercher \`a voir au sein de l\`a o\`u rien n'est d'avance visible.

\smallskip

Cette \'egalit\'e poss\`ede de nombreuses significations qui doivent
\^etre analys\'ees, dans le d\'etail de ses termes, dans la
signification de l'\'egalit\'e.

Commen\c{c}ons par la formule d'Euler. 
Elle poss\`ede une histoire qui date d'avant Euler. 
Il en existe plusieurs
dizaines de d\'emonstrations. C'est le premier grand
invariant topologique auquel un math\'ematicien a affaire. La somme
altern\'ee des c\^ot\'es, des faces et des sommets est invariante pour
tout polyg\^one r\'egulier. La plupart des d\'emonstrations proc\`edent
par visualisation des c\^ot\'es, faces et sommets, le poly\`edre
\'etant projet\'e dans le plan. Cette somme porte sur des objets qui
vont de la dimension topologique $1$ 
jusqu'\`a la dimension $2$. L'agencement des figures (poly\`edre) est un montage qui est une
reconstitution \`a partir de ces objets de dimension croissante, pris
dans leur `enti\`eret\'e'. Ce sont l\`a comme les constituants
minimaux d'un squelette de l'espace et donc de figures
(poly\`edres). Au lieu de chercher \`a calculer une longueur, une aire
ou un volume, quantit\'es g\'eom\'etriques, nous nous arr\^etons aux
\'el\'ements `non \'evalu\'es'. Leur combinaison peut \^etre l'objet
d'un r\'esultat alg\'ebrique. M\^eme si le terme de droite est un
nombre, les concepts qui donnent lieu \`a cette \'evaluation
num\'erique font partie des notions
classiques de topologie alg\'ebrique.
Essayons de les caract\'eriser plus pr\'ecis\'ement.

Tout poly\`edre r\'egulier peut \^etre reconstruit \`a partir de ses
constituants d'une mani\`ere qui est li\'ee \`a la structure de
l'espace qui les loge. La longue histoire des cinq poly\`edres
r\'eguliers qui remonte au Tim\'ee et passe par Kepler t\'emoigne de
la conception m\'etaphysique qui a pr\'esid\'e \`a son
\'elaboration. On sait qu'au XIX\textsuperscript{\`eme} 
si\`ecle, l'objet d'\'etude de la
topologie (devenue topologie alg\'ebrique) est non pas le poly\`edre,
mais le complexe simplicial. 
En allant au-del\`a de la dimension $3$, on d\'efinit
les complexes simpliciaux comme form\'es de simplexes soumis \`a
certaines conditions d'incidence. 
Et l'\'etude de ces complexes repose sur la
connaissance des nombres dits de Betti, attach\'es \`a
eux et invariants par toute \'equivalence topologique. 
Comme on le voit dans le cadre
de cette discipline on proc\`ede \`a l'\'etude de certaines
propri\'et\'es de l'espace qui vont dans une direction tr\`es
diff\'erente de celle de la g\'eom\'etrie diff\'erentielle.

Chacun de ces concepts doit \^etre comme moment de ce qui est appel\'e
une ascension, situ\'e en ces \'etages. Le concept de courbure
int\'egrale ou de courbure gaussienne est une d\'ecouverte de
Gauss. Comme nous l'avons not\'e, il change le point de vue sur la
g\'eom\'etrie. Comme `transcendantal', il rend possible {\em a priori}
plusieurs g\'eom\'etries. Il est remarquable que ce soit ce concept
qui soit l'un des op\'erateurs de liaison avec le point de vue de la
topologie (alg\'ebrique). Il est pourtant `perceptible' que la
courbure gaussienne puisse donner des indications sur les relations
entre des \'el\'ements simpliciaux, (pour les poly\`edres, E, F, V).
C'est l'aspect remarquable de ces disciplines, elles laissent 
comme un
\'echo intuitif des \'el\'ements de notre repr\'esentation perceptive
de l'espace.

\smallskip

Un \'etage interm\'ediaire plus accessible\,\,---\,\,lorsqu'on
ne le conna\^it pas d\'ej\`a\,\,---\,\,de cet \'enonc\'e est le cas
`particulier' (d\'ej\`a extr\^emement g\'en\'eral!) o\`u le domaine
$D = S$ est la surface enti\`ere, sans bord, de telle sorte
que deux termes disparaissent dans le membre de droite 
de l'\'egalit\'e. Cette \'etape interm\'ediaire poss\`ede une signification compl\`ete en elle-m\^eme.

 Plus qu'un corollaire, il s'agit d'un th\'eor\`eme
en tant que tel. 

\begin{Theoreme}
\label{theoreme-Gauss-Bonnet-compact}
{\bf [de Gauss-Bonnet Compact]}
La caract\'eristique d'Euler $\chi(S)$ d'une surface compacte
orient\'ee $S$ sans bord vaut:
\[
\chi(S)
\,=\,
\frac{1}{2\pi}
\int\!\!\int_S\,
K_{\sf G}\,
dA.
\]
\end{Theoreme}

Autrement dit, l'invariant topologique le plus \'el\'ementaire
$\chi(S)$ d'une surface vaut, \`a la constante $2\pi$ pr\`es,
la {\sl curvatura integra} de Gauss, \`a savoir l'int\'egrale
compl\`ete de la courbure $K_\GG$ sur $S$.

Ph\'enom\`ene surprenant, cette \'egalit\'e est vraie quelle que soit la
m\'etrique gaussienne lisse $ds^2$ choisie sur $S$\,\,---\,\,il y en a
un tr\`es grand nombre, de cardinal \'egal \`a la puissance du
continu. Si la surface est d\'eform\'ee, sa caract\'eristique
d'Euler, qui est un invariant topologique, ne change pas,
tandis que la courbure en certains points change.
Le Th\'eor\`eme~{\ref{theoreme-Gauss-Bonnet-compact}} 
\'etablit\,\,---\,\,ph\'enom\`ene quelque peu troublant\,\,---\,\,que 
l'int\'egrale de toutes les courbures ne change pas, quelle
que soit la d\'eformation. Par exemple, toute sph\`ere $S \subset
\R^3$ d\'efigur\'ee
par plusieurs `bosses' poss\`ede toujours 
une courbure totale \'egale \`a $4\pi = 2\pi \cdot \chi(S)$,
quel que soit le cabossage. 

Dans la formule g\'en\'erale, \`a cause des termes de bord
$\int_{C_\ell}\, k_g\, ds$ et des termes d'angle $\alpha_\ell$ qui
alourdissaient le travail mental de saisie intuitive compl\`ete
ad\'equate, cette ind\'ependance vis-\`a-vis du choix d'une m\'etrique
n'avait probablement pas \'et\'e vue ou th\'ematis\'ee dans l'horizon
de compr\'ehension. Autrement dit, particularisation et
sp\'ecialisation sont auxiliaires incontournables de mise en
perspective.

Pourquoi en est-il ainsi? Que signifie cette ind\'ependance
vis-\`a-vis d'une m\'etrique? D'abord nous devons remarquer que
tr\`es souvent, \`a cause des trajectoires multiples que permet
l'organisation math\'ematique, nous disposons d'une d\'emultiplication
des significations des concepts math\'ematiques. Ici, l'\'egalit\'e
nous indique qu'un concept de nature topologique exprim\'e par un
nombre caract\'eristique d'Euler est \'egal \`a un nombre qui exprime
un concept de nature g\'eom\'etrique. Cela veut dire que le nombre de
droite poss\`ede une signification topologique (= le nombre de gauche)
et certainement de nature globale. Est-ce l\`a la
raison qui rend la courbure gaussienne ind\'ependante de la
m\'etrique? On sait qu'une m\'etrique repr\'esente les conditions
fix\'ees sur une vari\'et\'e d'apr\`es lesquelles calculer 
une longueur, une aire. 
Et on le sait \'egalement, comment calculer une courbure.

Toute confrontation \`a un \'enonc\'e math\'ematique exige des {\sl
analyses d'int\'egration \`a la pens\'ee de v\'erit\'e}, moment
mentaux r\'eflexifs essentiels sans lesquels il serait impossible que
des individus en formation deviennent des math\'ematiciens, puisque
l'activit\'e `c\'er\'ebrale' intuitive permanente constitue le
n{\oe}ud r\'eel, peu transmissible par \'ecrit, du mouvement de la
pens\'ee. L'activit\'e intellectuelle d'un sujet, m\^eme mise entre
parenth\`eses au nom d'exigences transcendentales d'objectivit\'e
maintes fois promulgu\'ees par la philosophie des math\'ematiques,
demeure une n\'ecessit\'e de fait, pas seulement pratique, puisque la
th\'eorie manque, mais surtout constituante, et aussi partiellement
myst\'erieuse en tant qu'elle s'ex\'ecute en chaque autre sujet
pensant. 

Ce ph\' enom\`ene est l'effet d'une omnisubjectivit\'e dont nous
avons parl\'e. Et corr\'elatif de cette omnisubjectiv\'e, il y a
l'objet construit lui-m\^eme. 

Nous voudrions encore insister sur la
nature du concept de courbure gaussienne, dont nous verrons encore
quelques propri\'et\'es. La courbure gaussienne est un concept qui
`ach\`eve' la pr\'esentation de la sph\`ere, il est une forme de
r\'ealisation de la surface \`a un niveau sup\'erieur. On verra que
cette propri\'et\'e de la courbure gaussienne la rapproche de la
notion d\'ecouverte par Gauss \'egalement, du caract\`ere
intrins\`eque de celle-ci. 

Une \'ecriture \'equivalente du
th\'eor\`eme de Gauss-Bonnet Compact, mais plus `{\sl parfaite}' car
exprimant mieux le rapport aux mod\`eles directeurs, est:
\[
\frac{1}{2}\,
\chi(S)
\,=\,
\int\!\!\int_S\,
K_{\sf G}\,
\frac{dA}{4\pi},
\]
puisque cette \'ecriture signale, dans la fraction $\frac{dA}{4\pi}$,
le lien intime avec la sph\`ere unit\'e dans l'espace euclidien:
\[
S^2
\,:=\,
\big\{
(x,y,z)
\in
\R^3
\colon\,
x^2+y^2+z^2
=
1
\big\},
\]
dont l'aire vaut $4\pi$ d'apr\`es Archim\`ede:
\[
1
\,=\,
\int_{S^2}\,
K_\GG\,
\frac{dA}{4\pi}.
\]

Le lien avec la sph\`ere unit\'e dans l'espace euclidien, signal\' e ici qui montre cette aire \'egale \`a la caract\'eristique est l'expression du fait que (dans ce cas) l'aire est sur les deux terrains de la courbure et de la topologie de la caract\'eristique d'Euler. 
 De mani\`ere analogue, la mesure de
probabilit\'e $\frac{d\theta}{2\pi}$
sur le cercle unit\'e $\R \big/ 2\pi\Z$, de masse totale:
\[
1
\,=\,
\int_{-\pi}^\pi\,
\frac{d\theta}{2\pi},
\]
appara\^it naturellement
dans les formules de calcul des coefficients de Fourier:
\[
\widehat{f}(k)
\,:=\,
\int_{-\pi}^\pi\,
f(\theta)\,
\frac{d\theta}{2\pi}
\eqno
{\scriptstyle{(k\,\in\,\Z)}},
\]
d'une fonction $f \in L^1 \big( \R \big/ 2\pi\Z \big)$ int\'egrable
d\'efinie sur $\R$ et $2\pi$-p\'eriodique. 
M\^eme si le coefficient $\frac{1}{2}\, \chi(S)$ devant la
caract\'eristique d'Euler ne semble pas naturel pour l'instant, 
il deviendra clair plus tard qu'il a lui aussi un sens
normalisateur.

\smallskip

Avant d'entamer des commentaires explicatifs, il importe d'insister
sur la force d'embrassement des \^etres. Dans leur ensemble, les
math\'ematiques expriment un {\sl d\'esir de puissance synth\'etique}
sur les \^etres rationnels. La conceptualisation est une premi\`ere
force d'embrassement, mais les \^etres, toujours, s'y soustraient, en
partie, car le concept, en math\'ematiques, n'est souvent qu'une vue
partielle, d'autres sous-concepts demeurant non-vus, d'autres aspects
intrins\`eques demeurant non-conceptualis\'es. On doit voir la force
d'entra\^inement hors de soi que manifeste la pens\'ee math\'ematique
comme puissance d'expansion, ici propuls\'ee \`a travers cette unit\'e
\`a multiples ressorts de la topologie et de la g\'eom\'etrie.

{\sl Force d'embrassement} peut \^etre dite en termes de puissance
synth\'etique qui pousse au d\'eveloppement et \`a l'unit\'e. Un
ensemble d'unit\'es relationnelles, ind\'efiniment expansives et dont
les éléments sont toujours en voie de renouvellement, donne une forme
anticipatrice aux concepts ainsi obtenus. Caveing analyse cette mise
en {\oe}uvre de la pens\'ee relationnelle que poursuit toute activit\'e
math\'ematique. Nous emploierons \'egalement une notion de champ de
conscience qui a pour contenu la r\'egion du champ transcendantal qui
est th\'ematis\'ee par le math\'ematicien. 
Il est en devenir incessant, ouvert
et temporalis\'e.\footnote{\ Caveing, {\em Le probl\`eme des
objets dans la pens\'ee math\'ematique}, Vrin, Paris, 2004, p.167.}
Les termes comme `{\sl r\'egion}' sont m\'etaphoriques. Le 
{\sl champ de
conscience}, terme emprunt\'e \`a la ph\'enom\'enologie, est la
{\sl r\'egion} du champ transcendantal 
qui constitue son objet actuel, ce
que Cavaill\`es nommait le 
{\sl champ th\'ematique} du math\'ematicien, il
rel\`eve si l'on veut du pour soi, alors que le champ transcendantal
rel\`everait de l'en soi. C'est la conscience singuli\`ere du
math\'ematicien en tant qu'elle est habit\'ee par le langage
sp\'ecifique des math\'ematiques.  Et c'est elle qui s'engage dans
l'exp\'erience pens\'ee du relationnel expansif. C'est sans doute l'un
des aspects essentiels de l'activit\'e math\'ematique que cet
entra\^inement \`a progresser dans l'expansion synth\'etique du
relationnel.

\CITATION{
Retenons alors l'hypoth\`ese que de tout temps, quelle qu'ait \'et\'e
la place du d\'eveloppement historique, ce qui a guid\'e secr\`etement
et silencieusement la recherche math\'ematique, c'est en quelque sorte
la pr\'esence de sch\`emes relationnels \`a valeur anticipatrice
ouvrant en d\'efinitive sur l'explicitation ou l'invention de
relations dans le champ th\'ematique sous 
investigation.\footnotemark}

\footnotetext{\ {\em Ibidem}, p. 259.} 

L'ontologie physique offre un r\'eservoir in\'epuisable de surfaces
visibles qui convainc de l'extension ontologique premi\`ere du
th\'eor\`eme de Gauss-Bonnet: la notion de surface $S$ quelconque,
d\'ej\`a extr\^emement vaste. Ensuite, le choix d'un domaine
quelconque $D \subset S$ sur la surface renforce le sentiment
d'immensit\'e. Enfin, la variabilit\'e de la m\'etrique parach\`eve
cette extensionnalit\'e.

Le sentiment d'immensit\'e est de plusieurs ordres. Il y a celui qui
correspond au sublime dynamique de Kant. Nous sommes en proie \`a ce
sentiment esth\'etique analys\'e par Kant dans la {\em Critique du
Jugement} quand nous sommes math\'ematiciens des surfaces.

Il y a aussi celui que nous \'eprouvons quand nous consid\'erons la
vari\'et\'e des surfaces et l'infinit\'e d'une surface d\'ej\`a en
voie de math\'ematisation\,\,---\,\,ce, toujours dans le cadre des
formes de pens\'ee qui visent au contr\^ole math\'ematique des
infinit\'es. L'analyse de Kant est tr\`es complexe. Il distingue le
sublime math\'ematique et le sublime dynamique. Le sublime dynamique
correspond au sentiment du sublime que nous pouvons \'eprouver dans le
travail math\'ematique. Cette consid\'eration ne peut \^etre
compl\`etement d\'evelopp\'ee dans ce cadre. Ce sentiment est un
sentiment esth\'etique, subjectif. Il est en rapport avec la peur.

\CITATION{
Pour la facult\'e de juger esth\'erique r\'efl\'echissante la nature
ne peut poss\'eder une valeur <gelten> comme force, \^etre sublime
dynamiquement, que dans la mesure o\`u elle est consid\'er\'ee comme
objet de peur.\footnotemark}

\footnotetext{\ Kant, {\em Critique de la facult\'e de juger},
Trad. A. Philonenko, Vrin, Paris, 1968, p.~99.}

\noindent
Le champ math\'ematique pre\'sente ainsi une forme de sublime
apparent\'ee \`a celle de la nature. Il est certes pour une part
importante l'{\oe}uvre des math\'ematiciens, esprits humains. Mais il
ouvre aussi au spectacle de l'inach\`evement et de ses
d\'eflagrations. 

\CITATION{
Nous avons trouv\'e notre limite propre en ce qui est incommensurable
dans la nature et dans l'incapacit\'e de notre facult\'e \`a saisir
une mesure proportionn\'ee \`a l'\'evaluation esth\'etique de
la grandeur de son domaine, et cependant aussi d\'ecouvert en notre
raison, une autre mesure non sensible, qui comprend
sous elle comme unit\'e cette infinit\'e, par rapport \`a laquelle
tout dans la nature est petit, si bien que nous avons d\'ecouvert en
notre esprit une sup\'eriorit\'e sur la nature en son immensit\'e
[\dots]. Sa force irr\'esistible nous fait conna\^itre, en tant
qu'\^etres de la nature notre faiblesse physique, mais en
m\^eme temps, elle d\'evoile une facult\'e qui nous permet de nous
consid\'erer comme ind\'ependants par rapport \`a elle et une
sup\'eriorit\'e sur la nature sur laquelle se fonde une conservation
de soi-m\^eme toute diff\'erente de celle qui est attaqu\'ee par la
nature qui nous est ext\'erieure [\dots].
En ce sens, la nature n'est pas
consid\'er\'ee comme sublime dans notre jugement esth\'etique parce
qu'elle engendre la peur mais parce qu'elle constitue un appel
\`a la force qui est en nous (mais qui n'est pas
nature).\footnotemark}

\footnotetext{\ {\em Ibidem}.}

Les \'elaborations math\'ematiques dans leur
profondeur parfois inou\"{\i}e pr\'esentent des paysages mais aussi
toujours un au-del\`a de ces \'etendues induisant un sentiment
d'effroi. On pourrait qualifier ces \'etendues et leur au-del\`a du
vocable d'{\em immensit\'e sans cesse potentialis\'ee}.

\begin{These}
Les math\'ematiques \'elaborent des formes de pens\'ee
sur le Vaste des immensités
abstraites, elles-m\^emes \'evolutives.
\end{These}

Le rapport m\'etaphysique que les math\'ematiques d\'eveloppent
face au monde vaste du synth\'etique {\em a priori}
qu'elles cr\'eent de mani\`ere immanente, 
ne s'apparente pas seulement \`a l'inachevabilit\'e
de l'exploration astrophysique, il est engag\'e dans un destin
de m\'etamorphose de la nature des synth\`eses face \`a du Vaste
qui demeure, \`a cause de l'ind\'efinitude de l'infini,
insaisissable en totalit\'e. 

\Section{\bf La courbure locale intrins\`eque}
\label{courbure-locale-intrinseque}
\HEAD{\ref{courbure-locale-intrinseque}.~{\sf 
La courbure locale intrins\`eque}
}{
Jo\"el {\sc Merker} et Jean-Jacques {\sc Szczeciniarz}}

Un des concepts nouveaux de la g\'eom\'etrie contemporaine est celui
qui enveloppe la diff\'erence entre local et global. Il est le
r\'esultat d'une longue r\'eflexion sur les conditions de
l'\'elaboration math\'ematique, en particulier en topologie mais elle
concerne toutes les disciplines math\'ematiques. Le local d\'elimite
ou particularise l'espace de validit\'e d'un concept, pour le
pr\'eciser en mettant en suspens la question de sa globalisation \`a
tout le domaine o\`u il est construit. Il est au c{\oe}ur
de la signification du th\' eor\`eme de Gauss-Bonnet qui pr\'esente
une version locale et une version globale. Une fa\c{c}on de
l'expliciter consiste d'abord \`a faire appel \`a des concepts
alg\'ebriques qui structurent cette diff\'erence dans le cas de la
courbure des surfaces intrins\`eques.

Soit $J$ un transport parall\`ele (notion d\'efinie ci-dessous) sur
une surface abstraite $S$ munie d'une {\sl m\'etrique riemannienne},
laquelle associe \`a toute paire de champs de vecteur sur $S$ une
fonction:
\[
p
\,\,\longmapsto\,\,
\langle X_p,Y_p\rangle
\eqno
{\scriptstyle{(p\,\in\,S)}},
\]
r\'ealisant le {\sl produit scalaire} entre les vecteur
$X_p$ et $Y_p$ aux divers points $p \in S$. Bien entendu,
cette application doit \^etre sym\'etrique $\langle Y, X
\rangle = \langle X, Y \rangle$ et 
bilin\'eaire par rapport
\`a chacun de ses arguments:
\[
\big\langle
f_1X_1+f_2X_2,\,\,Y\big\rangle
\,=\,
f_1\,\langle X_1,\,Y\rangle
+
f_2\,\langle X_2,\,Y\rangle,
\]
quel que soit le choix de fonctions lisses $f_1, f_2 \in 
\mathcal{C}^\infty(S)$. 

Ainsi, une m\'etrique riemannienne est une fa\c{c}on de se donner des
conditions d'effectuation de mesures, qui lui sont donc
relatives. Dans cette mise en valeur des conditions de la mesure nous
poursuivons un travail d'objectivation et de constitution de
l'objectivit\'e. Cependant, 
si l'on consid\`ere une surface $S$ munie d'une
m\'etrique riemannienne $g$, on ne sait pas en g\'en\'eral comparer
des vecteurs tangents en des points diff\'erents de la surface. Une
structure suppl\'ementaire est alors n\'ecessaire.

On appelle {\sl transport parall\`ele sur $S$} la donn\'ee d'une
application:
\[
J
\colon
\ \ \
\big(\gamma,t,t'\big)
\,\,\longrightarrow\,\,
J^{\gamma}_{t't}
\]
qui, \`a toute courbe diff\'erentiable $\gamma: I \longrightarrow
V$ dans $S$ d\'efinie sur
un intervalle ouvert $I \subset \R$, 
et \`a tout couple $(t, t')$ de valeurs du param\`etre
dans l'intervalle $I$, associe une {\em isom\' etrie} des espaces
(vectoriels) tangents:
\[
J^{\gamma}_{t't}
\colon\ \ \
T_{\gamma(t)}S
\overset{\sim}{\,\,\longrightarrow} 
T_{\gamma(t')}S,
\]
lesquels sont munis des produits scalaires issus
de la m\'etrique riemannienne. 
Cet isomorphisme est alors appel\'e 
{\sl transport parall\`ele le long de $\gamma$}. 

\begin{center}
\includegraphics[width=9cm]{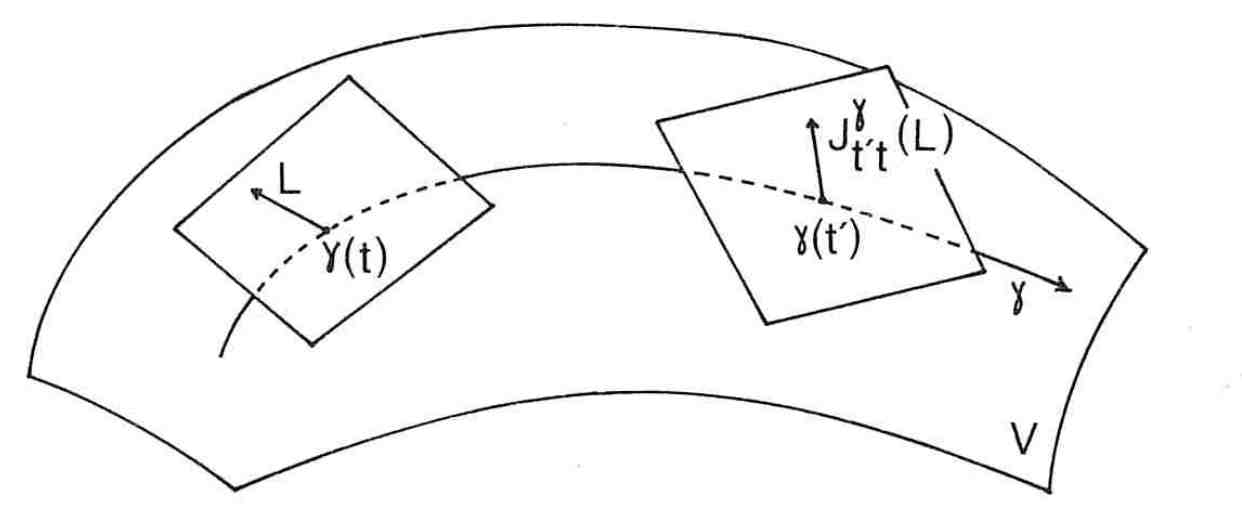}
\end{center}

C'est l\`a un concept-clef que nous devons voir comme une construction
de l'objectivit\'e de la courbure par son enregistrement dans le
d\'eplacement des espaces tangents. Nous sommes descendus un peu plus
profond\'ement dans l'\'enonciation de conditions de constitution de
l'objet g\'eom\'etrique, en \'etablissant une nouvelle contrainte
impos\'ee par la comparaison de vecteurs tangents en des points
distincts d'une vari\'et\'e. C'est ce que nous continuerons \`a
d\'evelopper en analysant le th\'eor\`eme de Gauss-Bonnet. Mais ici,
le transport parall\`ele nous montre comment des conditions de
parall\'elisme doivent \^etre construites quand nous {\em comparons}
des vecteurs tangents en des points distincts d'une vari\'et\'e. C'est
en ce sens que le transport explicite leur objectivit\'e.

On notera dor\'enavant $\mathcal{X}(S)$ l'espace des champs de
vecteurs diff\'erentiables sur $S$. Un champ de vecteurs, notion qui
vient certainement de la physique, est une application qui, \`a chaque
point d'un espace ou d'une vari\'et\'e, associe un vecteur. C'est une
mani\`ere de munir l'espace d'une potentialit\'e d'action 
r\'epartie de mani\`ere uniforme sur tout son domaine.

Pour tout champ $X \in \mathcal{X}(S)$ et
pour tout nombre $t_{0} \in I$, l'application $t \longmapsto
J^{\gamma}_{t_0, t}(X_{\gamma(t)})$ de $I$ dans
l'espace tangent fix\'e
$T_{\gamma(t_0)}(S)$ est diff\'erentiable.
On demande alors que sa d\'eriv\'ee:
\[
\Big(
\frac{d}{dt}
\Big)_{t_{0}}
\big[
J^{\gamma}_{t_{0}t}(X_{\gamma(t)})
\big]
\,\,\in\,\,
T_{\gamma_{(t_0)}}S
\]
ne d\'epende pas des distorsions \'eventuelles d'ordre sup\'erieur
de $\gamma$, mais seulement de son vecteur tangent 
$L := \frac{d\gamma}{dt}(t_0)$ au point $\gamma(t_{0})$. On note
$\nabla _{L}X$ cette d\'eriv\'ee, et on l'appelle {\sl d\'eriv\'ee
covariante} du champ de vecteurs $X$ par rapport \`a $L$.

\smallskip

Un th\'eor\`eme math\'ematique
pr\'ecis d\'emontre qu'\'etant donn\'e un tel transport parall\`ele
abstrait $J_{t',t}^\gamma(\cdot)$, pour tous champs de vecteurs
$X_\smallbullet, Y_\smallbullet \in \mathcal{X}(S)$
sur la surface, et pour toutes
fonctions lisses $f, g \in \mathcal{C}^\infty(S)$,
les $5 = 4 + 1$ propri\'et\'es suivantes sont
satisfaites:

\smallskip\noindent{\bf (1)}\,
$\nabla_X\big(Y_1+Y_2\big) = \nabla_X(Y_1) + \nabla_X(Y_2)$,

\smallskip\noindent{\bf (2)}\,
$\nabla_{X_1+X_2}(Y) = \nabla_{X_1}(Y) + \nabla_{X_2}(Y)$,

\smallskip\noindent{\bf (3)}\,
$\nabla_{gX}(Y) = g\, \nabla_X(Y)$,

\smallskip\noindent{\bf (4)}\,
$\nabla_X\big(f\, Y\big) = f\, \nabla_X(Y) + 
X(f)\,Y$,

\smallskip\noindent
o\`u $X(f)$ d\'esigne l'{\em action} du champ de vecteurs $X$,
envisag\'e comme {\em op\'erateur de diff\'erentiation d'ordre $1$},
sur la fonction $f$. Il est bien remarquable que le champ soit un
pourvoyeur d'action, ici de diff\'erentiation. Les concepts (ou les
structures) math\'ematiques agissent sur les objets. Et dans le cas
qui nous occupe, cette action sur les objets (les fonctions) est
indiff\'eremment r\'epartie sur tout l'espace.

La d\'eriv\'ee covariante doit respecter les r\`egles de d\'erivation,
car elle a int\'egr\'e \`a un niveau sup\'erieur le concept de
d\'erivation qu'elle adapte \`a une situation plus g\'en\'erale. La
nouvelle th\'eorie avec ses nouvelles situations est th\' eoriquement
homog\`ene, elle doit adapter tout l'arsenal
ancien de la d\'erivation
\`a cette situation, en particulier elle doit adapter un \'el\'ement
de lin\'earit\'e par rapport \`a la direction. 
Autrement dit, lorsque les concepts de la g\'eom\'etrie
diff\'erentielle se d\'eveloppent c'est une nouvelle forme de
d\'erivation qui prolonge et transforme ad\'equatement l'ancienne qui
se met en place. 

Ici, la condition {\small\bf (3)} exprime une
$\mathcal{C}^\infty(S)$-lin\'earit\'e par rapport \`a la {\em
direction} $X$ le long de laquelle $\nabla_X (\,
\centersmallbullet)$ d\'erive.

Mais comme c'est r\'eellement $Y$ que $\nabla_X
(\, \centersmallbullet)$ d\'erive, une loi de produit 
inspir\'ee de la r\`egle \'el\'ementaire de diff\'erentiation 
d'un produit\,\,---\,\,classiquement attribu\'ee \`a Leibniz\,\,---:
\[
\big(f\,g\big)'
\,=\,
f'\,g
+
f\,g',
\]
doit \^etre satisfaite, ce qu'exprime {\small\bf (4)} \`a un
niveau plus sophistiqu\'e.

De plus, l'hypoth\`ese faite ci-dessus que le transport
parall\`ele respecte la structure riemannienne, c'est-\`a-dire
les produits scalaires infinit\'esimaux, se r\'e-exprime comme une
{\em cinqui\`eme} condition de compatibilit\'e avec la m\'etrique:

\smallskip\noindent{\bf (5)}\,
$X\big( \langle X,\,Y \rangle \big) \,=\, \big\langle \nabla_X(Y),\,Z
\big\rangle + \big\langle Y,\,\nabla_X(Z)
\big\rangle$.

\smallskip\noindent
On appelle alors $\nabla_X(Y)$ {\sl d\'eriv\'ee covariante
de $Y$ dans la direction $X$}.

\smallskip

Inversement\,\,---\,\,et c'est toute la force du sous-bassement
synth\'etique et interrogatif des math\'ematiques, toujours
en qu\^ete de reformulations et d'{\sl \'equivalences avan\c{c}antes}\,\,---, 
on d\'emontre en g\'eom\'etrie diff\'erentielle qu'un concept abstrait
de d\'eriv\'ee covariante pour lequel les propri\'et\'es 
{\small\bf (1)} \`a {\small\bf (5)} sont transmu\'ees
en axiomes, engendre, gr\^ace \`a un proc\'ed\'e ordinaire d'int\'egration,
des transports parall\`eles $J_{t', t}^\gamma(\, \centersmallbullet)$ 
le long de toutes les courbes $\gamma \subset S$. 

L'\^etre abstrait d'un op\'erateur $\nabla_\smallbullet(\, 
\centersmallbullet)$ de d\'erivation covariante est manifestement
{\em diff\'erentiel}, donc en particulier local. En lui-m\^eme,
$\nabla_\smallbullet(\, 
\centersmallbullet)$ interpr\`ete la correction que l'on 
doit apporter, au plan infinit\'esimal, aux translations
euclidiennes na\"{\i}ves afin que le concept prenne tout son
sens math\'ematique. 

Parce qu'il ne font que diff\'erentier des objets fonctionnels
diff\'erentiables donn\'es,
les op\'erateurs de type $\nabla_\smallbullet(\, 
\centersmallbullet)$ sont intrins\`equement {\em plus simples}
que les applications de transport 
$J_{t', t}^\gamma(\, \centersmallbullet)$. 

Ainsi, la th\'eorie
math\'ematique \'etablit qu'il y a essentiellement \'equivalence,
dans leurs incarnations conceptuelles effectives, entre:
\[
\text{\footnotesize\sf
Id\'ee de transport parall\`ele}
\ \ \ \ \
\Longleftrightarrow
\ \ \ \ \
\text{\footnotesize\sf
Id\'ee de d\'erivation covariante}.
\]
De plus, elle d\'ecide d'{\em articuler} la progression th\'eorique
\`a l'\'el\'ementarit\'e r\'eelle, f\^ut-elle constat\'ee
{\em a posteriori}.

La th\'eorie math\'ematique v\'erifie aussi que le transport parall\`ele ne
d\'epend pas du choix d'un param\'etrage (quelconque) des
courbes. Cependant, en lui-m\^eme, le transport {\em d\'epend} du chemin
$\gamma$ emprunt\'e, et c'est
en cela que r\'eside la diff\'erence discriminante majeure avec la
th\'eorie euclidienne.

L'int\'er\^et des axiomes {\small\bf (1)} \`a {\small\bf (5)} est de
d\'efinir des r\`egles d'usage, transf\'er\'ees et ins\'er\'ees dans
le cadre des surfaces diff\'erentiables, bien au-del\`a du cadre
euclidien familier. Mais on peut d'ailleurs parfaitement voir comment
les r\`egles usuelles de d\'erivation dans $\R^2$ deviennent un cas
particulier restreint de ces conceptions. Il en va de m\^eme d'un
grand nombre de concepts \'el\'ementaires de la g\'eom\'etrie
euclidienne, qui deviennent, dans un regard th\'eorique a
posteriorique, cas particuliers de concepts sup\'erieurs. La pens\'ee
compl\`ete de l'objet math\'ematique exige en effet des
re-circulations synth\'etiques entre th\'eories engendr\'ees et
th\'eories g\'enitrices.

La g\'eom\'etrie initiale a ainsi jou\'e un double r\^ole. Elle a
d'une part \'et\'e la th\'eorie initiale \`a partir de laquelle r\'
efl\'echir \`a de nouvelles situations: elle port(ait) en elle des
g\'en\'eralisations potentielles de concepts dont elle permet
d'extraire une signification plus profonde. D'autre part, elle
<<\,retombe\,>> `a l'\'etat de cas particuliers 
des cas g\'en\'eraux qui l'ont port\'ee au-del\`a d'elle-m\^eme. 
C'est ce que signifie ici l'expression 
<<\,{\sl recirculations synth\'etiques}\,>>.

\Section{\bf $1$-forme diff\'erentielle associ\'ee \`a une diff\'erentiation
covariante}
\label{1-forme-differentiation-covariante}
\HEAD{\ref{1-forme-differentiation-covariante}.~{\sf $1$-forme 
diff\'erentielle associ\'ee \`a une diff\'erentiation
covariante}
}{
Jo\"el {\sc Merker} et Jean-Jacques {\sc Szczeciniarz}}

Maintenant, puisque la métrique riemannienne $\langle 
\centersmallbullet, \centersmallbullet \rangle$ sur $S$ a
\'et\'e conceptuellement cr\'e\'ee {\em pour} donner un sens au produit
scalaire entre paires de vecteurs, donc en particulier \`a leur
orthogonalit\'e, il est {\em naturel}, comme en g\'eom\'etrie
euclidienne, de rapporter les vecteurs \`a des bases
{\em orthonorm\'ees}. 

Sur un ouvert $U$ contenu dans $S$, envisag\'e comme `morceau'
(essentiellement local) de $S$, si deux champs de vecteurs
quelconques $X_1, Y_1 \in \mathcal{X}(S)$
sont donn\'es et sont lin\'eairement ind\'ependants, le proc\'ed\'e
de Gram-Schmidt permet de leur substituer deux autres
champs de vecteurs:
\[
\aligned
X
&
\,:=\,
\frac{1}{\sqrt{\langle X_1,X_1\rangle}}\,\,
X_1,
\\
Y
&
\,:=\,
\frac{\sqrt{\langle X_1,X_1\rangle}}{
\sqrt{
\langle X_1,X_1\rangle\,
\langle Y_1,Y_1\rangle
-
\langle X_1,Y_1\rangle^2}}\,\,
\bigg(
-\,
\frac{\langle X_1,Y_1,\rangle}{
\langle X_1,X_1\rangle}\,X_1
+
Y_1
\bigg),
\endaligned
\]
qui deviennent de norme $1$, et mutuellement orthogonaux,
en {\em tout} point de $U \subset S$:
\[
1
\,\equiv\,
\langle X,X\rangle,
\ \ \ \ \ \ \ \ \ \ \ \ \ \ \ \ \ \
0
\,\equiv\,
\langle X,Y\rangle
\,=\,
\langle Y,X\rangle
\ \ \ \ \ \ \ \ \ \ \ \ \ \ \ \ \ \
1
\,\equiv\,
\langle Y,Y\rangle.
\]

L'int\'er\^et de cette r\'eduction \`a un champ de rep\`eres
orthogonaux des espaces tangents:
\[
T_pS
\,=\,
\R\,X_p
\overset{\bot}{\,\oplus\,}
\R\,Y_p
\eqno
{\scriptstyle{(\forall\,p\,\in\,U)}},
\]
va au-del\`a d'une analogie entre la structure tangentielle de $S$ et
la g\'eom\'etrie euclidienne standard. En effet, la r\'eduction \`a
une orthogonalit\'e simplifie substantiellement la {\em forme
alg\'ebrique} du concept de d\'erivation covariante, et au-del\`a
encore, va m\^eme permettre de produire une formule int\'egrale
explicite du transport parall\`ele associ\'e.

Plus pr\'ecis\'ement, en utilisant le fait qu'une diff\'erentiation
des identit\'es:
\[
1
\,\equiv\,
\langle X,X\rangle,
\ \ \ \ \ \ \ \ \ \ \ \ \ \ \ \ \ \
0
\,\equiv\,
\langle X,Y\rangle,
\ \ \ \ \ \ \ \ \ \ \ \ \ \ \ \ \ \
1
\,\equiv\,
\langle Y,Y\rangle,
\]
au moyen d'un champ de vecteurs quelconque $Z \in \mathcal{X}(S)$
produit $0$ \`a gauche:
\[
0
\,\equiv\,
2\,
\big\langle
\nabla_Z(X),\,X
\big\rangle
\,\equiv\,
\big\langle
\nabla_Z(X),\,Y
\big\rangle
+
\big\langle
X,\,\nabla_Z(Y)
\big\rangle
\,\equiv\,
2\,
\big\langle
\nabla_Z(Y),\,Y
\big\rangle,
\]
on d\'emontre par un raisonnement math\'ematique rigoureux que la
diff\'erentiation covariante doit n\'ecessairement agir au niveau
infinit\'esimal comme une {\em rotation}, d'angle $-\frac{\pi}{2}$,
suivie d'une certaine dilatation qui d\'epend du point. 

C'est-\`a-dire que dans un champ de rep\`eres orthonorm\'es
$\big\{ X_p, Y_p \big\}_{p\in U}$, pour tout champ
$Z \in \mathcal{X}(U)$, on a des formules:
\[
\aligned
\nabla_Z(X)
&
\,=\,
-\,\omega(Z)\,
Y,
\\
\nabla_Z(Y)
&
\,=\,
\omega(Z)\,
X,
\endaligned
\]
qui expriment comment les vecteurs de la base 
ambiante sont diff\'erenti\'es.
Ici, la rotation d'angle $-\frac{\pi}{2}$ envoie
$X$ sur $-Y$ et $Y$ sur $X$. 
Ici, le facteur de dilatation en question est rep\'er\'e
au moyen d'une certaine {\sl $1$-forme diff\'erentielle:}
\[
\omega
\,\in\,
\Gamma\big(U,\,T^\ast X\big).
\] 

Une forme diff\'erentielle
générale $\alpha$ sur $U$, rappelons-le bri\`evement, 
est une application $\mathcal{C}^\infty(U)$-linéaire:
\[
\alpha
\colon\ \ \
\mathcal{X}(U)
\,\longrightarrow\,
\mathcal{C}^\infty(U)
\]
 du 
$\mathcal{C}^\infty$-module $\mathcal{X}(U)$ des champs de vecteurs
diff\'erentiables sur $U$ à valeurs dans la $\R$-alg\`ebre
$\mathcal{C}^\infty(U)$ des fonctions indéfiniment diff\'erentiables
sur $U$.  Dans cette premi\`ere \'etape de la pr\'esentation, nous
avons trouv\'e une synth\`ese interm\'ediaire en cette application
{\em linéaire} $\omega$ qui va des champs de vecteurs vers les
fonctions diff\'erentiables: ce que nous gagnons \`a travers cette
synth\`ese, c'est la lin\'earit\'e.

Cette $1$-forme de connexion $\omega$ est un nouveau concept
qui retranscrit alors une abstraction de la notion de connexion
comme collection de formes lin\'eaires sur les
espaces tangents \`a $S$. Autrement dit, en tout point $p \in U$,
la $1$-forme $\omega$ associe une forme lin\'eaire:
\[
\omega_p
\colon\ \ \ 
T_pS
\,\longrightarrow\,
\R,
\]
de telle fa\c{c}on que pour tout champ $Z \in \mathcal{X}(U)$, 
la fonction $p \longmapsto \omega_p(Z_p)$ soit 
diff\'erentiable.

Cette $1$-forme $\omega$, uniquement d\'etermin\'ee par 
$\nabla_\smallbullet (\, \centersmallbullet)$,
est appel\'ee {\sl forme} (locale) {\sl de connexion} associ\'ee
au champ de rep\`eres orthonorm\'es.

Ici, donc, le concept subit une transformation avantageuse.
Le nombre de fonctions au moyen desquelles il s'exprime
est r\'eduit. Sa signification g\'eom\'etrique s'en trouve
\'eclaircie.

Dans de telles bases locales orthonorm\'ees, les transports
parall\`eles:
\[
J_{t',t}^\gamma
\colon\ \ \ 
T_{\gamma(t)}S
\,\,\longrightarrow\,\,
T_{\gamma(t')}S
\]
deviennent alors de simples {\em rotations}\,\,---\,\,d\'ependant
de $\gamma$, de $t$, et de $t'$\,\,---, puisque la structure
m\'etrique est respect\'ee. 
Par cons\'equent, l'action lin\'eaire de $J_{t',t}^\gamma$
ne s'exprime {\em pas} comme une matrice $2 \times 2$
quelconque dans les bases:
\[
\big\{
X_{\gamma(t)},\,\,
Y_{\gamma(t)}
\big\}
\ \ \ \ \ \ \ \ \ \ \ \ \ \ \ \ \ \
\text{et}
\ \ \ \ \ \ \ \ \ \ \ \ \ \ \ \ \ \
\big\{
X_{\gamma(t')},\,\,
Y_{\gamma(t')}
\big\},
\]
mais comme une {\em matrice de rotation:}
\[
J_{t',t}^\gamma
\,\,=\,\,
\left(\!
\begin{array}{cc}
\cos\,\theta^\gamma(t',t)
&
-\sin\,\theta^\gamma(t',t)
\\
\sin\,\theta^\gamma(t',t)
&
\cos\,\theta^\gamma(t',t)
\end{array}
\!\right),
\]
pour un certain angle $\theta^\gamma(t',t) \in \R \big/
2\pi\Z$.

Les axiomes naturels d'inversion et de transitivit\'e
des sens de parcours imposent des relations de Chasles:
\[
\theta^\gamma(t,t)
\,=\,0
\ \ \ \ \ \ \ \ \ \ \ \ \ \ \ \ \ \
\text{et}
\ \ \ \ \ \ \ \ \ \ \ \ \ \ \ \ \ \
\theta^\gamma(t'',t')
\,=\,
\theta^\gamma(t'',t')
+
\theta^\gamma(t',t).
\]

Ensuite, une analyse de pens\'ee convainc ais\'ement
que si un transport parall\`ele $J_{\smallbullet, 
\smallbullet}^\smallbullet$ existe dans un champ de rep\`eres 
orthonorm\'es, alors son angle de rotation est
n\'ecessairement reli\'e \`a la forme diff\'erentielle fondamentale
par l'\'equation:
\[
\omega
\Big(
\frac{d\gamma}{dt}(t)
\Big)
\,=\,
-\,
\frac{\partial\theta^\gamma(t',t)}{\partial t'}
\bigg\vert_{t'=t}.
\]

Toujours, donc, avec une paire $\{ X, Y\}$ de champs
de vecteurs orthonorm\'es locaux d\'efinis dans un ouvert $U\subset S$:
\[
1
\,\equiv\,
\langle X,X\rangle,
\ \ \ \ \ \ \ \ \ \ \ \ \ \ \ \ \ \
0
\,\equiv\,
\langle X,Y\rangle
\,=\,
\langle Y,X\rangle
\ \ \ \ \ \ \ \ \ \ \ \ \ \ \ \ \ \
1
\,\equiv\,
\langle Y,Y\rangle.
\]
une simple int\'egration de cette relation diff\'erentielle
permet de reconstituer l'angle induit par un transport
parall\`ele \`a partir du constituant fondamental $\omega$
de l'op\'erateur de diff\'erentiation covariante:
\[
\theta^\gamma(t',t)
\,=\,
\int_t^{t'}\,
\omega
\Big(
\frac{d\gamma}{du}(u)
\Big)\,
du.
\]

\smallskip

En r\'esum\'e\,\,---\,\,et en conclusion transitoire\,\,---,
le parcours de pens\'ee math\'ematique a \'et\'e le suivant.

\smallskip\noindent\blue{$\square$}\,
\'Elaboration d'un concept abstrait de transport parall\`ele.

\smallskip\noindent\blue{$\square$}\,
Naissance d'un concept infinit\'esimal de diff\'erentiation covariante.

\smallskip\noindent\blue{$\square$}\,
Normalisation et simplification des \'equations math\'ematiques 
caract\'eristiques.

\smallskip\noindent\blue{$\square$}\,
R\'einterpr\'etation infinit\'esimale de l'\'ecart-transport.

\smallskip\noindent\blue{$\square$}\,
R\'eint\'egration globale de l'\'ecart-transport infinit\'esimal.

\medskip

En math\'ematiques, le savoir absolu, c'est la compr\'ehension
des synth\`eses, au-del\`a et par-del\`a le labyrinthe symbolique
des d\'emonstrations d\'etaill\'ees. Toujours, les concepts sont
travaill\'es par des \'equivalences potentielles ou actuelles,
puisqu'ils portent en eux un gradient transformationnel,
avec une grande part d'irr\'eversibilit\'e synth\'etique, 
t\'emoin d'avancement de connaissance. 

Ici, en g\'eom\'etrie {\em analytique et diff\'erentielle}, la plupart
des conceptions s'enracinent dans l'infinit\'esimal. Les
caract\'eristiques m\^eme des \^etres s'expriment selon certaines
modalit\'es fonctionnelles sp\'ecifiques qui int\`egrent ce qu'ils
sont r\'eellement. La pr\'esence de l'infinit\'esimal dans la
structuration th\'eorique, au fond, n'est due qu'\`a sa force de
donation instantan\'ee. L'objet, dans sa production intellectuelle, se
r\'eoriente, comme par un magn\'etisme, dans les directions
synth\'etiques ad\'equates.

Pourquoi l'infinit\'esimal est-il au fond de tous ces concepts?
D'abord, on sait qu'il est \`a la base m\'etaphysique de toute
conception de la notion de diff\'erentielle. Depuis le 
XVII\textsuperscript{ème} 
si\`ecle, l'infinit\'esimal est au c{\oe}ur de la r\'eflexion
physico-math\'ematique. Mais il s'est propag\'e \`a tous les concepts
introduits en g\'eom\'etrie {\em diff\'erentielle}. Et ce, m\^ eme
quand une contruction alg\'ebrique vient 
gouverner les structures qui
organisent la g\'eom\'etrie, comme c'est par exemple le cas pour les
formes diff\'erentielles. Bien plus, 
ces d\'eveloppements lui donnent
de nouvelles base d'appui, \`a travers lesquelles il se r\'ealise. On
doit noter la puissance synth\'etique du th\`eme infinit\'esimal qui
poss\`ede ce pouvoir de se lier \`a l'alg\'ebrique sous diverses
formes.

Il est \'etonnant que le transport parall\`ele puisse se r\'eduire,
comme nous l'avons vu plus haut, \`a une rotation
suivie d'une dilatation. Mais
pr\'ecis\'ement, une compr\'ehension en profondeur du transport
parall\`ele implique le concept de rotation. C'est la $1$-forme dite
{\em $1$-forme de connexion} qui a permis de l'introduire
effectivement. 

Au niveau de l'individualit\'e pensante, du sujet math\'ematique,
du je en apprentissage, tout cela demeure souvent invisible. 

\Section{\bf Holonomie du transport parall\`ele}
\label{holonomie-transport-parallele}
\HEAD{\ref{holonomie-transport-parallele}.~{\sf 
Holonomie du transport parall\`ele}
}{
Jo\"el {\sc Merker} et Jean-Jacques {\sc Szczeciniarz}}

En g\'en\'eral, le transport parall\`ele d'un vecteur $Z_p \in 
T_p S$ le long d'une courbe ferm\'ee partant d'un point 
$p \in S$ et y revenant donne des vecteurs diff\'erents,
y compris {\em a fortiori} lorsqu'on
change de courbe. Autrement dit, si $\lambda \colon [0, 1]
\longrightarrow S$ avec $\lambda(0) = p = \lambda(1)$
est une telle courbe,
alors en g\'en\'eral:
\[
J^\lambda\big(Z_p\big)
\,\neq\,
Z_p,
\]
pour la plupart des lacets $\lambda$ et des vecteurs $Z_p \in 
T_p S$.

Commen\c{c}ons par penser cette possibilit\'e au moyen d'une 
conceptualisation minimale, qui laisse implicites de nombreuses
questions. En math\'ematiques, fr\'equemment, les concepts
ne sont qu'expression de potentialit\'es, ils 
conceptualisent en anticipation, et non en cl\^oture: c'est
toute la force de leur flexibilit\'e, et aussi, la faiblesse
de leurs accomplissements, aspect d'ouverture 
fr\'equemment invisible
aux philosophies des math\'ematiques.

\`A cet effet, introduisons une notation
pour d\'esigner l'ensemble des lacets en question:
\[
\Omega_p(S)
\,:=\,
\Big\{
\text{\rm lacets diff\'erentiables}\,\,
\lambda\colon\,[0,1]
\longrightarrow
S,\,\,
\lambda(0)=p=\lambda(1)
\ \
\text{\rm modulo homotopie}
\Big\}.
\] 
Les lacets peuvent \^etre compos\'es. En effet, si $\lambda_1$ et
$\lambda_2$ sont deux lacets en $p$, on d\'efinit:
\[
\lambda_1
\cdot
\lambda_2
\,=:\,
\lambda,
\]
comme \'etant le lacet parcourant d'abord $\lambda_2$, puis
$\lambda_1$, et on renormalise le param\'etrage \`a \^etre encore
l'intervalle $[0,1]$. La relation d'{\sl homotopie}
considère alors comme équivalents deux lacets qui peuvent
être déformés continûment l'un vers l'autre.
D'ailleurs, pour des lacets (diff\'erentiables), 
$\Omega_p(S)$ est un 
groupe (abstrait), appel\'e {\sl groupe de Poincar\'e}
(diff\'erentiable) {\sl de $S$}. 

\`A tout lacet $\lambda \in \Omega_p(S)$ est associ\'e un
transport parall\`ele:
\[
J^\lambda
\colon\ \ \
T_pS
\underset{\bot}{
\overset{\sim}{\,\,\longrightarrow\,\,}}
T_pS,
\]
qui est un isomorphisme lin\'eaire respectant la structure
euclidienne infinit\'esimale de $T_p S$:
\[
\big\langle
J^\lambda(X),\,
J^\lambda(Y)
\big\rangle
\,=\,
\langle
X,\,Y
\rangle.
\]

On appelle {\sl groupe orthogonal} le groupe des transformations
qui respectent ce produit scalaire, et on le note:
\[
J^\lambda
\,\in\,
{\sf O}
\big(T_pS\big)
\,\,\cong\,\,
{\sf O}_2(\R)
\,\,:=\,\,
\big\{
A
\in
{\sf GL}_2(\R)
\colon\,
{}^\TT\!AA
=
I_{2\times 2}
\big\}.
\]
Il est \'evident que pour tous $\lambda_1, \lambda_2 \in
\Omega_p(S)$, on a:
\[
J^{\lambda_1}
\circ
J^{\lambda_2}
\,=\,
J^{\lambda_1\cdot\lambda_2},
\]
ainsi que:
\[
J^{\lambda^-}
\,=\,
\big(
J^\lambda
\big)^{-1},
\]
o\`u $\lambda^-$ est le lacet $\lambda$ parcouru en sens inverse.

Dans le groupe de Poincar\'e, certains lacets $\lambda$ sont dits
{\sl homotopes \`a z\'ero} lorsqu'ils peuvent \^etre
(diff\'erentiablement) d\'eform\'es jusqu'\`a devenir
{\em constants}\,\,---\,\,autrement dit, 
s'il existe une famille \`a $1$ param\`etre:
\reqnomode\usetagform{EngelLie}
\begin{align}
\big(
\lambda_s(t)
\big)_{s\in[0,1]}
\ \ \ \ \ \ \ \ \ \ \
\text{\rm avec}
\ \ \ \ \
\lambda_0(t)
&
\,\equiv\,
\lambda(t),
\notag
\\
\lambda_1(t)
&
\,\equiv\,
p
\tag{(\text{\rm lacet constant}).}
\end{align}

On v\'erifie alors que les $J^{\lambda^{\rm o}}$, pour des lacets
$\lambda^{\rm o}$ homotopes \`a z\'ero, pr\'eservent
toute orientation donn\'ee de $T_pS$ donn\'ee \`a l'avance,
c'est-\`a-dire: 
\[
J^{\lambda^{\rm o}} 
\,\in\,
{\sf SO}\big( T_p S\big)
\,\,\cong\,\,
{\sf SO}_2(\R)
\,\,:=\,\,
\big\{
A
\in
{\sf GL}_2(\R)
\colon\,
{}^\TT\!AA
=
I_{2\times 2},\,\,
\det\,A
=
1
\big\},
\]
ce qui conduit \`a introduire:
\[
\Omega_p^{\rm o}(S)
\,:=\,
\Big\{
\text{\rm lacets diff\'erentiables}\,\,
\lambda\colon\,[0,1]
\longrightarrow
S,\,\,
\text{\rm homotopes \`a z\'ero},\,\,
\lambda(0)=p=\lambda(1)
\Big\}.
\]
On appelle:
\[
\aligned
\big\{J^\lambda\big\}_{\lambda\in\Omega_p(S)}
\ \ \ \ \
&
\text{\sl groupe d'holonomie de}\,\,S\,\,\text{\rm en}\,\,p,
\\
\big\{J^\lambda\big\}_{\lambda\in\Omega_p^{\rm o}(S)}
\ \ \ \ \
&
\text{\sl groupe d'holonomie restreinte de}\,\,S\,\,\text{\rm en}\,\,p.
\endaligned
\]
Ces groupes d'holonomie sont intimement li\'es \`a une nature
quantitative qui est la {\em courbure} de la surface $S$
munie de sa m\'etrique.

Le commentaire philosophique de cette situation se d\'eveloppe sur
plusieurs plans. Dans les math\'ematiques r\'ecentes, la g\'eom\'etrie
est contr\^ol\'ee par une structure alg\'ebrique qui s'impose: la
structure de groupe. Par sa forme, un groupe r\'eussit \`a exprimer et
à contr\^oler un ensemble de gestes exploratoires d'une entit\'e
g\'eom\'etrique comme une surface. Cela est d\^u partiellement \`a la
forme de stabilit\'e qu'introduit la structure de groupe r\'egissant
les \'el\'ements g\'eom\'etriques qui sont en jeu. Les lacets, qui
sont des chemins dont les extr\'emit\'es sont confondues, 
permettent de pr\'elever juste ce qu'il faut
d'occupation de l'espace (surface) consid\'er\'e pour l'explorer. Une
sorte de <<\,{\sl squelette flexible}\,>>
peut rendre compte de propri\'et\'es
topologiques de la surface (ou de la vari\'et\'e). Et ces \'el\'ements,
en se composant suivant les lois du groupe qu'ils forment, d\'eploient
l'observation topologique.

Nous avons affaire ici \`a une forme d'abstraction 
<<\,interm\'ediaire\,>> 
dans laquelle un aspect de l'espace (surface ou vari\'et\'e) est
pr\'elev\'e (lacet), lui-m\^eme command\'e par 
des \'el\'ements de structure
alg\'ebrique.
C'est cette situation <<\,interm\'ediaire\,>> qui fait la force
conceptuelle du groupe de Poincar\'e. Nous sommes au plus pr\`es de
l'objet topologique. On remarquera que la propri\'et\'e topologique
principale \`a laquelle on fait souvent r\'ef\'erence est la
connexit\'e, qui nous indique si un espace est tout d'un bloc ou
non. Les \'el\'ements de base qui nous permettent d'examiner l'espace
sont suffisamment autonomis\'es pour que nous puissions exercer notre
contr\^ole et raisonner sur eux en usant et de l'alg\`ebre et de la
g\'eom\'etrie qu'ils comportent.

En construisant \`a partir des lacets le groupe d'holonomie, un
observateur se d\'eplace dans l'espace tangent \`a la vari\'et\'e, en
pr\'eservant l'orientation. Nous b\'en\'eficions des avantages de la
structure de groupe, mais cette fois en nous rapprochant de la
possibilit\'e d'une analyse de la courbure de la vari\'et\'e. Ce qu'il
faut comprendre tient au processus de la r\'eflexion math\'ematique:
nous voulons toujours comprendre ce qu'est la courbure, mais nous nous
sommes plac\'es \`a un niveau depuis lequel les concepts de
l'alg\`ebre et de la topologie peuvent contribuer \`a d\'eterminer le
concept de courbure. Et il faut comprendre comment ces nouvelles
d\'eterminations entrent en sc\`ene.

Il nous faut insister sur le fait qu'une fois la synth\`ese
d\'etermin\'ee effectu\'ee, les concepts ainsi imbriqu\'es continuent
d'agir ensemble, ici l'alg\`ebre et la topologie qui lui est
int\'egr\'ee. C'est la raison qui fait de cette synth\`ese une forme
d'irr\'eversibilit\'e.

Les deux concepts dont nous usons ici sont ceux de groupe d'holonomie
et de sa restriction. Consid\'erons cette dernière. Nous manipulons des
lacets homotopes \`a z\'ero que nous transportons parall\`element. Si
bien que nous synth\'etisons la th\'eorie du transport parall\`ele et
la th\'eorie du groupe de Poincar\'e. Le groupe d'holonomie qui
r\'eussit cette synth\`ese est devenu le nouvel analyseur de la
courbure d'une vari\'et\'e. Le groupe d'holonomie est un groupe
orthogonal, car il doit pr\'eserver le produit scalaire, la composition
des lacets correspondant \`a leur produit port\'e par le transport
parall\`ele.

\Section{\bf La courbure comme infinit\'esimalisation 
de l'holonomie}
\label{courbure-infinitesimalisation-holonomie}
\HEAD{\ref{courbure-infinitesimalisation-holonomie}.~{\sf La courbure 
comme infinit\'esimalisation de l'holonomie
}
}{
Jo\"el {\sc Merker} et Jean-Jacques {\sc Szczeciniarz}}

\begin{center}
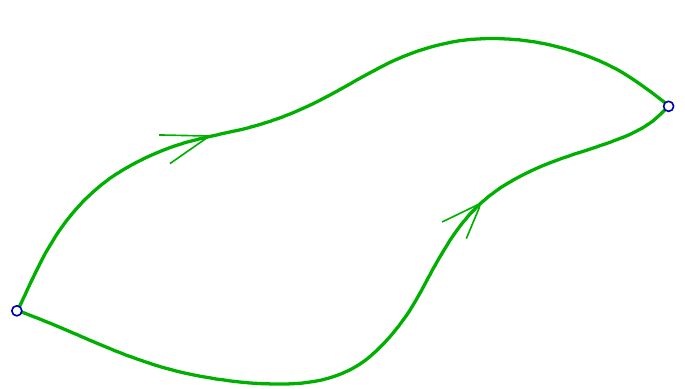
\end{center}

En g\'en\'eral, \'etant donn\'e deux points distincts quelconques sur une
surface:
\[
p,\,q\,\in\,S,
\ \ \ \ \ \ \ \ \ \ \ \ \ \ \ \ \ \
p
\,\neq\,
q,
\]
\'eventuellement arbitrairement proches, quand on prend
deux courbes diff\'erentes $\gamma$ et $\delta$ 
allant de $p$ \`a $q$, les transports parall\`eles associ\'es
{\em diff\`erent:}
\[
J_{q,p}^\gamma
\,\neq\,
J_{q,p}^\delta.
\]
Comment, alors, mesurer l'\'ecart?
Comment quantifier cette diff\'erence?

Une difficult\'e r\'eside dans la variabilit\'e infinie des
choix de courbes joignant $p$ \`a $q$: l'espace
des courbes est r\'eellement de dimension infinie!

Pour comprendre ce d\'efaut de co\"{\i}ncidence, \`a nouveau,
la th\'eorie math\'ematique va {\em infinit\'esimaliser les
consid\'erations}, et contre toute
attente, rien ne va \^etre perdu
\`a travers ce passage au microscopique.
Ainsi, le point $q$ sera suppos\'e
(tr\`es) proche de $p$. Ensuite, les courbes $\gamma$ et
$\delta$ vont, pour simplifier, \^etre envisag\'ees
comme c\^ot\'es de parall\'elogrammes (infiniment) petits,
afin de {\em normaliser} la forme de la d\'ependance
en le choix des courbes.

\begin{center}
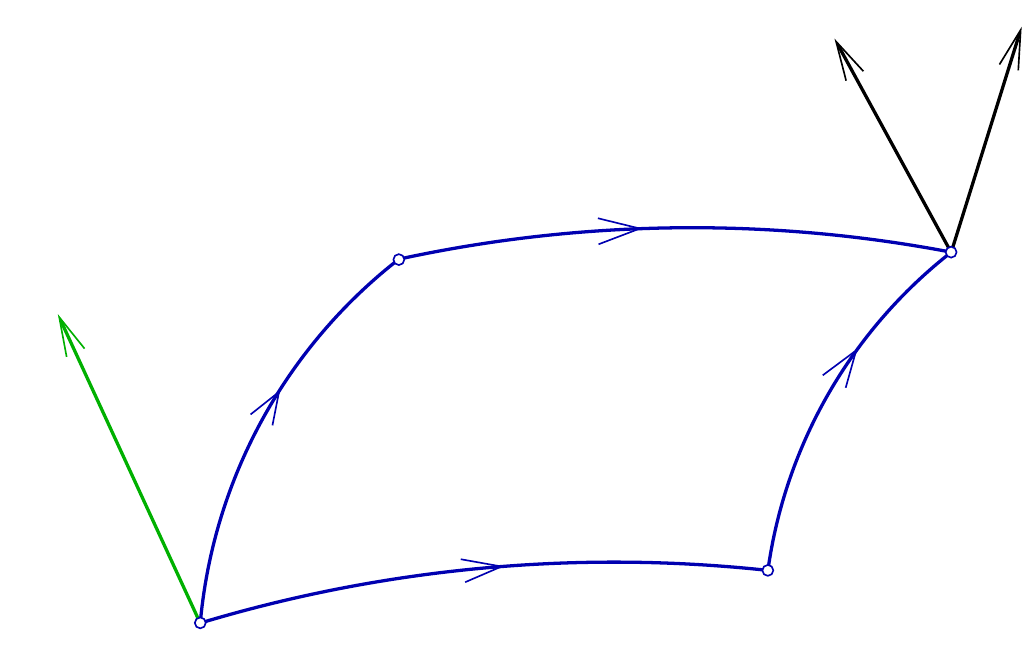
\end{center}

Plus pr\'ecis\'ement, dans des coordonn\'ees locales:
\[
(u,v)
\,\in\,
U
\,\subset\,
S
\]
sur un ouvert $U \subset S$, avec deux quantit\'es r\'eelles
$s, t \in \R$ petites qui tendront vers $0$, introduisons,
comme sur le diagramme, $4$ points qui sont sommets
d'un parall\'elogramme.

Si donc $Z_p = Z_{u,v}$ est un vecteur quelconque dans l'espace
tangent $T_pS$, et si $\gamma$ est la courbe inf\'erieure droite
du parall\'elogramme:
\[
\gamma
\colon\ \ \ 
(u,v)
\,\,\longrightarrow\,\,
\big(u+s,\,v\big)
\,\,\longrightarrow\,\,
\big(u+s,\,v+t\big),
\]
tandis que $\delta$ est la courbe sup\'erieure gauche:
\[
\delta
\colon\ \ \ 
(u,v)
\,\,\longrightarrow\,\,
\big(u,\,v+t\big)
\,\,\longrightarrow\,\,
\big(u+s,\,v+t\big),
\]
on n'a en g\'en\'eral pas co\"{\i}ncidence des transports le long
de $\gamma$ et le long de $\delta$:
\[
J^\gamma\big(Z_p\big)
\,\neq\,
J^\delta\big(Z_p\big)
\eqno
{\scriptstyle{(\text{\rm en g\'en\'eral})}}.
\]

En fait, comme le point $q = (u+s,\, v+t)$ va varier avec
$(s,t) \longrightarrow (0,0)$ lorsque le parall\'elogramme va
se r\'eduire, il vaut mieux\,\,---\,\,car cela
est essentiellement \'equivalent\,\,---\,\,mesurer
un \'ecart de transport de vecteurs en $q$ vers
l'espace tangent {\em fixe} $T_p S$ en $p$.

On consid\`ere donc un champ de vecteurs $Z \in \mathcal{X}(U)$
d\'efini au voisinage de $p$, et on effectue les deux 
transports parall\`eles de la valeur $Z_q = Z_{u+s, v+t}$ de $Z$
en $q$ le long des deux courbes $\gamma^-$ et
$\delta^-$ parcourues {\em en sens inverse:}
\[
J^{\gamma^-}\big(Z_q\big)
\,\neq\,
J^{\delta^-}\big(Z_q\big)
\eqno
{\scriptstyle{(\text{\rm en g\'en\'eral})}}.
\]

\begin{center}
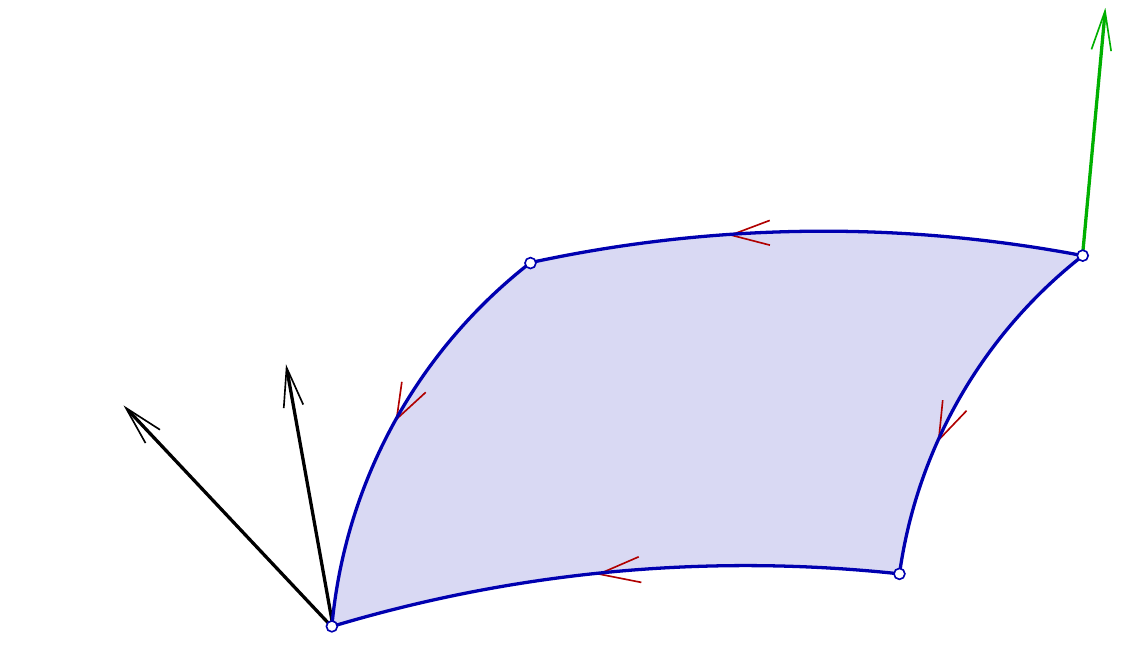
\end{center}

Par divers arguments heuristiques test\'es sur des exemples
o\`u des calculs complets sont accessibles, on se convainc
que cette diff\'erence est quantitativement de l'ordre
de l'{\em aire} du parall\'elogramme en question, laquelle
est approximativement \'egale \`a une constante que multiplie
le produit:
\[
s\,t.
\]
Comme on l'a d\'ej\`a fait pour d\'efinir la diff\'erentiation
covariante:
\[
\aligned
\nabla_{X_0}(Y)
\,:=\,
&\,
\frac{d}{dt}
\Big\vert_{t=t_0}
J_{t_0,t}^\gamma
\big(
Y_{\gamma(t)}
\big)
\\
\,=\,
&\,
\underset{t\to t_0}{\lim}\,\,
\frac{J_{t_0,t}^\gamma\big(Y_{\gamma(t)}\big)
-
Y_{\gamma(t_0)}}{t-t_0},
\endaligned
\]
on va ici examiner le comportement du {\em quotient:}
\[
\underset{(s,t)\to(0,0)}{\lim}\,\,
\frac{J^{\gamma^-}(Z)
-
J^{\delta^-}(Z)}{s\,\,t},
\]
qui converge toujours vers une quantit\'e finie, comme le d\'evoile
un th\'eor\`eme fondamental qui suit. Dans cet \'enonc\'e, on observera
que:

\smallskip\noindent$\square$\,\,
$\gamma$ suit $\frac{\partial}{\partial u}$ pendant le temps
$s$, puis $\frac{\partial}{\partial v}$ pendant le temps $t$;

\smallskip\noindent$\square$\,\,
$\delta$ suit $\frac{\partial}{\partial v}$ pendant le temps
$t$, puis $\frac{\partial}{\partial u}$ pendant le temps $s$.

\begin{Theoreme}
Pour tout champ de vecteurs $Z \in \mathcal{X}(U)$:
\[
\underset{(s,t)\to(0,0)}{\lim}\,\,
\frac{J^{\gamma^-}(Z)
-
J^{\delta^-}(Z)}{s\,\,t}
\,\,=\,\,
\nabla_{\!\!\frac{\partial}{\partial u}}
\Big(
\nabla_{\!\!\frac{\partial}{\partial v}}
(Z)
\Big)
-
\nabla_{\!\!\frac{\partial}{\partial v}}
\Big(
\nabla_{\!\!\frac{\partial}{\partial u}}
(Z)
\Big).
\eqno\qed
\]
\end{Theoreme}

Ainsi, la quantification du d\'efaut de commutativit\'e
infinit\'esimale du transport parall\`ele se r\'e-exprime-t-elle comme
d\'efaut de commutativit\'e des deux op\'erateurs de diff\'erentiation
covariante $\nabla_{\!\!\frac{\partial}{\partial u}}$ et
$\nabla_{\!\!\frac{\partial}{\partial v}}$ dirig\'es par les c\^ot\'es
du parall\'elogramme infinit\'esimal.

Ce r\'esultat remarquable manifeste un mouvement th\'eorique
fr\'equent en math\'ematique: le passage d'une propri\'et\'e d'un
domaine d'objets \`a une propri\'et\'e du m\^eme type des
op\'erateurs.

Cette expression peut ensuite \^etre int\'egr\'ee pour retrouver
des \'ecarts globaux dans les transports parall\`eles,
comme nous allons le voir dans peu de temps.
Pour l'instant, ajoutons que cette \'equation-limite
devrait se g\'en\'eraliser \`a des paires de champs de vecteurs
quelconques:
\[
X,\,Y\,\in\,
\mathcal{X}\big(U\subset S\big)
\ \ \ \ \ \ \ \ \ \ \ \ \ \ \ \ \ \
\text{\rm au lieu de}\ \
\frac{\partial}{\partial u},\,\,
\frac{\partial}{\partial v},
\]
mais une diff\'erence de taille appara\^it. En effet, pour deux
champs quelconques:
\[
\aligned
X
&
\,=\,
A(u,v)\,
\frac{\partial}{\partial u}
+
B(u,v)\,
\frac{\partial}{\partial v},
\\
Y
&
\,=\,
C(u,v)\,
\frac{\partial}{\partial u}
+
D(u,v)\,
\frac{\partial}{\partial v},
\endaligned
\]
on n'a pas en g\'en\'eral 
co\"{\i}ncidence des deux points-extr\'emit\'es
de courbes int\'egrales:
\[
\exp\,
\Big(
tY
\big(
\exp\,(sX)
\big)(p)
\Big)
\,\,\neq\,\,
\exp\,
\Big(
sX
\big(
\exp\,(tY)
\big)(p)
\Big),
\]
ce qui \'etait le cas avec $X = \frac{\partial}{\partial u}$ et
$Y = \frac{\partial}{\partial v}$, car:
\[
\big(
(u,v)
+
(s,0)
\big)
+
(0,t)
\,\,=\,\,
\big(
(u,v)
+
(0,t)
\big)
+
(s,0).
\]
Donc puisque les points-extr\'emit\'es de $\gamma$ et de $\delta$
ne co\"{\i}ncident pas en g\'en\'eral, on ne peut pas vraiment parler
du quotient:
\[
\frac{J^{\gamma^-}(Z)-J^{\delta^-}(Z)}{s\,\,t}
\overset{\text{\bf ?}}{\,\,=\,\,}
\nabla_X
\big(
\nabla_Y(Z)
\big)
-
\nabla_Y
\big(
\nabla_X(Z)
\big).
\]
En fait, la non-co\"{\i}ncidence des deux extr\'emit\'es ci-dessus
est quantifi\'ee par le {\sl crochet de Lie} entre
les deux champs de vecteurs:
\[
\big[X,\,Y\big]
\,:=\,
\big(
X(C)
-
Y(A)
\big)\,
\frac{\partial}{\partial u}
+
\big(
X(D)
-
Y(B)
\big)\,
\frac{\partial}{\partial v},
\]
au sens o\`u en un point $p \in U \subset S$, la valeur
de ce crochet intervient comme
terme d'ordre $2$ dans le {\em commutateur}
entre les flots respectifs:
\[
\exp\,(-tY)
\circ
\exp\,(-tX)
\circ
\exp\,(tY)
\circ
\exp\,(tX)
(p)
\,\,=\,\,
p
+
t^2\,
\big[X,Y\big](p)
+
{\rm O}(t^3),
\]
commutateur qui ne peut donc \^etre identiquement nul
que lorsque le crochet de Lie s'annule identiquement.

\begin{center}
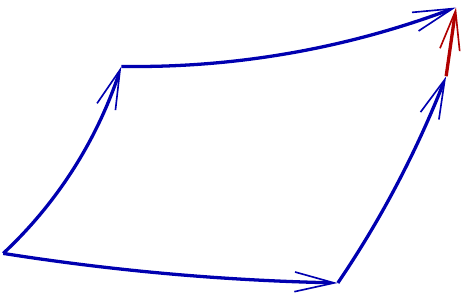
\end{center}

\Section{\bf Courbure associ\'ee \`a une diff\'erentiation 
covariante}
\label{courbure-differentiation-covariante}
\HEAD{\ref{courbure-differentiation-covariante}.~{\sf Courbure 
associ\'ee \`a une diff\'erentiation covariante}
}{
Jo\"el {\sc Merker} et Jean-Jacques {\sc Szczeciniarz}}

Ceci conduit n\'ecessairement \`a tenir compte du crochet de Lie dans le
quotient ci-dessus en ajoutant un terme correctif.
Le bon choix s'exprime alors comme un nouveau concept:
\[
R(X,Y)\,Z
\,:=\,
\nabla_X\big(\nabla_Y(Z)\big)
-
\nabla_Y\big(\nabla_X(Z)\big)
-
\nabla_{[X,Y]}\,Z.
\]

\begin{center}
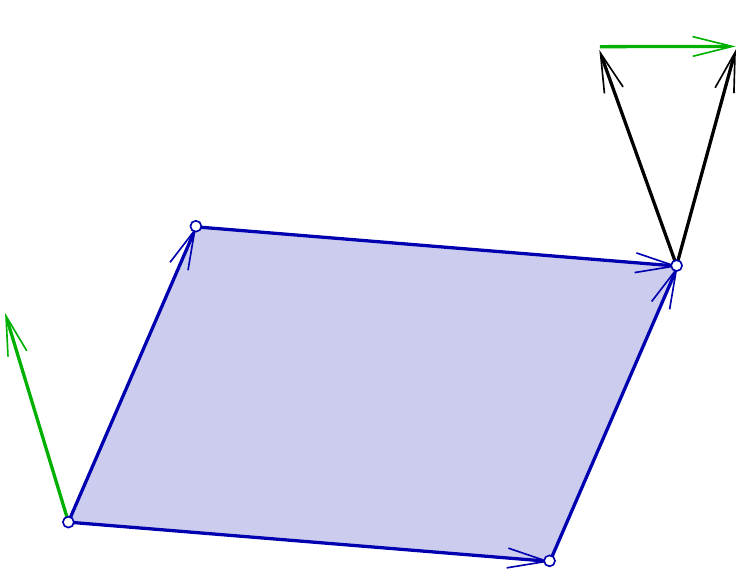
\end{center}

Un diagramme, dessin\'e dans le cas o\`u $[X,Y] = 0$ commutent,
illustre ce que fait cet op\'erateur $R(X,Y)\, \centersmallbullet$
d\'efini pour une 
paire de champs $(X,Y)$ qui d\'ecoupent un parall\'elogramme:
il transforme tout champ de vecteurs $Z$ en 
le champ de vecteurs $R(X,Y)\, Z$ obtenu en 
calculant la diff\'erence des valeurs de ses deux transports
parall\`eles le long des deux c\^ot\'es du parall\'elogramme.

\begin{Theoreme}
L'op\'erateur de courbure $R(\centersmallbullet,
\centersmallbullet) \centersmallbullet$ est un {\sl tenseur}
au sens o\`u l'application:
\[
\big(X,Y,Z\big)
\,\,\longmapsto\,\,
R(X,Y)\,Z,
\]
antisym\'etrique par rapport \`a ses deux premiers arguments, est
trilin\'eaire, y compris pour la multiplication par des fonctions
diff\'erentiables $f$, $g$, $h$ quelconques sur la surface:
\[
R\big(fX,\,gY\big)\,hZ
\,=\,
fgh\,R(X,Y)\,Z.
\]
Vis-\`a-vis de la m\'etrique, il satisfait de plus: 
\[
\big\langle
R(X,Y)Z,\,T
\big\rangle
\,\,=\,\,
-\,
\big\langle
R(X,Y)T,\,Z
\big\rangle.
\eqno\qed
\]
\end{Theoreme}

Les calculs pass\'es sous silence
ici montrent qu'il y a n\'ecessit\'e du terme correctif
$-\nabla_{[X,Y]}Z$ pour garantir ce caract\`ere tensoriel.  Une double
non-commutativit\'e est donc \`a l'{\oe}uvre ici: celle du crochet de
Lie, et celle de la non-holonomie.

En conclusion, on a infinit\'esimalis\'e l'holonomie, envisag\'ee
comme \'ecart de commutativit\'e, et on l'a ensuite conceptualis\'ee
au moyen d'un op\'erateur simple qui d\'epend de parall\'elogrammes
associ\'es \`a des paires de champs de vecteurs. En pla\c{c}ant les
uns \`a c\^ot\'es des autres de nombreux parall\'elogrammes
infinit\'esimaux qui reconstitueraient des domaines macroscopiques de
la surface $S$, on pourra s'imaginer reconstituer une vision plus
globale de l'holonomie.

Le commentaire philosophique porte ici sur la non-commutativit\'e. De
fa\c{c}on g\'en\'erale, il est facile de comprendre comment la 
non-commutativit\'e est comprise comme un \'ecart \`a l'holonomie et
pourquoi il en est ainsi. Lorsque deux vecteurs sont transport\'es
parall\`element le long de deux c\^ot\'es suivant deux parcours
diff\'erents d'un
parall\`elogramme infinit\'esimal trac\'e sur une surface, ces deux
vecteurs ne sont pas les m\^emes en direction. On se repr\'esente
facilement pourquoi cette diff\'erence mesure la courbure de la
surface. L'\'ecart en direction est exprim\'ee par l'holonomie. Le
d\'eveloppement conceptuel \'etend l'explication \`a un domaine plus
g\'en\'eral qui est celui des rotations. Cette nouvelle synth\`ese
explique et d\'eveloppe en m\^eme temps le concept de courbure. Quand
l'espace est courbe il n'est nullement indiff\'erent de se d\'eplacer
suivant des parcours (infinit\'esimaux) pour aller d'un point \`a un
autre, et il devient possible de concevoir et de mesurer la
diff\'erence. Du m\^eme coup, j'ai concr\'etis\'e les d\'eterminations
de la courbure, comme modalit\'e de la spatialisation de la quantit\'e
ou quantification de l'espace du d\'eplacement.

Il est tr\`es remarquable que tous les arguments qui cherchent \`a
montrer quelle est la structure explicite de la courbure reposent sur
des \' enonc\'es portant eux-m\^emes sur l'examen de la {\em
commutativit\'e} d'un parcours de trajectoires. On a vu plus haut
comment cette r\'eflexion s'effectue \`a l'aide de la structure de
groupe de Lie. Celle-ci est une autre synth\`ese qui nous introduit
\`a un autre champ de r\'eflexion. Nous en dirons un minimum pour
\^etre \`a m\^eme d'\'enoncer un commentaire philosophique.

Sur une vari\'et\'e $M$, on d\'efinit une d\'eriv\'ee de Lie par
rapport à un champ de vecteurs $\xi \in \mathcal{X}(M)$.  La
d\'eriv\'ee de Lie $\mathcal{L}_\xi (\centersmallbullet)$ est
d\'efinie comme un op\'eration effectu\'ee par rapport \`a ce champ de
vecteurs $\xi$. Appliqu\'ee \`a une grandeur géométrique $Q$, elle
mesure les variations de cette grandeur, par comparaison \`a ce
qu'elle serait si elle \'etait simplement entra\^in\'ee par le champ
$\xi$. C'est une mesure de variation relative \`a un champ de
vecteurs. On peut la noter $\mathcal{L}_\xi (Q)$. Elle s'applique dans
des contextes tr\`es g\'en\'eraux. Consid\'erons la d\'eriv\'ee de Lie
d'un champ de vecteurs $\eta$ par rapport au champ $\xi$.

On a vu qu'un champ de vecteurs s'interpr\`ete comme un op\'erateur
diff\'erentiel agissant sur des champ scalaires et satisfaisant aux
lois d'un op\'erateur diff\'erentiel. Il en est de m\^eme de
l'op\'erateur
\[
L(\Phi)
\,:=\,
\xi\big(\eta(\Phi)\big)
- 
\eta\big(\xi(\Phi)\big),
\]
qui satisfait \`a ces lois si $\xi$ et $\eta$ y satisfont aussi. Ce {\sl
commutateur} des deux op\'erations $\xi$ et $\eta$ s'\'ecrit donc
sous la forme d'un crochet de Lie:
\[
L
\,=\,
\xi\circ\eta
-
\eta\circ\xi
\,:=\,
\big[\xi,\eta\big],
\]
et nous l'avons utilis\'e ci-dessus.

De nouveau, ce commutateur a une signification g\'eom\'etrique. Nous
recourons de nouveau \`a un parall\`elogramme de fl\`eches construit
\`a l'aide de $\xi$ et de $\eta$.  Les c\^ot\'es sont construits \`a
l'ordre ${\rm O}(\varepsilon)$, et $L$ mesure le
d\'efaut de commutativit\'e \`a l'ordre ${\rm O}(\varepsilon^2)$. Nous
travaillons, comme nous l'avons indiqué, pour aller au-del\`a de la
simple structure de l'espace tangent, dans une exploration {\em \`a
l'ordre $2$}. Un d\'eplacement infinit\'esimal d'un \'el\'ement du
groupe le long du champ de vecteurs $\xi$ repr\'esente la
transformation correspondant \`a la multiplication \`a gauche de tout
\'el\'ement du groupe par l'\'el\'ement infinit\'esimal repr\'esent\'e
par $\xi$. Si on multiplie par $\varepsilon$, on peut concevoir
$\varepsilon\, \xi$ comme un d\'eplacement \`a l'ordre ${\rm
O}(\varepsilon)$ le long du champ de vecteurs $\xi$. Pour obtenir
une information sur la structure de groupe, correspondant \`a la
structure additive de l'espace tangent \`a l'\'el\'ement neutre du
groupe, on doit aller comme on l'a fait ci-dessus jusqu'\`a l'ordre
${\rm O} (\varepsilon^2)$, en regardant le commutateur qui est {\em
un essai dans l'abstrait de la commutativit\'e}. Et $\varepsilon^2
[\xi, \eta]$ correspond \`a la fermeture du parall\'elogramme dont les
premiers c\^ot\'es sont $\varepsilon\, \xi$ et $\varepsilon\, \eta$ \`a
l'origine. La plong\'ee dans l'infinit\'esimal s'est faite \`a des
ordres diff\'erents qui ont permis de d\'ecrire les propri\'et\'es de
commutativit\'e.

Si l'on se rapporte aux d\'eveloppements ci-dessus, on comprend la
relation de cette op\'eration de commutation avec les \'el\'ements
infinit\'esimaux d'un groupe de Lie. L'exploration d'une surface ou
d'une vari\'et\'e a consist\'e au sein d'une infinit\'esimalisation de
base \`a tester la commutativit\'e de trajectoires (de courbes) sur la
surface en question. Comme les trajectoires doivent pouvoir se
composer et se d\'ecomposer elles ont \'et\'e elles-m\^emes
ins\'er\'ees dans une structure de groupe, elles-m\^emes supports
pour des propri\'et\'es topologiques ou diff\'erentielles. Ces
constructions dans lesquelles objets et fonctions \'echangent leur
place approfondissent \`a chaque \'etape les formes de r\'eflexivit\'e
synth\'etique, installent {\em nos} productions comme des objets
d'analyse {\em en soi}.

\Section{\bf Invariances de la forme et de la fonction de courbure}
\label{invariance-courbure}
\HEAD{\ref{invariance-courbure}.~{\sf Invariances 
de la forme et de la fonction de courbure
}
}{
Jo\"el {\sc Merker} et Jean-Jacques {\sc Szczeciniarz}}

\`A pr\'esent, comment s'exprime l'action de cet op\'erateur 
$R(\centersmallbullet, \centersmallbullet)\, 
\centersmallbullet$ dans un
rep\`ere orthonorm\'e local $\{ X, Y\}$ sur un ouvert $U \subset S$? Un
calcul qui utilise ce qui pr\'ec\`ede conduit aux expressions:
\[
\aligned
R(Z,T)\,X
&
\,=\,
-\,
d\omega(Z,T)\,Y,
\\
R(Z,T)\,Y
&
\,=\,
d\omega(Z,T)\,X,
\endaligned
\]
o\`u $d\omega = d\omega_{X,Y}$ est la diff\'erentielle
ext\'erieure de la $1$-forme de connexion 
$\omega = \omega_{X,Y}$ associ\'ee\,\,---\,\,une 
$2$-forme diff\'erentielle absolument fondamentale!
Car la cl\'e de vo\^ute de toute la th\'eorie est
{\em l'invariance absolue de cette $2$-forme diff\'erentielle}. 

En effet, on d\'emontre que pour tout autre choix $\{ X', Y'\}$
de rep\`ere orthorm\'e sur l'ouvert $U \subset S$, qui s'exprime
alors n\'ecessairement sous la forme:
\[
\aligned
X_p'
&
\,=\,
\cos\,\theta(p)\,X_p
+
\sin\,\theta(p)\,Y_p,
\\
Y_p'
&
\,=\,
-\,\sin\,\theta(p)\,X_p
+
\cos\,\theta(p)\,Y_p,
\endaligned
\]
au moyen d'une certaine fonction diff\'erentiable unique
$p \longmapsto \theta(p)$ d\'efinie sur $U$, 
la nouvelle $1$-forme diff\'erentielle de connexion associ\'ee 
$\omega_{X',Y'}$ est reli\'ee \`a l'ancienne par la formule simple:
\[
\omega_{X',Y'}
\,=\,
\omega_{X,Y}
-
d\theta,
\]
et alors la relation de Poincar\'e $d \circ d = 0$ montre l'\'egalit\'e:
\[
\aligned
d\omega_{X',Y'}
&
\,=\,
d\omega_{X,Y}
-
d\circ d\theta
\\
&
\,=\,
d\omega_{X,Y}
-
0.
\endaligned
\]

En conclusion, la $2$-forme diff\'erentielle, dite
{\sl $2$-forme de courbure:}
\[
\Omega
\,:=\,
d\omega,
\]
est parfaitement bien d\'efinie, de mani\`ere interne \`a la surface,
ind\'ependamment de tout choix de rep\`ere orthonorm\'e local
$\{X,Y\}$, et par cons\'equent, elle est m\^eme d\'efinie
{\em globalement} sur $S$. 

Cette simple relation $\omega' = \omega - d\theta$ qui fait
appara\^itre en l'esp\`ece du terme additionnel $-d\theta$ une
correction n\'ecessaire, ou {\sl changement de jauge}, lorsqu'on
effectue un changement de r\'ef\'erentiel, produit une relation
invariante:
\[
\Omega'
\,=\,
\Omega,
\]
gr\^ace \`a un {\em proc\'ed\'e d'\'elimination} universel 
en alg\`ebre et en g\'eom\'etrie: l'{\sl effectuation suppressive}. 
Les objets d\'ependent de certains choix? 
Cherchons alors \`a faire dispara\^itre ces choix par 
approfondissement du calcul\,\,---\,\,ici, diff\'erentiation.
Ce proc\'ed\'e est un {\sl leitmotiv} des math\'ematiques
dans leur ensemble, et c'est lui qui, par le calcul,
engendre des conceptualisations. L'{\oe}uvre
d'\'Elie Cartan en t\'emoigne.

Ensuite, faisons agir cette $2$-forme diff\'erentielle
$\Omega$ sur la paire de r\'ef\'erence:
\[
k
\,:=\,
\Omega(X,Y)
\,=\,
d\omega(X,Y),
\]
ce qui produit une {\em fonction globale} sur la surface,
appel\'ee {\sl courbure de Gauss}, en tant qu'une v\'erification
th\'eorique convainc qu'elle co\"{\i}ncide bien avec le
concept dont Gauss a d\'emontr\'e l'invariance d'une mani\`ere
totalement diff\'erente. 

Qui plus est, le raisonnement \'el\'ementaire qui pr\'ec\`ede 
permet de d\'eduire une version g\'en\'erale du {\sl Theorema
Egregium} de Gauss, \`a savoir: Si une \'equivalence $\Phi$
entre deux ouverts de deux surfaces
munies de deux m\'etriques est donn\'ee:
\[
S
\,\supset\,
U
\overset{\Phi}{\,\,\longrightarrow\,\,}
\overline{U}
\,\subset\,
\overline{S},
\] 
et si $\big\{ \overline{X}, \overline{Y} \big\}$ est le
rep\`ere orthonorm\'e transf\'er\'e par $\Phi$ sur la surface-image:
\[
\overline{X}
\,=\,
\Phi_\ast(X)
\ \ \ \ \ \ \ \ \ \ \ \ \ \ \ \ \ \
\text{et}
\ \ \ \ \ \ \ \ \ \ \ \ \ \ \ \ \ \
\overline{Y}
\,=\,
\Phi_\ast(Y),
\]
alors la $2$-forme de courbure et la fonction de courbure
co\"{\i}ncident aux points correspondants par $\Phi$:
\[
\Omega
\,=\,
\overline{\Omega}
\ \ \ \ \ \ \ \ \ \ \ \ \ \ \ \ \ \
\text{et}
\ \ \ \ \ \ \ \ \ \ \ \ \ \ \ \ \ \
k
\,=\,
\overline{k}.
\]

Nous avons fait appara\^itre un op\'erateur, une nouvelle forme
synth\'etique a abstrait le tenseur de courbure. Cet op\'erateur
appara\^it dans le calcul des d\'erivations covariantes et exprime le
fait que l'holonomie restreinte emp\^eche
d'avoir un th\'eor\`eme d'interversion des d\'eriv\'ees secondes,
lorsque nous comparons la d\'efinition du crochet avec l'expression
suivante des d\'eriv\'ees covariantes:
\[
R(X, Y) 
\cdot 
Z
\,=\,
\nabla_{X}\nabla_{Y}Z
-
\nabla_{Y}\nabla_{X}Z
-
\nabla_{[X, Y]}Z,
\]
d\'efinie pour tout triplet $X, Y, Z$ de champs de
vecteurs. Cette expression est $C^{\infty}$-trilin\'eaire en
$X$, $Y$, $Z$, et antisym\'etrique en $X, Y$. Cette
trilin\'earit\'e est encore une mani\`ere dont se 
<<\,d\'eploie\,>> la
courbure en fonction des champs de vecteurs lorsqu'ils sont
multipli\'es par une fonction et ce, pour chacun d'eux. 
La multiplication se transporte d'un
champ \`a toute l'expression, de m\^eme que l'expression est additive
en chacun des trois champs. La courbure est devenue un concept qui
permet de manipuler trois champs de vecteurs en respectant les
contraintes ci-dessus. Nous y avons gagn\'e une modalit\'e de
r\'eflexion de l'espace.

En particulier nous avons not\'e qu'une d\'etermination r\'eciproque
entre l'holonomie restreinte et la forme de courbure est mise en
\'evidence.  Nous avons \'egalement, en effet, pu donner une
expression abstraite de la forme de courbure. Nous disposons d'abord
d'un op\'erateur de courbure. Cet op\'erateur est une nouvelle
synth\`ese, mais en tant qu'op\'erateur. Il est unique. S'il est
d\'efini sur un certain type d'espace local, il ne
d\'epend pas de l'orientation qui y est choisie. Cette remarque pour
insister sur le fait que des propri\'et\'es essentielles de l'espace
(des ouverts de cet espace) sont incluses dans l'op\'erateur
lui-m\^eme. 

Que signifie ce fait th\'eorique que la courbure dont nous
avons pu voir les expressions soit concentr\'ee et exprim\'ee dans un
op\'erateur? Cela signifie avant tout qu'un objet, ou un ph\'enom\`ene
g\'eom\'etrique appr\'ehendable dans une repr\'esentation
caract\'erisant des modes de pr\'esentation (des reliefs, des
creusements, des d\'eformations du planaire,
{\em etc.}) se traduit dans un
mode de connaissance conceptuelle et calculable. Ce mode de
connaissance fortement synth\'etique est devenue une forme (qui
int\`egre \`a la r\'eflexion qu'il propose les avantages fournis par
les formes diff\'erentielles), pourvue d'une tr\`es grande capacit\'e
d'int\'egration. Nous pouvons y voir les conditions de possibilit\'e
{\em a priori} tout autant intuitives que conceptuelles de la
quantification de l'espace ou, comme dit ci-dessus, de la spatialisation
de la quantification.

\Section{\bf R\'eint\'egration de l'infinit\'esimalis\'e: 
mesure de courbure}
\label{reintegration-infinitesimalise}
\HEAD{\ref{reintegration-infinitesimalise}.~{\sf R\'eint\'egration de 
l'infinit\'esimalis\'e: mesure de courbure
}
}{
Jo\"el {\sc Merker} et Jean-Jacques {\sc Szczeciniarz}}

Toutes ces quantit\'es $X$, $Y$, 
$\omega$, $\Omega$, $k$ sont infinit\'esimales
et ponctuelles. \`A un certain moment, nous avons fait observer
que l'infinit\'esimalisation \'etait une simplification d\'ecid\'ee
pour les besoins d'une exploration.
Et rien ne nous garantissait que notre {\em conception} du
d\'efaut de commutativit\'e du transport parall\`ele
serait compl\`ete et achev\'ee de la sorte. Au contraire,
nous pouvions craindre que l'emploi de formes g\'eom\'etriques
aussi simples que les parall\'elogrammes infinit\'esimaux
risqu\^at de nous faire perdre l'espoir de comprendre
ce qu'il se passe au niveau macroscopique.
Contre toute attente, nous allons maintenant argumenter
que rien du concept r\'eel n'a \'et\'e perdu en chemin, 
et ce, gr\^ace \`a des proc\'ed\'es d'{\em int\'egration}. 

Soit $D \subset S$ un domaine bord\'e par un lacet
diff\'erentiable $\lambda := 
\partial D$. Tout ce que nous avons pr\'esent\'e montre que
la variation d'angle d'un transport parall\`ele le long de 
$\lambda$ est donn\'ee par la formule: 
\[
\theta(\lambda)
\,=\,
\int_\lambda\,
\omega
\eqno
{\scriptstyle{(\modsmall\,2\pi)}}.
\]
Mais alors le th\'eor\`eme de Green-Riemann-Stokes permet de
la r\'e-exprimer comme int\'egrale de surface, dans le domaine
$D$, de la diff\'erentielle ext\'erieure de $\omega$:
\[
\aligned
\theta(\lambda)
&
\,=\,
\int\!\!\int_D\,
d\omega
\\
&
\,=\,
\int\!\!\int_D\,
\Omega,
\endaligned
\]
qui n'est autre que la $2$-forme de courbure invariante
$\Omega$, cette int\'egrale n'\'etant en g\'en\'eral pas nulle, 
puisque la courbure n'est pas toujours identiquement nulle.

Nous savions que le changement d'angle est en g\'en\'eral
non nul le long d'un lacet, et nous avions implicitement
compris aussi qu'il ne d\'ependait pas du choix
de rep\`eres orthonorm\'es $\{X, Y\}$ ou $\{X', Y'\}$, puisqu'une
int\'egration de la relation:
\[
\omega'
\,=\,
\omega
-
d\theta,
\]
donnait instantan\'ement la co\"{\i}ncidence:
\[
\aligned
\int_\lambda\,\omega'
&
\,=\,
\int_\lambda\,\omega
-
\int_\lambda\,d\theta
\\
&
\,=\,
\int_\lambda\,\omega
-
0,
\endaligned
\]
sachant que sur une courbe ferm\'ee, l'int\'egrale d'une d\'eriv\'ee
est forc\'ement nulle puisque les valeurs de la primitive 
sont les m\^emes aux extr\'emit\'es \'egales. 

Mais maintenant, d'un point de vue conceptuel, tout s'\'eclaire:
$\Omega$ est une quantit\'e diff\'erentielle $2$-dimensionnelle
{\em invariante}, et gr\^ace \`a la transformation de Green-Riemann-Stokes,
le changement d'angle s'exprime en fonction d'une int\'egrale
que Gauss appelait la {\sl mesure de courbure:}
\[
\theta(\lambda)
\,=\,
\int\!\!\int_D\,
\Omega
\eqno
{\scriptstyle{(\modsmall\,2\pi)}},
\]
o\`u l'orientation de $\lambda$ est choisie de telle sorte
que $D$ se situe du c\^ot\'e gauche quand on parcourt $\lambda = 
\partial D$. 
En cons\'equence, l'holonomie est non r\'eduite \`a l'identit\'e
d\`es qu'il y a de la courbure, \`a savoir, d\`es que la courbure
de Gauss $k$ est non identiquement nulle. 

L'invariance du concept de courbure s'exprime d'ailleurs
\`a un niveau tellement \'elev\'e que ces int\'egrales:
\[
\int\!\!\int_D\,
\Omega
\,\,=\,\,
\int\!\!\int_{-D}\,
-\,\Omega
\]
ont une valeur inchang\'ee lorsque $S$ est munie d'une orientation
oppos\'ee.

\Section{\bf Théorème de Gauss-Bonnet dans le plan euclidien $\R^2$}
\label{Gauss-Bonnet-plan-R-2}
\HEAD{\ref{Gauss-Bonnet-plan-R-2}.~{\sf Théorème de 
Gauss-Bonnet dans le plan euclidien $\R^2$}
}{
Jo\"el {\sc Merker} et Jean-Jacques {\sc Szczeciniarz}}

Lorsqu'on munit le plan $\R^2$ de la métrique pythagoricienne standard:
\[
\langle
X,\,Y
\rangle
\,:=\,
X_1Y_1
+
X_2Y_2,
\]
un calcul simple montre que la forme de courbure
associée $\Omega \equiv 0$ est identiquement nulle.

En ce moment présent, nous proc\'edons d'une mani\`ere qui est
fr\'equente en math\'ematiques. Du point de vue d'un d\'eveloppement
plus conceptuel et plus complet, nous venons en effet de r\'eexprimer
une situation et un r\'esultat qui était connus par ailleurs, et qui
paraissent maintenant \'el\'ementaires, rétrospectivement. Cette
proc\'edure est une fa\c{c}on théorique de traiter et d'envisager {\sl
par l'apr\`es}, ou {\sl par le haut}, le fondement même d'un
r\'esultat \'el\'ementaire. Ici donc à un niveau supérieur, nous
concevons dorénavant mieux ce que signifie {\sl \^etre planaire}, ou
même {\sl \^etre euclidien}. Précisément, nous pouvons maintenant nous
représenter qu'il s'agit d'une surface pour laquelle la forme de
courbure agissant sur deux champs de vecteurs est identiquement nulle.

Essayons alors de mesurer le chemin conceptuel que nous avons
accompli, et qui s'est accompli de lui-même par un mouvement
intrinsèque autonome. Le plan est maintenant vu \`a partir d'un cadre
conceptuel beaucoup plus ample, qui permet d'explorer en m\^eme temps
que lui d'autres surfaces. Il a d'abord \'et\'e la surface de
r\'ef\'erence, mais il est devenu un cas particulier parmi les
surfaces de diff\'erentes courbures possibles. Or un ensemble de
concepts a \'et\'e d\' egag\'e qui n'apparaissait pas ou ne pouvait
pas appara\^itre au seul examen du plan. Et maintenant, le plan se
trouve pourvu de ces propri\'et\'es nouvelles qui n'apparaissent pas
spontanément, tout comme elles n'apparaissaient pas initialement. Ces
propriétés ont donc chang\'e de statut: 
{\em leur inexistence n'est plus la m\^eme}.

Avant d'exposer le théorème de Gauss-Bonnet général
dont l'essence réelle {\em est} de tenir compte
de la présence d'une courbure de Gauss
$k \neq 0$ non nulle, il est instructif de présenter
ce qu'il se passe quand $k \equiv 0$, avec des hypothèses
géométriques assez riches.

Soit donc $\gamma \subset \R^2$ une courbe différentiable
fermée simple, à tangente nulle part nulle.
Le théorème de Jordan différentiable montre que $\gamma$
borde un domaine $P \subset \R^2$ borné, ayant
la topologie d'un disque (déformé).
\footnote{\ Les figures qui suivent sont extraites
du livre de Lehmann-Sacr\'e, 
{\em G\'eom\'etrie et topologie des surfaces}, 
Presses Universitaires de France, Paris, 1982, 348~pp.}
 
\begin{center}
\includegraphics[scale=0.75]{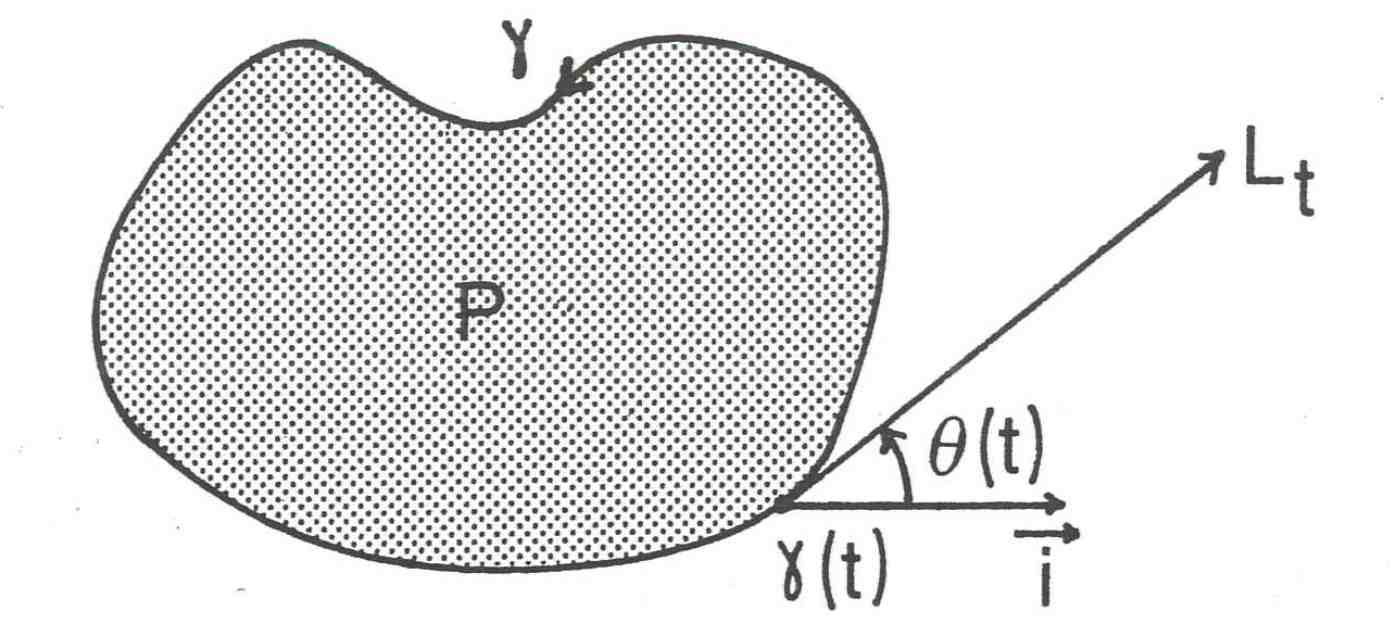}
\end{center} 

Paramétrons $\gamma \colon [0, 1] \longrightarrow \R^2$,
avec $\gamma(0) = \gamma(1)$. En tout point $t \in [0,1]$,
notons:
\[
\theta(t)
\,:=\,
\Angle\,
\Big(
\frac{d\gamma}{dt}(t),\,\,
\text{axe}\,0x
\Big)
\]
l'angle orienté que fait le vecteur tangent avec la direction
horizontale représentée par l'axe $0x$. Le résultat
élémentaire qui constitue le germe fécond de la théorie
de Gauss-Bonnet est un

\begin{Theoreme}
La variation totale de l'angle $\theta(t) = \Angle\,
\big( \frac{d\gamma}{dt}, 0x \big)$ le long de la courbe
fermée simple $\gamma$ vaut toujours:
\[
2\pi
\,=\,
\int_0^1\,
\frac{d\theta}{dt}(t)\,
dt
\,=\,
\int_\gamma\,
d\theta.
\eqno\qed
\]
\end{Theoreme}

Ainsi, lorsqu'on parcourt $\gamma$, l'angle effectue un tour
comple $2\pi = 360^{\text o}$. Après reparamétrisation
de $\gamma$ par longueur d'arc:
\[
\gamma(s),
\ \ \ \ \ \ \ \ \ \ \ \ \ \ \ \ \ \
s\,\in\,[0,\LL],
\ \ \ \ \ \ \ \ \ \ \ \ \ \ \ \ \ \
\LL
\,=\,
\longueur(\gamma),
\]
la dérivée:
\[
\frac{d\theta}{ds}(s)
\,=:\,
\kappa(s)
\]
n'est autre que la {\sl courbure} de $\gamma$ au point $\gamma(s)$.

\begin{center}
\includegraphics[scale=0.75]{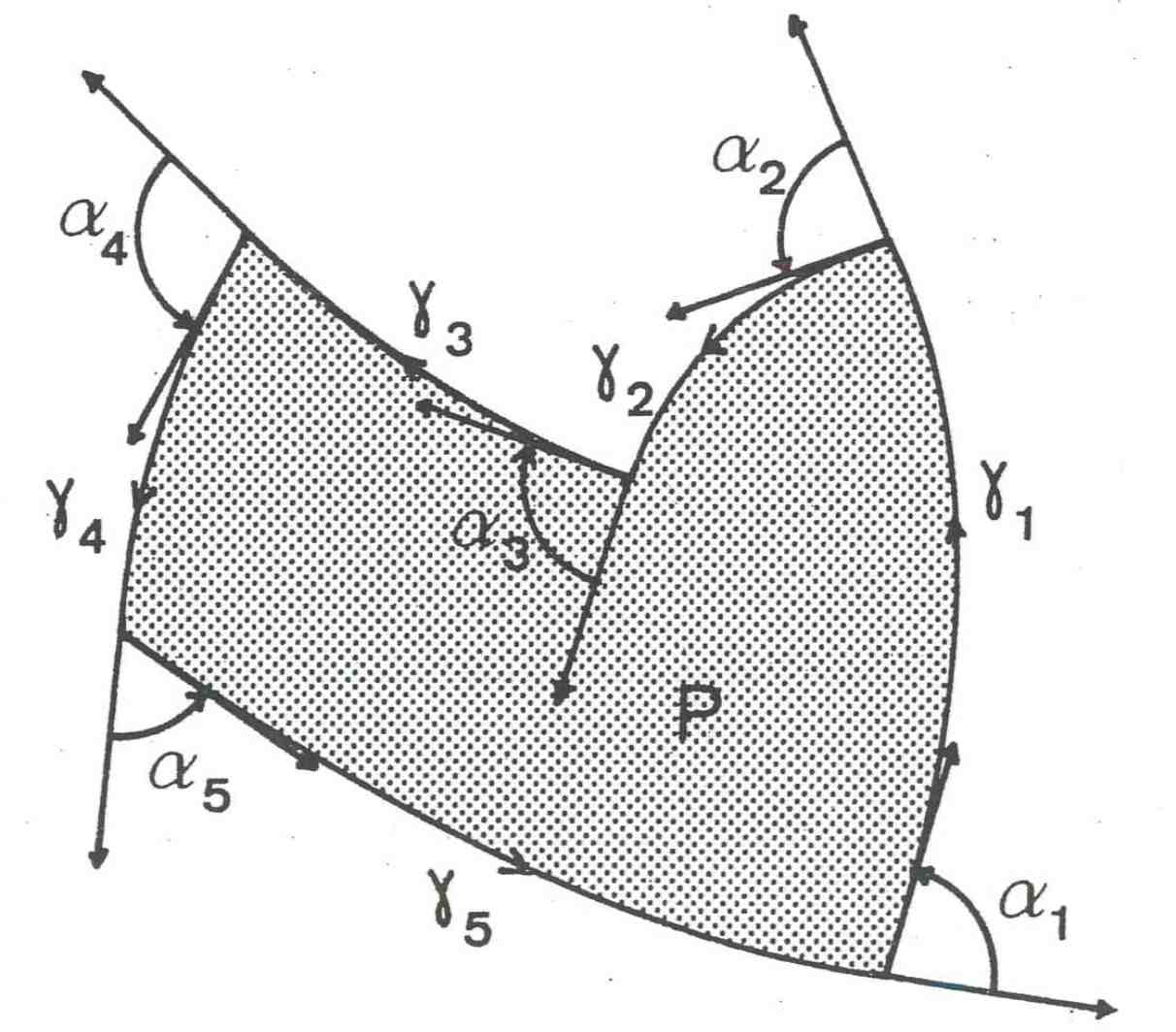}
\end{center} 

On a consid\'er\'e la variation d'un angle le long de la courbe, la
variation \'etant d\'etermin\'ee par la mani\`ere dont l'angle suit
son contour. Cette variation totale vaut $2\pi$, comme le th\'eor\`eme
l'a indiqu\'e. Il s'agit de comprendre cette variation angulaire au
cours du parcours. C'est aussi un changement en direction d'une
tangente \`a la courbe qui fait un angle avec un axe. La courbure
d'une courbe est bien perceptible, comme propri\'et\'e d'un
parcours. De mani\`ere quasi anthropomorphique, on dit que la tangente
`{\sl \'eprouve}' des discontinuit\'es angulaires, au point de se
d\'edoubler. D'abord, il faut noter que la courbure est ce qui fait
qu'une courbe est {\sl courbe}. 
Elle est comme immanente \`a la courbe,
qu'elle anime de l'int\'erieur. Et aux points anguleux,
c'est comme si la courbe se perdait dans sa nature m\^eme de
courbe. Une fois
d'un c\^ot\'e et une fois de l'autre en le m\^eme point:
il y a en effet
deux tangentes en ces points. Il faudrait sans doute dire que
dans ce domaine de la virtualit\'e, il existe en ces points deux
mani\`eres d'\^etre nul qui se cumulent ou se concentrent.
 
En suivant la remarque que nous avons faite plus haut sur l'ontologie
du virtuel\,\,---\,\,l'ontologie de l'inexistence\,\,---, il se peut
que la concentration de courbure fasse \'echo en sym\'etrique au
calcul habituel de la courbure en d\'eriv\'ee seconde.

Plus généralement, en supposant que $\gamma$ est seulement
différentiable par morceaux, la formule devient:
\leqnomode\usetagform{default}
\begin{align}
\label{2-pi-kappa-alpha}
2\,\pi
\,=\,
\int_\gamma\,\kappa(s)\,ds
+
\sum_{1\leqslant m\leqslant\MM}\,
\alpha_m,
\end{align}
où $\alpha_1, \dots, \alpha_\MM$ sont les angles extérieurs
entre deux morceaux lisses consécutifs de $\gamma$.
\`A ces endroits, la tangente éprouve des discontinuités
angulaires. D'une certaine manière, il se produit
en ces points une {\em concentration infinie de courbure}, 
comme nous allons l'expliquer et le diagrammatiser.

\begin{center}
\includegraphics[scale=0.75]{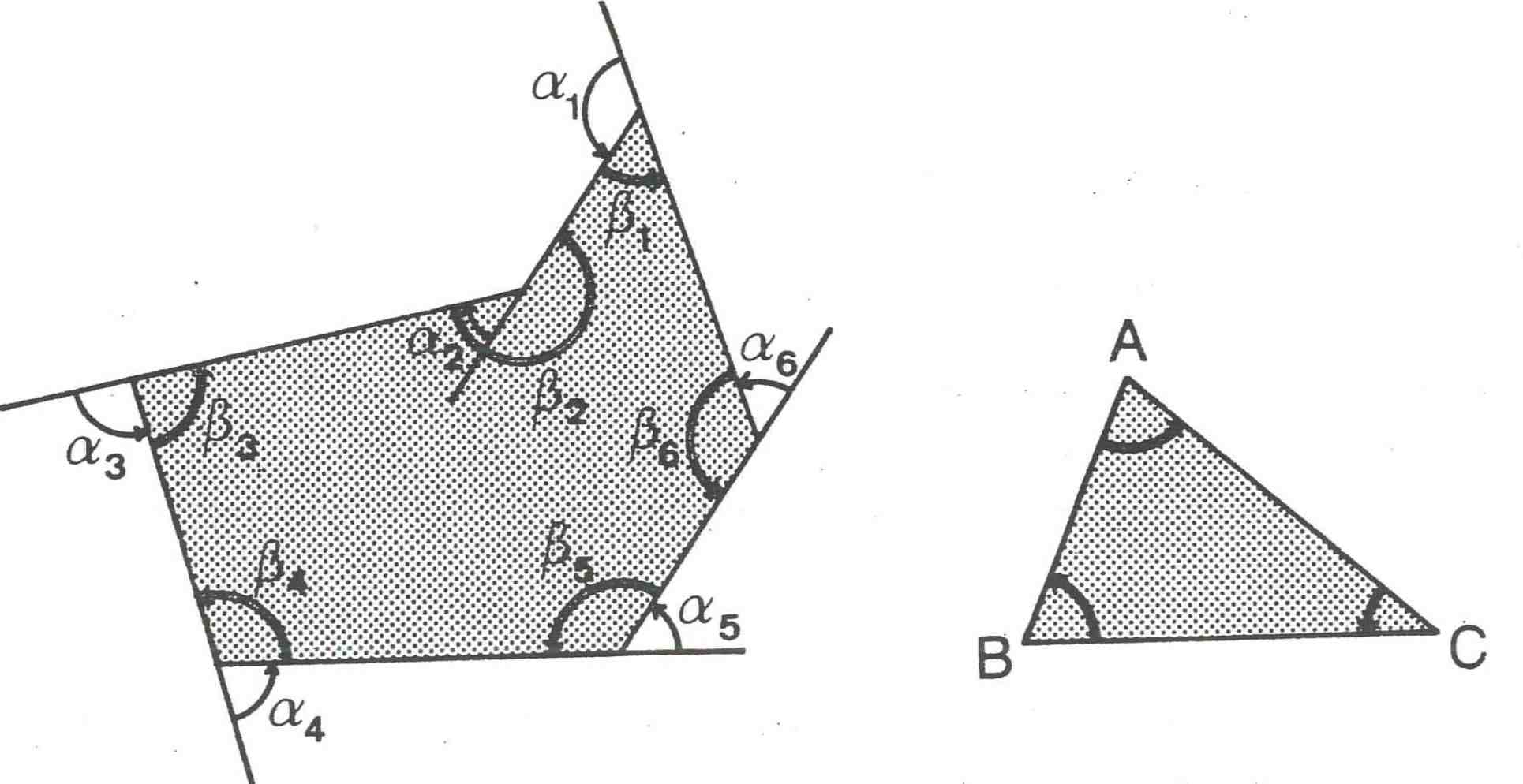}
\end{center} 

Pour un polygone $P \subset \R^2$ dans le plan euclidien,
dont les côtés, rectilignes, n'ont en eux-mêmes aucune
courbure, les intégrales de $\kappa(s)$ sur les morceaux
lisses du bord $\gamma = \partial P$ sont toutes nulles,
et il ne reste plus que:
\[
2\pi
\,=\,
0
+
\sum_{1\leqslant m\leqslant\MM}\,
\alpha_m.
\]
En introduisant les angles intérieurs:
\[
\alpha_m
\,:=\,
\pi
-
\beta_m,
\]
on en déduit (on retrouve) un théorème classique de géométrie
euclidienne: la somme des angles intérieurs à un polygone
fermé simple vaut:
\[
\big(\MM-2\big)\,
\pi
\,=\,
\sum_{1\,\leqslant\,m\,\leqslant\,\MM}\,
\beta_m.
\]

En particulier, lorsque $\MM = 3$, à savoir pour un triangle,
on retrouve la célèbre formule d'Euclide:
\[
\pi
\,=\,
\widehat{A}
+
\widehat{B}
+
\widehat{C}.
\]

\begin{center}
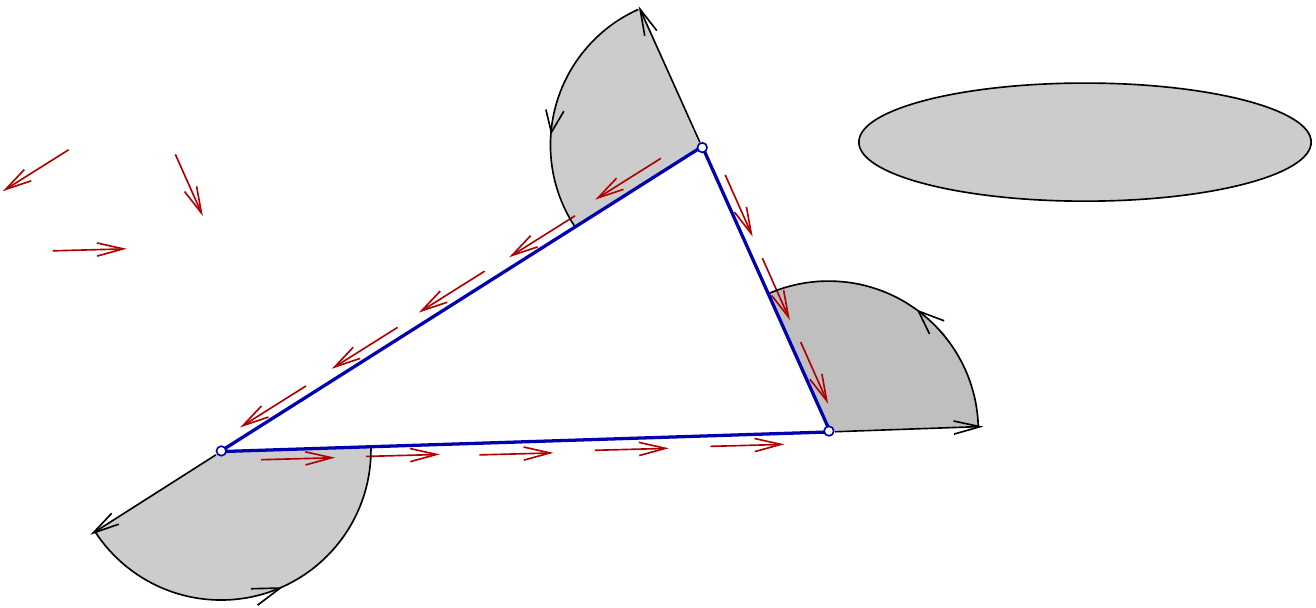
\end{center}

Rétrospectivement, à la lumière de la formule 
générale~({\ref{2-pi-kappa-alpha}}) que nous venons d'écrire,
il faut comprendre que la courbure 
est {\em concentrée} en les trois sommets $A$, $B$, $C$
du triangle. Cette observation de compréhension
est bien entendu valable aussi pour des polygones simples
à un nombre quelconque $\MM \geqslant 3$ de côtés.

Cette premi\`ere pr\'esentation pr\'eliminaire au th\'eor\`eme de
Gauss-Bonnet g\'en\'eral est en fait d\'ej\`a un \'enonc\'e de ce
th\'eor\`eme. Nous avons au d\'epart une autre sorte de
g\'en\'eralisation du th\'eor\`eme d'Euclide.

Nous avons consid\'er\'e essentiellement un concept de courbure qui
correspond imm\'ediatement \`a la mani\`ere dont l'angle $\theta$
varie en fonction de l'abscisse curviligne, le long de la
courbe parcourue.

Maintenant, regardons comment le concept de
courbure va se mettre \`a l'{\oe}uvre.

\Section{\bf Théorème de Gauss-Bonnet local sur une surface
riemannienne $S$}
\label{Gauss-Bonnet-local-S-riemannienne}
\HEAD{\ref{Gauss-Bonnet-local-S-riemannienne}.~{\sf Théorème 
de Gauss-Bonnet local sur une surface
riemannienne $S$}
}{
Jo\"el {\sc Merker} et Jean-Jacques {\sc Szczeciniarz}}

Qu'en est-il sur une surface $S$ munie d'une
connexion métrique quelconque $\nabla$? La fonction 
de courbure $k$, ou, de manière équivalente, 
la $2$-forme de courbure $\Omega$ sur $S$, va introduire (insérer)
un terme nouveau et central dans la 
formule~({\ref{2-pi-kappa-alpha}}).

Comme précédemment, soit $\gamma \colon\, [0, \LL] \longrightarrow S$
une courbe différentiable fermée simple à tangente nulle part
nulle, paramétrée par longueur d'arc:
\[
1
\,\equiv\,
\Bignorm
\frac{d\gamma}{ds}(s)
\Bignorm
\eqno
{\scriptstyle{(\forall\,s\,\in\,[0,\LL])}}.
\]
Cette courbe $\gamma = \partial D$ borde un certain domaine
$D \subset S$ difféomorphe à un disque. On peut supposer que
$D$ se situe à gauche d'un point mobile parcourant $\gamma$.
Abrégeons ces vecteurs tangents:
\[
\tau(s)
\,:=\,
\frac{d\gamma}{ds}(s).
\]
Soit aussi:
\[
\nu(s)
\,:=\,
\Rotation_{\frac{\pi}{2}}
\big(\tau(s)\big)
\]
le champ de vecteurs unitaire le long de $\gamma$
rentrant dans le domaine $D$.

Alors l'invariance de la métrique sous l'action de $\nabla$
assure une relation d'orthogonalité qui implique
que $\nabla_{\tau(s)}(\tau)$ est orthogonal à $\tau(s)$, 
et par conséquent on peut écrire:
\[
\nabla_{\tau(s)}\big(\tau\big)(s)
\,=\,
k_g(s)
\cdot
\nu(s),
\]
au moyen d'une certaine fonction invariante $s \longmapsto k_g(s)$, 
appelée {\sl courbure géodésique} de $\gamma$.

\begin{center}
\includegraphics[scale=0.75]{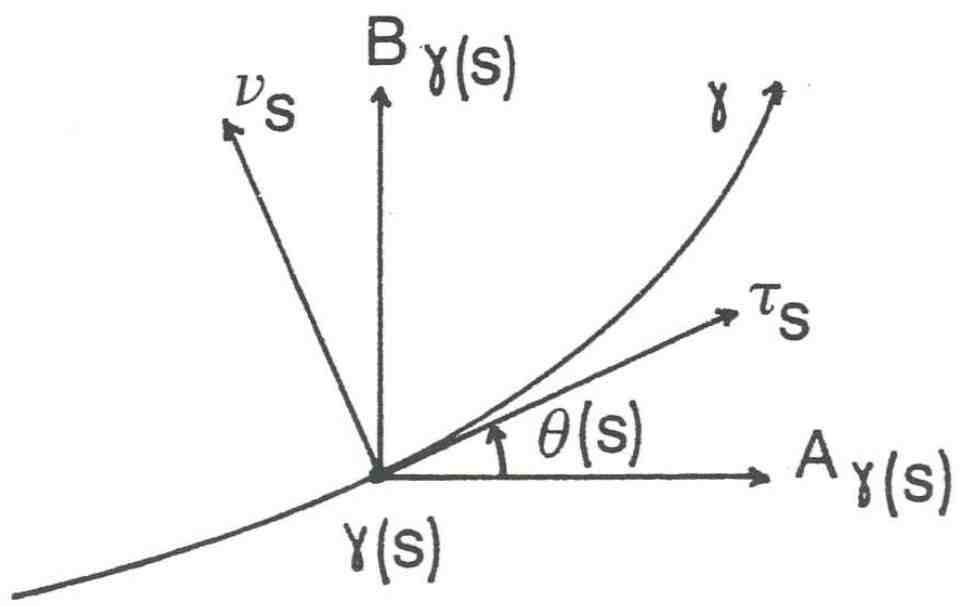}
\end{center} 

Pour comprendre ce concept
d'une manière géométrique et intuitive, il
vaut mieux envisager la surface
$S \subset \R^3$ comme {\em extrins\`eque}, {\em i.e.}
plong\'ee dans l'espace euclidien
tridimensionnel $\R^3$. Donnons-nous un champ de vecteurs
unitaires $s \longmapsto X(s)$ sur une courbe $\gamma 
\colon [0,\LL] \rightarrow S \subset \R^3$ param\'etr\'ee 
par longueur d'arc. 
Comme tous les $X(s)$ pour $s \in [0, \LL]$
sont de norme $1$, les vecteurs
dérivés $\frac{dX}{ds}(s)$ sont orthogonaux 
aux vecteurs $X(s)$. Si donc $N(s)$ désigne un champ
de vecteurs unitaires {\em orthogonaux}
à la surface $S$ aux points $X(s)$ pour
une certaine orientation de $S$, alors 
les propriétés standard du produit
vectoriel $(\centersmallbullet) \wedge (\centersmallbullet)$
dans l'espace $\R^3$ assurent qu'il existe
des réels $\lambda(s)$ tels que la projection
orthogonale sur l'espace tangent
$T_{\gamma(s)}S$ de la dérivée dans $\R^3$ de $X(s)$ s'exprime comme:
\[
\frac{DX}{ds}
(s)
\,:=\,
\proj_{TS}^\bot
\bigg(
\frac{dX}{ds}(s)
\bigg)
\,=\,
\lambda(s)\,
N(s)
\wedge
X(s)
\eqno
{\scriptstyle{(\forall\,s\,\in\,[0,\,\LL])}}.
\]
Pour $s$ fixé, le nombre réel $\lambda(s)$ est alors appelé
{\sl valeur alg\'ebrique de la d\'eriv\'ee covariante} 
de $X$ en $s$.

Maintenant, si $X(s) := \frac{d\gamma}{ds}(s)$ est le champ de
vecteurs tangents à une courbe $\gamma \colon [0, \LL] \longrightarrow
S$ paramétrée par longueur d'arc comme ci-dessus, alors la valeur
alg\'ebrique $\lambda(s)$ de la d\'eriv\'ee covariante de $\gamma'(s)$
en $\gamma(s)$ est appel\'ee {\sl courbure g\'eod\'esique} de la
courbe $\gamma$ en son point $\gamma(s)$. Pour la connection $\nabla$
naturellement associée à $S \subset \R^3$, on démontre que ces nombres
{\em extrinsèques:}
\[
\lambda(s)
\,=\,
\pm\,
k_g(s)
\]
sont {\em égaux}\,\,---\,\,éventuellement au signe près\,\,---\,\,aux
nombres intrinsèques $k_g(s)$ introduits ci-dessus. 
D'un point de vue externe
\`a la surface, la valeur de la
courbure g\'eod\'esique $k_g(s)$ est la valeur
de la composante tangentielle du vecteur 
accélération $\gamma''(s)$.

C'est donc ici
tout une force de cohérence de la théorie mathématique qui s'exprime
dans les coïncidences conceptuelles modulées par des perspectives
ontologiques complémentaires.

En proc\'edant comme l'exige le travail math\'ematique, nous avons
ainsi introduit pour une courbe sur une surface l'analogue de la
courbure d'une courbe plane. L'ambiguïté sur le signe $\pm$
s'explique simplement par le fait que le concept de courbure
g\'eod\'esique que nous considérons d\'epend des
deux choix possibles pour une orientation de $S$.

Dans sa g\'en\'eralit\'e, le th\'eor\`eme de Gauss-Bonnet doit faire
intervenir la courbure g\'eod\'esique, comme nous l'avions \'enonc\'e
au d\'ebut de ce travail, puisqu'il consid\`ere des courbes sur une
surface qui ne sont pas des g\'eod\'esiques.

Maintenant, en revenant à la théorie intrinsèque plus
générale et plus puissante,
si le domaine $D \subset U$ est de plus contenu dans un ouvert
de carte $U \subset S$, dans lequel est donné un champ de
repères orthonormés des espaces tangents:
\[
T_pS
\,=\,
\R X_p
\overset{\bot}{\oplus}
\R Y_p,
\]
de telle sorte qu'il existe une fonction angulaire $s \longmapsto
\theta(s)$ satisfaisant:
\[
\aligned
\tau(s)
&
\,=\,
\cos\,\theta(s)\,
X_{\gamma(s)}
+
\sin\,\theta(s)\,
Y_{\gamma(s)},
\\
\nu(s)
&
\,=\,
-\,\sin\,\theta(s)\,
X_{\gamma(s)}
+
\cos\,\theta(s)\,
Y_{\gamma(s)},
\endaligned
\]
on démontre une

\begin{Proposition}
{\bf [Cruciale]}
Pour tout paramètre $s \in [0, \LL]$, la courbure géodésique vaut:
\[
k_g(s)
\,=\,
\frac{d\theta}{ds}(s)
-
\omega
\Big(
\frac{d\gamma}{ds}(s)
\Big),
\]
où $\omega = \omega_{X,Y}$ est la $1$-forme différentielle
associée à la différentiation covariante $\nabla$,
dans le champ de repères orthonormés $\big\{ X_p, Y_p\big\}_{p 
\in U \subset S}$.
\end{Proposition}

\proof
Avec toutes les données qui précèdent, un bref calcul suffit
à éclairer cette formule:
\begin{align}
\nabla_{\tau(s)}\big(\tau\big)(s)
&
\,=\,
\nabla_{\tau(s)}
\Big(
\cos\,\theta(s)\,
X_{\gamma(s)}
+
\sin\,\theta(s)\,
Y_{\gamma(s)}
\Big)
\notag
\\
&
\,=\,
\Big(
-\,\sin\,\theta\,
X_{\gamma}
+
\cos\,\theta\,Y_\gamma
\Big)\,
\frac{d\theta}{ds}
+
\cos\,\theta\,
\nabla_\tau(X)
+
\sin\,\theta\,
\nabla_\tau(Y)
\notag
\\
&
\,=\,
\Big[
\frac{d\theta}{ds}(s)
-
\omega\big(\tau(s)\big)
\Big]\,
\nu(s).
\qedhere
\end{align}
\endproof

\begin{center}
\includegraphics[scale=0.75]{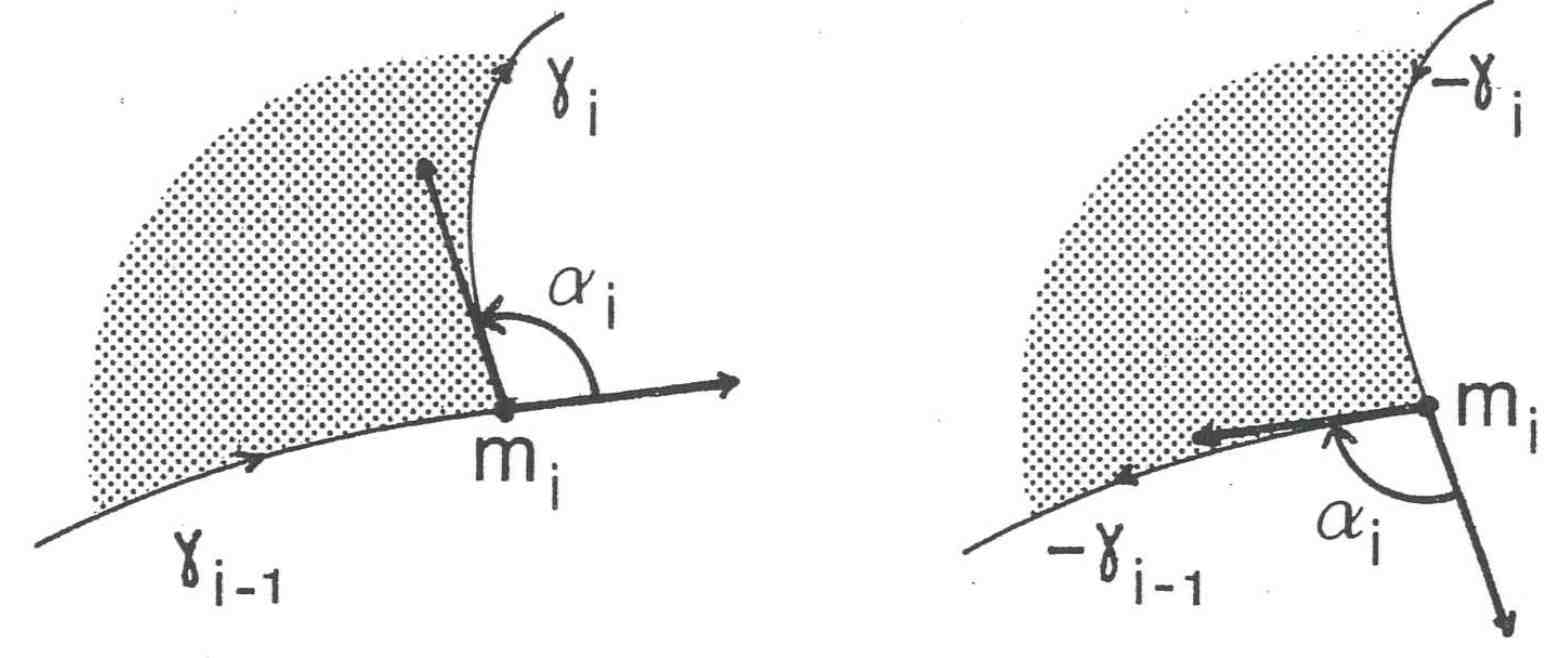}
\end{center} 

Nous pouvons enfin énoncer, sans entrer dans les détails de
la démonstration, le 

\begin{Theoreme}
{\bf [de Gauss-Bonnet local]}
Soit une courbe fermée simple $C = \partial D$ bord orienté
lisse par morceaux d'un domaine $D \subset U$ contenu
dans un ouvert de carte $U \subset S$ d'une surface
riemannienne abstraite $S$. Alors:
\[
2\pi
-
\int\!\!\int_D\,
\Omega
\,\,=\,\,
\int_{\partial D}\,
k_g(s)\,ds
+
\sum_{1\leqslant m\leqslant\MM}\,
\alpha_m.
\]
\end{Theoreme}

\begin{center}
\includegraphics[scale=0.75]{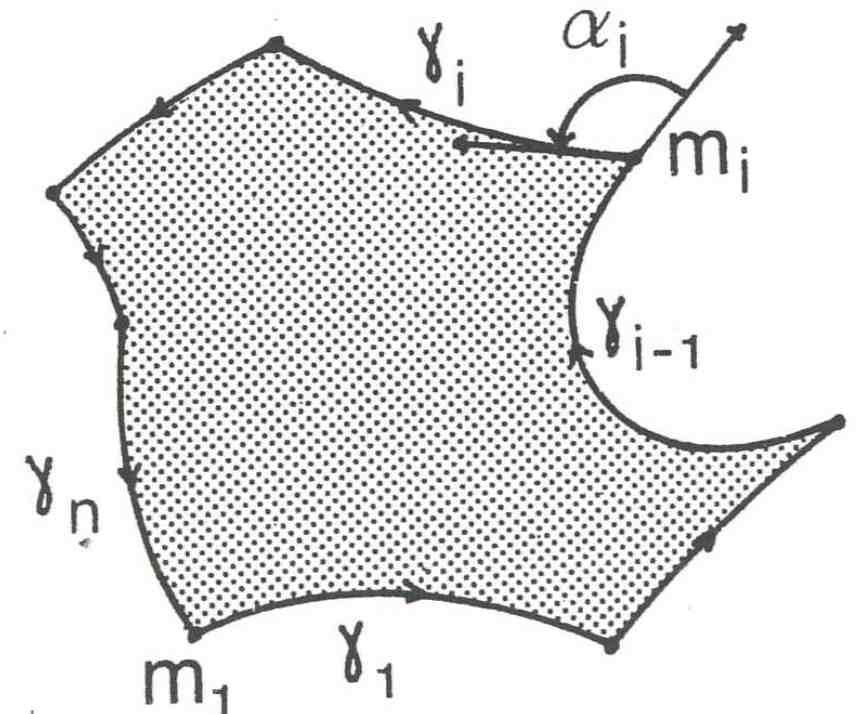}
\end{center} 

On doit remarquer que la complication de la d\'efinition qui va de
pair avec son insertion dans des champs diff\'erents du corpus de la
g\'eom\'etrie diff\'erentielle et de la topologie alg\'ebrique,
arrache certains \'el\'ements de l'appr\'ehension initiale (intuition)
de la notion de courbure pour les r\'einvestir dans de nouveaux
\'el\'ements de d\'efinition. De ce fait notre connaissance de la
g\'eom\'etrie et de la topologie de la surface ou de la vari\'et\'e
s'accro\^it. Il reste remarquable que nous parlions encore de
<<\,courbure\,>>, tout comme nous parlons de l'<<\,espace\,>>.

Si nous consid\'erons ce r\'esultat, plus puissant et synth\'etique de
plusieurs domaines des math\'ematiques mettant en relation (nous
avons insist\'e sur ces relations d\'efinitionnelles), des concepts
qui se con\c coivent dans chacune des disciplines dont ils sont issus,
essentiellement, g\'eom\'etrie, topologie, topologie alg\'ebrique,
nous pouvons essayer de voir s'il s'agit bien de la m\^eme notion de
courbure qui est en jeu.

Il faut noter, et c'est non n\'egligeable, que dans le cas de notions
du m\^eme type issues de la g\'eom\'etrie \'el\'ementaire plus
intuitive, l'objet initial qu'elles rappellent comme en \'echo
continue d'insister sous cette premi\`ere forme que d'aucuns
pourraient nommer {\sl allusive}. 
Non pas seulement parce que nous pouvons y
revenir en particularisant la notion g\'en\'eralis\'ee. Mais parce que
la forme initiale a acquis une puissance abstraite qu'elle alimente ou
qu'elle peut alimenter. C'est l\`a une des mani\`eres dont agit
l'imaginaire math\'ematique. Il s'agit de ce que nous avons appel\'e
du {\sl symbolico-intuitif}. Nous pouvons, par exemple, user de
repr\'esentations perceptives, voir intuitives, pour travailler \`a un
niveau plus \'elev\'e d'abstraction, en ne conservant de cette
perception que ce qui avec bonheur continue de s'adapter \`a ce
niveau.

Et c'est pr\'ecis\'ement \`a travers ces formes de synth\`eses ainsi
fortement diversifi\'ees que nous acqu\'erons une force de propulsion.

Le th\'eor\`eme  de Gauss Bonnet
pr\'esente une g\'eom\'etrie assez simple, 
et la difficult\'e de la preuve tient \`a certains faits
topologiques qui lui sont int\'egr\'es. C'est donc par cette
synth\`ese m\^eme qu'il est relativement difficile et profond.
Comme l'exprime Spivak cit\'e en exergue au début de cet article,
le th\'eor\`eme  de Gauss Bonnet
est probablement le th\'eor\`eme le plus profond de la
g\'eom\'etrie diff\'erentielle. On se souvient de sa premi\`ere
version chez Gauss. L'exc\`es sur $\pi$ de la somme des angles
int\'erieurs $\alpha_{1}, \alpha_{2}, \alpha_{3}$ d'un triangle
g\'eod\'esique $T$ est \'egal \`a l'int\'egrale de la courbure
gaussienne $K_\GG$ sur $T$:
\[
\alpha_1
+
\alpha_2
+
\alpha_3
-
\pi
\,=\,
\int\!\!\int_{T} 
K_\GG\,
dA.
\]

On peut voir directement sur cette formule simple une extension du
th\'eor\`eme de Thal\`es et d'Euclide dans le cas o\`u $K_\GG
\cong 0$. On
obtient le cas de la sph\`ere unit\'e, la somme des angles
int\'erieurs d'un triangle g\'eod\'esique est sup\'erieur \`a $\pi$, et
l'exc\`es sur $\pi$ est exactement l'aire de $T$. 

Il existe des
pr\'esentations
classiques du th\'eor\`eme de Gauss-Bonnet local qui n'utilisent
pas
le formalisme de le forme de courbure. 
Ces pr\'esentations insistent sur la topologie de la surface
consid\'er\'ee, en usant du th\'eor\`eme des tangentes tournantes: la
variation totale de l'angle d'un vecteur tangent \`a une courbe
$\gamma \colon [0, \LL]\rightarrow S$ donn\'ee par une fonction
continue allant d'un intervalle ferm\'e $[0, \LL]$ sur une surface
différentiable $S$. 
Et le th\'eor\`eme de Gauss-Bonnet est \'enonc\'e
pour une r\'egion de $S$ born\'ee.

Le th\'eor\`eme \'enonc\'e ci-dessus est plus abstrait et use du
formalisme des formes différentielles. 
L'int\'egrale de la mesure de courbure selon
Gauss est l'int\'egrale de la forme de courbure.

La courbure que nous pouvons percevoir dans sa ph\'enom\'enologie,
comme les modalit\'es d'un contour d'une courbe sur une surface a
d'abord \'et\'e conceptualis\'ee comme un exc\`es sur $\pi=$ somme des
angles int\'erieurs d'un triangle euclidien (g\'eod\'esique). Cette
conceptualisation est d\'ej\`a un saut important, la courbure a
\'et\'e m\'etriquement traduite. Notre appr\'ehension de l'espace a
trouv\'e un e modalit\'e de quantification. Notre attention perceptive
s'est focalis\'ee sur une premi\`ere courbe ferm\'ee (un triangle)
dont nous essayions de capter les formes de parcours en tenant compte
des angles, jusqu'\`a sa fermeture. On peut dire que sous cette
analyse la courbure a \'etendu et r\'ealis\'e un domaine
d'application. Mais du m\^eme coup elle a ouvert ou r\'eouvert la
r\'eflexion topologique dont elle est porteuse. Dans un mouvement
th\'eorique qui semble aller \`a rebours de l'appr\'ehension kantienne,
la courbure s'est int\'egr\'ee \`a notre appr\'ehension de l'espace,
dont elle se fait une condition.

Avec le concept de forme de courbure, nous sommes pass\'es \`a un stade
sup\'erieur d'abstraction. On dispose d'une $2$-forme qui d\'etermine
enti\`erement la courbure $K_\GG$. 
Elle est construite \`a partir d'une $1$-
forme associ\'ee \`a deux champs de vecteurs. La
différentielle extérieure $d
(\centersmallbullet)$ de cette 1-forme
nous donne la forme de courbure

La $2$-forme de courbure a \'et\'e pr\'esent\'ee 
dans la Section~{\ref{reintegration-infinitesimalise}}.
Elle
nous fait b\'en\'eficier des propri\'et\'es des $2$-formes
diff\'erentielles renvoyant \`a la th\'eorie des formes
diff\'erentielles. Et elle est directement reli\'ee \`a l'holonomie
restreinte que nous avons pr\'esent\'ee. Dans chacune de ces
d\'eterminations, c'est une fa\c{c}on de concevoir un aspect de la
courbure qui trouve une explication. Et c'est aussi une nouvelle
appr\'ehension des anciennes formes d'appr\'ehension du concept de
courbure: extr\^eme continuit\'e, extr\^eme nouveaut\'e, dans ces
synth\`eses expansives.
 
La construction conceptuelle fond\'ee sur une $2$-forme permet de
conserver l'organisation des champs de vecteurs \`a laquelle renvoie
la $1$-forme toujours mobilisable dans notre compr\'ehension. Et cette
r\'ef\'erence aux champs nous permet de suivre une forme conceptuelle
qui compose la courbure.  Une nouvelle ph\'enom\'enologie
transcendantale peut \^etre d\'evelopp\'ee.  En ce sens la
compr\'ehension primitive (intuitive) de la courbure a \'et\'e
absorb\'ee et continue d'agir.

Nous avons \'egalement montr\'e que la question des fondements de la
g\'eom\'etrie (euclidienne) \'etait enti\`erement transform\'ee dans
ces nouveaux cadres th\'eoriques. La g\'eom\'etrie est fond\'ee par le
moyen de ses extensions synth\'etiques, qui en sont autant
d'accroissements de connaissances.\footnote{\
Une seconde partie de ce texte, consacrée aux aspects topologiques et
réellement globaux du théorème de Gauss-Bonnet,
paraîtra ultérieurement.}

\medskip\noindent
{\bf Remerciements.}
Christian Houzel a relu le manuscrit.

\linestop


\bibliographystyle{amsplain}

\begin{thebibliography}{10}

{\scriptsize

{\bf\bibitem{Berger-Gostiaux-1992} 
{\sc Berger}}, M.; {\sc Gostiaux}, B.:
{\em Géométrie différentielle: variétés, courbes et surfaces}.
Deuxi\`eme édition. Math\'ematiques. 
Presses Universitaires de France, Paris, 1992, vi+513~pp. 

\smallskip

{\bf\bibitem{Cartan-1967} 
{\sc Cartan}}, H.:
{\em Formes diff\'erentielles. Applications \'el\'ementaires 
au calcul des variations et \`a la th\'eorie des courbes 
et des surfaces}. Hermann, Paris 1967, 186~pp.

\smallskip

{\bf\bibitem{Cavailles-1994} 
{\sc Cavaill\`es}}, J.:
{\em M\'ethode axiomatique et formalisme}. 
{\em Sur la logique et la th\'eorie de la science}.
In:
{\em {\OE}uvres compl\`etes de philosophie des sciences}, 
Paris, Hermann, 1994, vii+686~pp.

\smallskip

{\bf\bibitem{Caveing-2004} 
{\sc Caveing}}, M.:
{\em Le probl\`eme des objets dans la pensée mathématique}. 
Probl\`emes et Controverses. 
Librairie Philosophique J. Vrin, Paris, 2004.

\smallskip

{\bf\bibitem{Chern-1967} 
{\sc Chern}}, M.:
{\em Curves and surfaces in Euclidean space}. 
1967 Studies in Global Geometry and Analysis, 
pp.~16--56. Math. Assoc. Amer. (distributed by Prentice-Hall, 
Englewood Cliffs, N.J.).

\smallskip

{\bf\bibitem{do-Carmo-2016} 
{\sc do Carmo}}, M.P.:
{\em Differential geometry of curves \& surfaces}.
Revised \& updated second edition. Dover Publications, Inc., 
Mineola, NY, 2016, xvi+510~pp.

\smallskip

{\bf\bibitem{Dombrowski-1979} 
{\sc Dombrowski}}, P.:
{\em 150 years after Gauss' "{\sl Disquisitiones generales circa 
superficies curvas}". 
With the original text of Gauss}. 
Astérisque, 62. Société Mathématique de France, 
Paris, 1979, iv+153~pp.

\smallskip

{\bf\bibitem{Euclide-1998} 
{\sc Euclide d'Alexandrie}}:
{\em Les éléments}. Vols. I, II, III, IV. 
Traduits du texte de I.L. Heiberg par B. Vitrac, 
Biblioth\`eque d'Histoire des Sciences. 
Presses Universitaires de France, Paris, 1998, viii+433~pp.

\smallskip

{\bf\bibitem{Geenberg-1993} 
{\sc Greenberg}}, M.J.:
{\em Euclidean and non-Euclidean geometries. Development and history}. 
Third edition. W. H. Freeman and Company, New York, 1993, xviii+483~pp.

\smallskip

{\bf\bibitem{Hilbert-Cohn-Vossen-1952} 
{\sc Hilbert}}, D.; {\sc Cohn-Vossen}, S.:
{\em Geometry and the imagination}. 
Translated by P. Neményi. 
Chelsea Publishing Company, New York, N. Y., 1952, ix+357~pp.

\smallskip

{\bf\bibitem{Kant-1987} 
{\sc Kant}}, E.:
{\em Critique de la raison pure.} 
Trad. A. Tremeseygues \& B. Pacaut B. Presses
Universitaires de France, Paris,
1944. 
Trad. J. Barni, revue par P. Archambault, Flammarion, Paris, 1987.

\smallskip

{\bf\bibitem{Kant-1995} 
{\sc Kant}}, E.:
{\em Critique de la faculté de juger}. 
Trad. A. Philonenko, Vrin, Paris, 1968.
Trad. A. Renaut, Aubier, Paris, 1995.

\smallskip

{\bf\bibitem{Lehmann-Sacre-1982} 
{\sc Lehmann}}, D.; {\sc Sacr\'e}, C.:
{\em G\'eom\'etrie et topologie des surfaces}, 
Presses Universitaires de France, Paris, 1982, 349~pp.

\smallskip

{\bf\bibitem{Spivak-1979-I} 
{\sc Spivak}}, M.:
{\em A comprehensive introduction to differential geometry}, 
Vol. I. Second edition. Publish or Perish, Inc., Wilmington, Del., 
1979, xiv+668~pp.

\smallskip

{\bf\bibitem{Spivak-1979-II} 
{\sc Spivak}}, M.:
{\em A comprehensive introduction to differential geometry}, 
Vol. II. Second edition. Publish or Perish, Inc., 
Wilmington, Del., 1979, xv+423~pp.

\smallskip

{\bf\bibitem{Spivak-1979-III} 
{\sc Spivak}}, M.:
{\em A comprehensive introduction to differential geometry}, 
Vol. III. Second edition. Publish or Perish, Inc., 
Wilmington, Del., 1979, xii+466~pp.

\smallskip

{\bf\bibitem{Spivak-1979-V} 
{\sc Spivak}}, M.:
{\em A Comprehensive introduction to differential geometry}, 
vol. V. Second edition. Publish or Perish, Inc., 
Wilmington, Del., 1979, viii+661~pp.

}\end{thebibliography}


\vfill\end{document}